\documentclass[12pt,compress]{elsarticle}

\usepackage{ulem}
\usepackage{orcidlink}
\usepackage{amsmath}
\usepackage{xcolor}%
\usepackage{hyperref}
\usepackage{cleveref}
\hypersetup{
    colorlinks = true,
    linkcolor = blue
    }
\usepackage{ragged2e}
\usepackage{algorithm}
\usepackage{algorithmic}
\usepackage{caption}
\usepackage{subcaption}
\crefname{algorithm}{Algorithm}{Algorithms}
\newcommand{\algref}[2]{\hyperref[#1]{Algorithm #2}}
\usepackage{soul}   
\usepackage{adjustbox}
\usepackage{multirow}

\usepackage{array} 
\usepackage{booktabs} 

\usepackage{amssymb}
\usepackage{amsthm}

\newtheorem{theorem}{Theorem}[section]
\newtheorem{assumption}{Assumption}[section]

\newtheorem{lemma}{Lemma}[section] 

\newcommand{\R}{\mathbb{R}}
\newcommand{\Ek}{\mathbb{E}_k}
\newcommand{\sbb}{S_{BB}}
\newcommand{\nosbb}{\overline{S_{BB}}}
\newcommand{\thalf}{\frac{1}{2}}
\newcommand{\Pbk}{\mathbb{P}_k}
\newcommand{\eqdef}{\stackrel{\rm def}{=}}
\newcommand{\ev}{{\cal{E}}_k}
\newcommand{\E}{\mathbb{E}}
\newcommand{\cunosbb}{C_{1,k}^{BB}}
\newcommand{\cunonosbb}{\overline{C_{1,k}^{BB}}}
\newcommand{\cduesbb}{C_{2,k}^{BB}}
\newcommand{\cduenosbb}{\overline{C_{2,k}^{BB}}}

\newcommand{\trishuno}{TRishBB\_v1 }
\newcommand{\trishunop}{TRishBB\_v1}
\newcommand{\trishdue}{TRishBB\_v2 }
\newcommand{\trishduep}{TRishBB\_v2}
\newcommand{\trishtre}{TRishBB\_v3 }
\newcommand{\trishtrep}{TRishBB\_v3}

\journal{Journal of Computational and Applied Mathematics}

\begin{document}

\begin{frontmatter}



\title{Fully stochastic trust-region methods with 
Barzilai-Borwein steplengths}
\author[label1]{Stefania Bellavia\,\orcidlink{0000-0002-3691-7836}}
\ead{stefania.bellavia@unifi.it}
\author[label1]{Benedetta Morini\corref{cor1}\,\orcidlink{0000-0002-9213-3622}}
\ead{benedetta.morini@unifi.it}
\author[label1]{Mahsa Yousefi\,\orcidlink{0000-0002-2937-9654}}
\ead{mahsa.yousefi@unifi.it}
\affiliation[label1]{organization={Department of Industrial Engineering, University of Florence},
addressline={ Viale G.B Morgagni 40, Florence, 50134 (FI)},
            country={ Italy }
}
\cortext[cor1]{Corresponding author}
\begin{abstract}
We investigate stochastic gradient methods and stochastic counterparts of the Barzilai-Borwein steplengths and their application to finite-sum minimization problems. Our proposal is based on the Trust-Region-ish (TRish) framework introduced in   [F. E. Curtis, K. Scheinberg, R. Shi, {\it  A stochastic trust region algorithm based on careful step normalization}, Informs Journal on Optimization, 1, 2019]. The new framework, named TRishBB, aims to enhance the performance of TRish and  at reducing the computational cost of  the second-order TRish variant.  We propose three different methods belonging to the TRishBB framework and present the convergence analysis  for possibly nonconvex objective functions, 
considering biased and unbiased gradient approximations. Our analysis requires  neither diminishing step-sizes nor full gradient evaluation. The numerical experiments in machine learning applications demonstrate the effectiveness of applying the Barzilai-Borwein steplength with stochastic gradients and show improved testing accuracy compared to the  TRish method.
\end{abstract}

\begin{keyword}
finite-sum minimization\sep stochastic trust-region methods\sep Barzilai-Borwein method, machine learning.
\MSC[2010]{65K05 \sep90C30}
\end{keyword}

\end{frontmatter}

\section{Introduction}\label{IntroSec}
We consider the following optimization problem

\begin{equation}\label{eq1}
    \min_{x\in \mathbb{R}^n}\, f(x) = \frac{1}{N}\sum_{i=1}^{N}\, f_i(x),
\end{equation}
which minimizes a finite-sum function $f$ with large $N$. We assume that each function $f_i: \mathbb{R}^n\to\mathbb{R}$, $i\in \mathcal{N}\eqdef\{1,\ldots, N\}$, is continuously differentiable. Several important problems can be stated in this form, e.g., classification problems in Machine Learning (ML) or Deep Learning (DL), data fitting problems, and sample average approximations of an objective function given in the form of mathematical expectation. In the class of first-order optimization methods,  Gradient Descent
(GD) methods generate  iterations of the form
\begin{equation}\label{eq2}
    x_{k+1} = x_k - \mu_k \nabla f(x_k),
\end{equation}
for some positive scalar $\mu_k$ and $k\ge 0$. Convergence and efficiency of the GD method and more generally of gradient-based methods depend heavily 
on the sequence of the steplengths \cite{barzilai1988two, nesterov, nocedal2006numerical}. In \cite{barzilai1988two} Barzilai and Borwein  proposed to use a parameter $\mu_k$  which attempts to capture second-order information without computing the Hessian matrix of $f$ nor its inverse. In other words, the  Barzilai-Borwein (BB) method can be regarded as a quasi-Newton method where the Hessian is approximated by 
$\mu_k^{-1} I$, $I\in \mathbb{R}^{n\times n}$ being the identity matrix, and   $\mu_k$ minimizes the least-squares error of a secant equation involving the difference between two successive iterates and the difference between two successive gradients.

Computing the gradient $\nabla f(x_k)$ to perform the update \eqref{eq2} is practically inconvenient in ML and DL applications where $f$ represents a loss function on a given dataset and $N$ is large. The Stochastic Gradient (SG) method \cite{robbinsmonro, bottou2004large} and the mini-batch SG method \cite{bottou2018optimization} are iterative procedures where the full gradient $\nabla f(x_k)$ is replaced   by either a single gradient $\nabla f_i(x_k)$, $i\in \mathcal{N}$ selected at random, or by a mini-batch of samples of the form
$\nabla f_{\mathcal{N}_k}(x_{k})=\frac{1}{|\mathcal{N}_k|}\sum_{i\in \mathcal{N}_k}\nabla f_i(x_k)$ where the set $\mathcal{N}_k \subseteq \mathcal{N}$ is selected at random and $|\mathcal{N}_k|$ denotes its  cardinality, i.e., the mini-batch size. A limitation of such methods is the need to choose problem-dependent steplengths to guarantee  convergence. An attempt to reduce such dependence was made in \cite{curtis2019stochastic} where the steplengths depend on  a Trust-Region-ish (TRish) rule  and normalized steps  are selected  whenever the norm of the stochastic gradient is within a prefixed interval whose endpoints need to be tuned. The TRish method was designed for both first-order models \cite{curtis2019stochastic,  bellavia2024trish} and quadratic models \cite{curtis2022fully} and, 
contrary to a standard trust-region approach, does not test for actual  vs.  predicted reduction and thus does not call for $f$-evaluations.
Results from \cite{curtis2019stochastic, bellavia2024trish} show that TRish can outperform a standard stochastic gradient approach and it was recently extended to  equality-constrained problems \cite{fang2024fully}.

The aim of this work is to combine the TRish methodology and stochastic BB steplengths in order to employ second-order information at a low computational cost; the selection of BB steplengths is ruled by TRish automatically.  We devise the computation of BB steplengths in a stochastic setting and provide a theoretical analysis without requiring either diminishing steplengths nor unbiased stochastic gradients. In addition, we design three different procedures belonging to the above framework and show numerical results where our procedures enhance the performance of the first-order TRish method \cite{curtis2022fully}.

This paper is organized as follows. The remainder of this section is devoted to a review of the relevant literature providing background, motivation, and contribution to our work. We outline the proposed stochastic algorithms
in \Cref{MethodSec} and present the theoretical analysis in \Cref{ConvergenceSec}. In \Cref{NumericalSec} we present the implementation details, experimental setup, and an extensive evaluation of the method's performance. Finally, in \Cref{ConSec} we provide an overview of the obtained results and perspectives.

\paragraph{\textbf{Notations}}
The Euclidean norm is denoted as $\|\cdot\|$.
We use the symbol $|\cdot|$ to denote the absolute value of a scalar, the cardinality of a set and the number of columns of a matrix, the meaning will be clear from the context. The operation $\lfloor\cdot \rfloor$ gives the nearest integer less than or equal to the operand.
The modulo operation $\text{\textbf{mod}}(p,t)$ returns the remainder after division of $p$ by $t$.  
$\Ek[\cdot]$
and $\Pbk[\cdot]$ denote, respectively, the conditional expectation and probability given that the algorithm has reached $x_k$ at the $k$-th iteration. $\E[\cdot]$ denotes the expected value taken with the joint distribution of all random variables.

\subsection{Related works and motivations} 
The use of  BB steplengths combined with  stochastic gradients calls for the definition of the stochastic counterparts of  the standard BB steplengths and for convergence analysis in expectation or high probability of the resulting procedures.
Several attempts have been made \cite{tan2016barzilai,wang2021stochastic,liang2019barzilai,bellavia2023slises,kkj2018} 
and typically  originate from the class of stochastic quasi-Newton methods \cite{berahas2020robust, bollapragada2018progressive, gower2016stochastic, mokhtari2015global, wang2017stochastic, erway2020trust, yousefi2023deep}.

In \cite{tan2016barzilai}, BB steplengths were incorporated   into the Stochastic Variance Reduction Gradient (SVRG) method \cite{johnson2013accelerating} and into the SG method. The proposed methods are inner-outer loop-based frameworks. Since SVRG provides a full gradient at each outer loop, the BB steplength was computed at outer loops using the standard formula for the BB steplength; convergence of  the resulting  procedure was proved for   strongly convex objective functions. On the other hand, single stochastic  gradients $\nabla f_i$ were employed at each iteration in the SG method and, in the inner loop, a  moving averaging rule over $N$ stochastic gradients at different points was computed; then, such average was used at each outer loop to evaluate a stochastic BB parameter. A  theoretical support for this algorithm, denoted SGD-BB algorithm, was not provided. In \cite{wang2021stochastic}, deterministic BB steplengths from  \cite{dai2002modified}, based on exact gradient and function evaluations,  were used in the SVRG method; convergence was  analyzed  for strongly convex objective functions. 
In \cite{liang2019barzilai}, stochastic BB steplengths  were formed by a moving average rule
on mini-batch gradients evaluated at different points and selected if they lied within an interval with endpoints tending to zero. Using diminishing stepsizes provides convergence due to well-known results for SG method but, in principle, BB steplengths may be discarded eventually  since they are not supposed to tend to zero.
A further method in which stochastic  BB steplengths were formed and diminishing stepsizes were imposed is 
given  in \cite{bellavia2023slises}; the distinguishing feature is that the mini-batch of sample is taken constant for a prefixed number of iterations with the aim to explore the underlying approximate Hessian spectrum.  Finally, the paper \cite{kkj2018} investigates on  the use of BB steplengths combined with projected stochastic  gradients for constrained optimization problems and with the objective function in the form of mathematical expectation.

The above mentioned papers show limitations such as use of diminishing stepsizes, use of increasing batch sizes $|{\cal{N}}_k|$, use of full gradient information, need for $f$-evaluations. Most of the theoretical results hold for convex or strongly convex objective functions. Our proposal overcomes the above drawbacks employing 
the TRish paradigm given in 
\cite{curtis2019stochastic, curtis2022fully} which allows to use a simple stochastic quadratic model inspired by the BB method and to select  adaptively the step by a trust-region methodology not invoking $f$-evaluations. The algorithms proposed in this work are cheap to implement and supported by theoretical results for possibly nonconvex objective functions and for biased or unbiased gradients.

\section{New Stochastic Trust-Region-ish Algorithms}\label{MethodSec}
\begin{algorithm}[t]
\caption{{\bf TRishBB}}\label{alg.TRishBB}
\smallskip

\begin{algorithmic}[1]
\justifying
\STATE{Choose an initial iterate $x_0\in\mathbb{R}^n$  and a parameter $\mu_0>0$.}
\STATE{Choose   $(\alpha,  \gamma_{1}, \gamma_{2})$ such that  $\alpha>0$, 
 $0<\gamma_{2}\leq\gamma_{1}<\infty$ and $\mu_{min}>0$.}
\FOR{$k=0,1,\ldots$}
\smallskip

\STATE{Set $H_k= \mu_k^{-1} I$.}
\STATE{Compute a stochastic gradient $g_k$.}
\STATE{Compute the solution $p_k$ of the trust-region subproblem} 
\begin{equation}\label{eq.tr_pb}
\min_{p\in \mathbb{R}^n} m_k(p)= g_k^{T}p + \frac{1}{2} p^TH_kp \quad \text{s.t.}\quad \|p\|\leq \Delta_k,
\end{equation}
where
\begin{equation}\label{eq.radius_tr}
\Delta_{k}=
\begin{cases}
\alpha\gamma_{1}\|g_k\|
& \quad \text{if } \|g_k\| \in \left [0, \frac{1}{\gamma_{1}}\right)\\
\alpha  
& \quad \text{if } \|g_k\| \in \left[\frac{1}{\gamma_{1}}, \frac{1}{\gamma_{2}}\right]\\
\alpha\gamma_{2}\|g_k\|  
& \quad \text{if } \|g_k\| \in \left(\frac{1}{\gamma_{2}}, \infty\right)
\end{cases}.
\end{equation}
\STATE{Set $x_{k+1}=x_k+p_k$.}
\STATE{Form a parameter $\mu_{k+1}\ge\mu_{min}$.}
\ENDFOR
\end{algorithmic}
\end{algorithm}
In this section we present our trust-region-ish algorithm employing stochastic BB steplengths.
First, we provide  a brief overview of the BB gradient method in the deterministic setting \cite{barzilai1988two, Raydan93, BirginMartRaydan93,di2018steplength},  second  we 
describe the main features of the method denoted as TRishBB, finally we present three algorithms belonging to the  TRishBB frameworks.

At $k$-th iteration, the standard BB method minimizes the second-order model
$m_k(p)= \nabla f(x_k)^{T}p + \frac{1}{2} p^TH_kp$ where $H_k= \mu_k^{-1} I$ for some positive $\mu_k$; 
thus the iterates have the form
\eqref{eq2}.  Regarding the computation of the steplength, given $x_{k+1}$, the main idea of the BB method is to form the parameter   $\mu_{k+1}^{\text{BB}}$ which   solves
\begin{equation}\label{stepBB}
\mu_{k+1}^{\text{BB}}=\text{argmin}_{\mu} \  \|\mu^{-1} s_{k}-y_{k}\|^2 = \frac{s_{k}^{T}s_{k}}{s_k^{T}y_k},
\end{equation}
with 
$s_{k}=x_{k+1}-x_{k}$, and $y_{k}= \nabla f(x_{k+1}) - \nabla f(x_{k})$ commonly denoted as correction pair. If $f$ is strictly convex, then the scalar $\mu_{k+1}^{\text{BB}}$ is positive.

If $f$ is twice continuously differentiable then  $\mu_{k+1}^{\text{BB}}$
is related to the spectrum of the Hessian matrix of $f$ as follows.  Given a symmetric matrix $A\in \mathbb{R}^{n\times n} $, 
let $q(A,v)$ be the Rayleigh quotient associated to the vector $v\in \mathbb{R}^n$;  in addition, let $\lambda_{\min}(A)$ and $\lambda_{\max}(A)$ denote the smallest and largest eigenvalue of $A$, respectively. Then, by the mean value theorem \cite[Th. 3.1.3]{conn2000trust}, we have
\begin{eqnarray}
\mu_{k+1}^{\text{BB}} &=&\frac{s_{k}^{T}s_{k}}{s_{k}^{T}\left(\int_0^1 \nabla^2 f(x_{k}+ts_{k}) dt\right )s_{k}}=\frac{1}{q\left( \int_0^1 \nabla^2 f(x_{k}+ts_{k}) dt, s_{k}\right)} \label{eq.Rayleigh}
\end{eqnarray} 
Thus,
$\mu_{k+1}^{\text{BB}}\in  \left[\frac{1}{\lambda_{\max}\left(\int_0^1 \nabla^2 f(x_{k}+ts_{k}) dt\right)},
\frac{1}{ \lambda_{\min}\left(\int_0^1 \nabla^2 f(x_{k}+ts_{k}) dt\right)}\right]$  if $f$ is strictly convex. On the other hand, for a general function $f$, $\mu_{k+1}^{\text{BB}}$  is related to the spectrum of the matrix $\int_0^1 \nabla^2 f(x_{k}+ts_{k}) dt$ by means of \eqref{eq.Rayleigh} but  it may be negative. 

In the BB method,  $H_{k+1}=\mu_{k+1}^{-1} I$ is formed using 
$|\mu_{k+1}^{\text{BB}}|$ 
and thresholding $|\mu_{k+1}^{\text{BB}}|$ by means of  two predefined scalars $0<\mu_{min}<\mu_{max}$, i.e., setting
$\mu_{k+1}=\max\{\mu_{min}, \min\{|\mu_{k+1}^{\text{BB}}|,\mu_{max}\} $  to prevent excessively small or large values; in the case of convex problems it trivially holds $\mu_{k+1}=\max\{\mu_{min}, \min\{\mu_{k+1}^{\text{BB}},\mu_{max}\} $.
We refer to 
\cite{di2018steplength} for a survey on both practical rules for fixing $\mu_{k+1}$ and  step acceptance rules in the deterministic  convex and nonconvex settings.

Now we analyze \Cref{alg.TRishBB} which is based on the second-order version of TRish in \cite{curtis2022fully} and specialized to the case where the Hessian matrix is approximated by the diagonal matrix   $H_k = \mu_k^{-1} I$, $\forall k$. Our goal is to use steplengths $\mu_k$'s inspired by the BB method and  to enhance 
the first-order version of TRish \cite{curtis2019stochastic} where $H_k=0$ for all $k\ge 0$. However, since stochastic gradients are available, the correction pair $(s_k, y_k)$ is not uniquely defined and $\mu_{k+1}^{\text{BB}}$  in (\ref{stepBB})  is a random parameter.
The choice of the  pair $(s_k, y_k)$ is crucial
for the practical behaviour of the method and we postpone the discussion of three options for the 
definition of   $(s_k, y_k)$ and  $\mu_{k+1}$   to Subsections~\ref{Sub2.2}-\ref{Sub2.4}.

At the generic  $k$-th iteration of TRishBB, the stochastic gradient $g_k$ is computed. Then, given the trust-region radius $\Delta_k$ and  the scalar $\mu_k>\mu_{min}$, we set  $H_k = \mu_k^{-1} I$ and let \eqref{eq.tr_pb}
be the trust-region problem with strictly convex quadratic model. The trust-region solution is cheap to compute and takes the form
\begin{equation}\label{eq.Solpk}
    p_k=\begin{cases}
 -\mu_kg_k  & \quad \text{if } \|\mu_kg_k\| < \Delta_k, \\
\displaystyle -\frac{\Delta_k}{\|g_k\|}g_k  & \quad \text{otherwise. }
\end{cases}
\end{equation}
In fact, when  $p_k$ takes the form $p_k = -H_k^{-1}g_k=-\mu_k g_k$,
it holds  $\|p_k\|< \Delta_k$ and $p_k$  is the unconstrained minimizer of the model $m_k$. We refer to such a step as  
 the unconstrained solution of \eqref{eq.tr_pb}. Otherwise, the solution lies on the boundary of the trust-region and  it is the constrained Cauchy point; we will refer to this step as the constrained solution. In both cases,  $p_k$ is a step along $(-g_k)$ and the cost for solving the trust-region problem is low. On the contrary, using a second-order trust-region methods and a general $H_k$ requires a linear algebra phase such as the application of the conjugate gradient method \cite[\S 7.5.1]{conn2000trust}.

The updating rule \eqref{eq.radius_tr} for $\Delta_k$ is a distinctive feature of TRish. It is ruled
by three positive and predefined scalars, $\alpha$, $\gamma_1$ and $\gamma_2$ and the magnitude of the norm of $g_k$.
The insight of this rule is given in  \cite[\S 2]{curtis2019stochastic}; briefly, adopting the standard rule, i.e., setting 
$\Delta_k = \alpha$, and thus normalizing $(-g_k)$ if the trust-region solution is constrained,  may cause the algorithm to compute ascent directions more likely  than descent directions  and then fail to progress in expectation. We note that in \cite{curtis2019stochastic, curtis2022fully} the algorithm admits the use of  positive sequences  $\{\alpha_k\}$, $\{\gamma_{1,k}\}$ and $\{\gamma_{2,k}\}$ but here we focus on the choice of constant parameters.

Due to  \eqref{eq.radius_tr} and \eqref{eq.Solpk}, the constrained trust-region solution has the form 
\begin{equation}\label{eq.step_TRish_first_order}
p_k=
  \begin{cases}
    -\alpha \gamma_{1} g_k     
    & \quad \text{if } \|g_k\| \in \left[0, \frac{1}{\gamma_{1}}\right) \\
    - \displaystyle \alpha \frac{g_k}{\|g_k\|} & \quad \text{if } \|g_k\| \in \left[\frac{1}{\gamma_{1}}, \frac{1}{\gamma_{2}}\right]\\
    -\alpha \gamma_{2} g_k  & \quad \text{if } \|g_k\| \in \left( \frac{1}{\gamma_{2}}, \infty\right)
  \end{cases},
\end{equation}
analogously to the first-order TRish method in \cite{curtis2019stochastic}, whereas 
the unconstrained trust-region solution is $p_k=-\mu_k g_k$ when
\begin{equation}\label{eq.step_TRishBB}
\mu_{k}<
  \begin{cases}
  \alpha \gamma_{1}  & \quad \text{if } \|g_k\| \in \left [0, \frac{1}{\gamma_{1}}\right)\\
      \displaystyle \frac{\alpha}{\|g_k\|} & \quad \text{if } \|g_k\| \in \left[\frac{1}{\gamma_{1}}, \frac{1}{\gamma_{2}}\right]\\
     \alpha \gamma_{2}   & \quad \text{if } \|g_k\| \in \left( \frac{1}{\gamma_{2}}, \infty\right)
  \end{cases}.
\end{equation}

Importantly, once the step $p_k$ is formed, the new iterate $x_{k+1}$ is computed without performing the standard test 
of trust-region methods on predicted vs. actual reduction  \cite{conn2000trust}; thus the Algorithm TRishBB does not require the evaluation of the objective function 
$f$.

We conclude  this   section  noting that each step $p_k$ taken by TRishBB is along $(-g_k)$ as in a first-order method
but the unconstrained trust-region solution aims to exploit some second-order information
via $\mu_k$. 

In the following subsections we explore three ways to form stochastic  BB steplengths;
each of them provides a variant  of the general Algorithm TRishBB. 

\subsection{\trishuno} \label{Sub2.2}
\begin{algorithm}[t]
\caption{{\bf\trishuno (Iteration $k$) }}\label{alg.TRishBB(1)}
\smallskip

\begin{algorithmic}[1]
\justifying
\STATE{Given $x_{k}\in\mathbb{R}^n$, $(\alpha, \gamma_1, \gamma_2)$, $\mu_k>0$,
$0<\mu_{min}<\mu_{max}$,  $m\geq 1$.}
\smallskip

\STATE{Choose $\mathcal{N}_k\subseteq \mathcal{N}$, compute $g_k=\nabla f_{\mathcal{N}_k}(x_{k})$ 
and set $H_k= \mu_k^{-1} I$.}
\STATE{Obtain $p_k$ from \eqref{eq.Solpk} and set $x_{k+1}=x_k+p_k$.}
\IF{$\text{\textbf{mod}}(k,m) = 0$}
\STATE{$s_k = p_k$ and $y_k = \nabla f_{\mathcal{N}_k}(x_{k+1}) - g_k$}
\STATE{$\mu_{k+1} = \max\left\{\mu_{min}, \min\left\{\left|\frac{s_k^{T}s_k}{s_k^{T}y_k}\right|, \mu_{max}\right\} \right \}$}
\ELSE
\STATE{$\mu_{k+1}=\mu_k$}
\ENDIF
\end{algorithmic}
\end{algorithm}
The displayed Algorithm \trishuno is a straightforward implementation of the BB method employing stochastic gradients and  guidelines from stochastic quasi-Newton methods.
Specifically, once the trust-region solution $p_k$ in \eqref{eq.Solpk} is computed, the successive steplength
$\mu_{k+1}$ is formed using the correction pair $(s_k, y_k)$ at line 5 of this algorithm.
The vector $s_k$ coincides with  $p_k$, i.e., it is the step between $x_{k+1}$ and $x_k$ whereas
the vector $y_k$ is the difference between the mini-batch gradients
$\nabla f_{\mathcal{N}_k}(x_{k+1})$ and $\nabla f_{\mathcal{N}_k}(x_{k})$. We underline that using the same 
set of indexes  $\mathcal{N}_k$ in the computation of $y_k$ is suggested in \cite{schraudolph2007stochastic, berahas2020robust,
bollapragada2018progressive} and implies that \eqref{eq.Rayleigh} holds with the matrix $\int_0^1 \nabla^2 f_{\mathcal{N}_k}(x_k+ts_k)$. Thus,
the BB steplength inherits spectral information as long as it is positive and  $f$ is twice continuously differentiable.
Following a heuristic common rule, see e.g.,  \cite{cyclicBB}, given an integer $m\geq 1$, the parameter $\mu_{k+1}$ is updated every $m$ iterations.

Compared to TRish, the extra cost of this algorithm amounts to the computation of $\nabla f_{\mathcal{N}_k}(x_{k+1})$ and two scalar products to evaluate $\mu_{k+1}$.

\begin{algorithm}[t]
\caption{{\bf \trishdue} (Iteration $k$)}
\label{alg.TRishBB(2)}
\smallskip

\begin{algorithmic}[1]
\justifying
\STATE{Given $x_{k}\in\mathbb{R}^n$, $(\alpha, \gamma_1, \gamma_2)$, $\mu_k >0$,
$\bar \mu>0$, $0<\mu_{min}< \mu_{max}$, $m\ge 1$,  $\beta=\frac{m- 1}{m}$, $\eta \in (0,1)$, $\bar{x}_{old}$, $\bar{g}_{old}$, $\bar{g}$\, (if $k=0$, then $\bar{x}_{old} = x_0$, $\bar{g}_{old}=0$, $\bar{g} =  0$, and $\bar{\mu} = {\mu}_0$):}
\smallskip

\STATE{Choose $\mathcal{N}_k\subseteq \mathcal{N}$, compute 
$g_k = \nabla f_{\mathcal{N}_k}(x_k)$ and set $H_k = \mu_k^{-1} I$;}
\STATE{Obtain $p_k$ from \eqref{eq.Solpk} and set $x_{k+1} = x_k + p_k$;}
\STATE{Compute $\bar g =  \beta \bar g + (1-\beta) g_k$;}
\IF{$k> 0$\, and\, $\text{\textbf{mod}}(k,m)=0$}
\STATE{$\bar{x} = x_{k+1}$}
\STATE{$s_{k} = \bar{x}  - \bar{x}_{old}$ and $y_{k}= \bar{g} - \bar{g}_{old}$}
\STATE{ $\widehat{\mu}_k =\frac{1}{m}\left|\frac{s_k^{T}s_k}{s_k^{T}y_k}\right|$, and $\bar{\mu} = \eta\bar{\mu} + (1-\eta)\widehat{\mu}_k$}
\STATE{ $\mu_{k+1} = \max\{\mu_{min}, \min\{\bar{\mu}, \mu_{max}\}\}$}
\STATE{$\bar{x}_{old} =  \bar{x}$, $\bar{g}_{old} =  \bar{g}$, $\bar{g} =  0$}
\ELSE
\STATE{$\mu_{k+1}=\mu_k$}
\ENDIF
\end{algorithmic}
\end{algorithm}
\subsection{\trishdue}\label{Sub2.3}
In Algorithm \trishuno  the pair $(s_k, y_k)$ contains information  
from one single iteration and $y_k$ is  obtained using stochastic gradients, possibly evaluated on 
a small sample ${\mathcal{N}_k}$.
On the contrary, Algorithm \trishdue forms a correction pair that accumulates information over a prefixed number $m$ of iterations and is
described below.

The vector $s_k$ is the difference between two, not consecutive, iterates;
the indexes of such iterates differ by $m$, see lines 7 and 10 of the algorithm.  
The vector  $y_k$ is the difference between vectors computed via a moving average rule which accumulates stochastic gradients over $m$ iterations, see lines 4, 7 and 10 of the algorithm; in particular $\bar g$ is formed 
 over the last $m$ iterations while $\bar g_{old}$ is 
the  analogous vector formed in the previous cycle of $m$ iterations.  The parameter $\beta=(m-1)/m$ is fixed in accordance with \cite{kingma2017adam} to avoid hyper-parameter tuning.

The form of $s_k$ and $y_k$ and $\widehat \mu_k$ is borrowed from the stochastic gradient descent algorithm with BB step size (SGD-BB) algorithm proposed in \cite{tan2016barzilai} where $|\mathcal{N}_k|=1$, $\forall k$, and  $m=N$, i.e., $N$ single gradient evaluations are performed before updating the steplength. In \trishdue we allow $|\mathcal{N}_k|>1$; letting $N_b= \lfloor N/|\mathcal{N}_k|\rfloor$ be (roughly) the number of disjoint mini-batches of size  $|\mathcal{N}_k|$ in $\mathcal{N}$ and $m=N_b$, the cost for computing $\bar g$ through $m$ iterations is approximately the cost for one full gradient evaluation. 

The steplength $\mu_{k+1}$ is updated every $m$ iterations and computed using accumulated information, see lines 8 and 9 of the algorithm. Given  $\eta\in (0,1)$, the scalar $\bar \mu$ is computed as the  convex combination of the current stochastic BB steplength  $\widehat \mu_k$ and the current accumulated parameter $\bar \mu$.  The larger 
$\eta$ is, the higher   the penalization to $\widehat \mu_k$ is and this may soften the effect of intrinsic noise arising from stochastic gradients.

Compared to TRish, the extra cost of this algorithm amounts to two scalar products to evaluate $\mu_{k+1}$; in terms of storage,  three vectors of size $n$ are required.

\begin{algorithm}[t]
\caption{{\bf \trishtre} (Iteration $k$)}\label{alg.TRishBB(3)}
\smallskip

\begin{algorithmic}[1]
\justifying
\small{
\STATE{Given $x_{k}\in\mathbb{R}^n$, $(\alpha, \gamma_1, \gamma_2)$, $\mu_k >0$, $\bar \mu>0$, $0<\mu_{min}< \mu_{max}$, $m\ge 1$, 
$m_F\in \mathbb{N}$
$F_k\in \mathbb{R}^{n\times q}$ with $q\le m_F$, $x_{avg}$, and $\bar{x}_{old}$, \, (if $k=0$, then $F_0=[\, ]$, $x_{avg} = 0$, and $\bar{x}_{old} =0 $, and $\bar{\mu} = {\mu}_0$).} \smallskip

\STATE{$x_{avg} = x_{avg} + x_k$}
\STATE{Choose $\mathcal{N}_k\subseteq \mathcal{N}$, compute  $g_k = \nabla f_{\mathcal{N}_k}(x_k)$ 
and set $H_k = \mu_k^{-1} I$. }
\STATE{Obtain $p_k$ from \eqref{eq.Solpk} and set $x_{k+1} = x_k + p_k$.}
\IF{$k\leq m_F$}
\STATE{Store $g_k$ in $k$-th column of $F_k$.}
\ELSE
\STATE{Remove the first column of $F_k$, and store $g_k$ in its $m_F$-th column.}
\ENDIF
\IF{$\text{\textbf{mod}}(k,m) = 0$}
\STATE{$\bar{x} = \frac{x_{avg}}{m}$, and $x_{avg} = 0$}
\IF{$k>0$}
\STATE{${s}_k = \bar{x} - \bar{x}_{old}$, and ${y}_k = \frac{1}{|F_k|} F_k(F_k^T s_k)$}
\STATE{ $\widehat{\mu}_k =\frac{1}{m}\left|\frac{s_k^{T}s_k}{s_k^{T}y_k}\right|$, and $\bar{\mu} = \eta\bar{\mu} + (1-\eta)\widehat{\mu}_k$;}
\STATE{ $\mu_{k+1} = \max\{\mu_{min}, \min\{\bar{\mu}, \mu_{max}\}\}$;}
\ENDIF
\STATE{$\bar{x}_{old} = \bar{x}$}
\ELSE
\STATE{$\mu_{k+1}=\mu_k$}
\ENDIF
}
\end{algorithmic}
\end{algorithm}
\subsection{\trishtre}\label{Sub2.4} Algorithm \trishtre is originated by the arguments in \cite{byrd2016stochastic} to avoid both the potentially harmful effects of the gradient difference when $\|s_k\|$ is small and the noise in stochastic gradient estimates. In \cite{byrd2016stochastic}, given $s_k=\bar x-\bar x_{old}$ for some $\bar x, \, \bar x_{old}$, the vector $y_k$ is defined as  $y_k=\nabla^2 f_{\widetilde{\mathcal{N}_k}}(\bar{x})s_k$ where  $\nabla^2 f_{\widetilde{\mathcal{ N}_k}}(\bar{x}) $ is a subsampled Hessian approximation and $|\widetilde{\mathcal{N}_k}|$ is supposed to be larger than  the batch size used for gradients. If the objective function $f$ is the cross-entropy loss function or the least-squares function, then
$\nabla^2 f_{\widetilde{\mathcal{ N}_k}}(\bar{x}) $ can be approximated by the empirical Fisher Information Matrix (eFIM) 
\cite[\S 11]{martens2020new} defined as
$$\frac{1}{|\widetilde{\mathcal{N}_k}|}\sum_{i\in \widetilde{\mathcal{N}_k}} \nabla f_i(\bar{x}) \nabla f_i(\bar{x})^T.$$

In order to reduce the cost, the eFIM was further simplified by introducing the 
\textit{accumulated} Fisher Information Matrix (aFIM) \cite{keskar2016adaqn, faria2023fisher}.
In Algorithm \trishtrep, we update $s_k$ and $y_k$ every $m$ iterations, and $y_k$ is computed using the aFIM and the first-in, first-out (FIFO) approach \cite{keskar2016adaqn, faria2023fisher}. In particular, as displayed in line 13 of the algorithm, $s_k$ is formed using the averages $\bar x$ and $\bar x_{old}$; $\bar x$ is the average of the last $m$ iterates, and $\bar x_{old}$ is the average of the previous $m$ iterates, see lines 2 and 11 of the algorithm.
Average iterates are meant to represent more stable approximations than a  single iterate for computing the difference $s_k$. Concerning $y_k$, let $F_k$ be   the limited memory matrix whose columns contain at most the last $m_F$ computed stochastic gradients for some positive integer $m_F$. Then, the curvature vector $y_k$
is defined as ${y}_k = \frac{1}{|F_k|} F_k(F_k^T s_k)$, see  line 13 of the algorithm, where $|F_k|$ denotes the number of columns of $F_k$.  Since $F_kF_k^T$ is the sum of $|F_k|$ rank-one matrices,  \eqref{eq.Rayleigh}  implies
$$\widehat \mu_{k}= \frac{|F_k|}{m} \frac{1}{\sum_{t=k-|F_k|+1}^{k} q(  \nabla f_{\mathcal{N}_t}(x_t) ( \nabla f_{\mathcal{N}_t}(x_t))^T,  {s}_k)},$$ 
where $\nabla f_{\mathcal{N}_t}(x_t) ( \nabla f_{\mathcal{N}_t}(x_t))^T$ is a rank-one matrix.
Finally, the update of $\mu_k$ is performed every $m$ iterations and it is the same as in \trishduep.

Compared to TRish, the extra cost of this algorithm amounts to $m_F+2$ scalar products to evaluate $\mu_{k+1}$; in terms of storage, two  vectors of size $n$ and one  matrix of size $n \times m_F$ are required.

\section{Convergence Analysis}\label{ConvergenceSec}
We present the convergence analysis of TRishBB considering the use of both biased and unbiased stochastic gradients. Further, we   
investigate the case where the function $f$ satisfies the Polyak-Lojasiewicz (PL) condition  and the case where  $f$ is a general nonconvex function. 

The assumption of the objective function  made in our analysis and the PL condition are presented below.
\begin{assumption}
\label{ass.f}
 Function  $f:\mathbb{R}^{n} \to \mathbb{R}$ is continuously differentiable, bounded below by $f_*$, and  there exists a positive $L$ such that
  \begin{equation}\label{eq.g_Lip}
    f(x) \leq f(y) + \nabla f(y)^T(x - y) +  \frac{L}{2} \|x - y\|^2\ \ \text{for all}\ \  x, y  \in \mathbb{R}^{n}.
  \end{equation}
\end{assumption}

\begin{assumption}\label{ass.PL}
 For any $x \in \R^{n}$, the Polyak-Lojasiewicz (PL) condition holds with $c \in (0,\infty)$, i.e.,
\begin{equation}\label{eq.PL}
  2c(f(x) - f_*) \leq \|\nabla f(x)\|^2\ \ \text{for all}\ \ x \in \R^{n}.
\end{equation}
\end{assumption}
In our analysis we characterize the iterations of TRishBB by three cases:  ``Case 1''$: \|g_k\| \in \left[0,\tfrac{1}{\gamma_{1}}\right)$, ``Case 2'': $\|g_k\| \in \left[\tfrac{1}{\gamma_{1}},\tfrac{1}{\gamma_{2}}\right]$, and ``Case 3'': $\|g_k\| \in \left(\tfrac{1}{\gamma_{2}},\infty\right)$.   We denote by $C_{i,k}$ the events where Case $i$ occurs, for $i \in \{1,2,3\}$.

\subsection{Analysis with Biased Gradients}\label{sec_biased}
The results presented in this subsection concern biased gradients satisfying the following assumption.
\begin{assumption}
\label{hp_g}
There exist positive $\omega, \, M_1, M_2$ such that the stochastic gradient $g_k$ satisfy
\begin{eqnarray}
& & \nabla f(x_k)^T\Ek[g_k]\ge \omega\|\nabla f(x_k)\|^2, \label{g_biased} \\
& &  \Ek[\|g_k\|^2]\le M_1+M_2 \|\nabla f(x_k)\|^2. \label{second_momentum}
\end{eqnarray}
\end{assumption}
\noindent
The inequality \eqref{second_momentum} implies that
\begin{equation}\label{eq.ce_fg_upperbound}
\Pbk[\ev]\E_k[\nabla f(x_k)^T g_k | \ev] \leq h_1 + h_2\|\nabla f(x_k)\|^2,
\end{equation}
with  $h_1 =\thalf\sqrt{M_1}$ and $h_2= h_1 + \sqrt{M_2}$; see  \cite[Lemma 2]{curtis2019stochastic}.

We also make the following assumption of the parameter $\mu_{min}$.
\begin{assumption}
\label{hp_mu}
Given positive $\mu_{min}$, for all $k\ge 0$, it holds
$$
\mu_{min} \le \gamma_2\alpha.
$$
\end{assumption}
The results presented are based on the analysis developed by \cite{curtis2019stochastic}. The following two lemma are technical results. As for the notation, for the events of taking the unconstrained and the constrained solutions \eqref{eq.Solpk}  we use the notations $\sbb$ and $\nosbb$, respectively. Moreover, for the events $C_{i,k}$,  $i \in \{1,2,3\}$, that occur within each defined case  we set ${C_{i,k}^{BB}}\eqdef C_{i,k}\cap \sbb$ and $\overline{C_{i,k}^{BB}}\eqdef C_{i,k}\cap \nosbb$. Finally, we denote with 
 $\ev$ the event that $\nabla f(x_k)^Tg_k \geq 0$ and with $\overline{\ev}$  the event that $\nabla f(x_k)^Tg_k < 0$.

\begin{lemma}\label{lem.basic} 
Under Assumptions~\ref{ass.f} and \ref{hp_mu}, the iterates  generated by Algorithm TRishBB satisfy
\begin{eqnarray}
\left. \begin{aligned}
 \Ek[f(x_{k+1})] - f(x_k) \leq & (\gamma_{1} \alpha -\mu_{min}) \Pbk[\ev] \Ek[\nabla f(x_k)^Tg_k | \ev]\\
+ & \frac{L}{2} \gamma_{1}^2 \alpha^2 \Ek[\|g_k\|^2] - \gamma_{1} \alpha\E_k[\nabla f(x_k)^Tg_k],
\end{aligned}
\right.
\end{eqnarray}\label{eq.cases} 
for all $k$, where $\ev$ is the event that $\nabla f(x_k)^Tg_k \geq 0$.
\end{lemma}

\begin{proof}
See \Cref{AppB1}.
\end{proof}
\begin{lemma}\label{lem.basic2}
Under Assumptions~\ref{ass.f},\ref{hp_g}, 
and \ref{hp_mu}, the iterates generated  by Algorithm TRishBB satisfy
\begin{align}\label{eq.cases.final}
& \Ek[f(x_{k+1})] - f(x_k) \leq  (\gamma_{1}\alpha- \mu_{min})h_1   +\frac{L}{2}\gamma_{1}^2 \alpha^2  M_1\\
& - \gamma_{1}\alpha \left(\omega -
\frac{L}{2} \gamma_{1} M_2 \alpha\right)\|\nabla f(x_k)\|^2 + 
(\gamma_{1}\alpha- \mu_{min})h_2\|\nabla f(x_k)\|^2, \nonumber
\end{align}
for all $k$ with $h_1, h_2$ given in  \eqref{eq.ce_fg_upperbound}.
\end{lemma}

\begin{proof}
It is a direct consequence of Assumption \ref{hp_g}, \Cref{lem.basic} and \eqref{eq.ce_fg_upperbound}.
\end{proof}

For the rest of analysis, we define 
\begin{equation}\label{theta}
\theta_1=\thalf \omega \gamma_1-{h_2}(\gamma_1 - \rho \gamma_2), \quad \theta_2=h_1(\gamma_{1}\alpha\ - \mu_{min})
+\thalf \gamma_{1}^2 LM_1 \alpha^2,
\end{equation}
where $\rho=\frac{\mu_{min}}{\alpha \gamma_2}\in (0,1]$, $h_1$ and $h_2$ are given in \eqref{eq.ce_fg_upperbound}, $\omega$ and $M_1$ are given in Assumption \ref{hp_g}. 
Note that $\theta_1$ is positive for any choice of $(\gamma_1,\gamma_2)$  in the case 
$h_2-\frac{1}{2}\omega\le 0$, while 
condition
\begin{equation}\label{gamma1gamma2}
   \frac{\gamma_1}{\gamma_2}< \frac{\rho h_2}{h_2-\frac{1}{2}\omega} 
\end{equation}
ensures $\theta_1> 0$ in the case $h_2-\frac{1}{2}\omega>0$. Further, note that $\theta_2$ is positive by Assumption $\ref{hp_mu}$.

First we consider the case where the objective function satisfy the Polyak-Lojasiewicz condition \eqref{ass.PL}, then we  analyze the general case.

\begin{theorem}\label{th.bounded_variance_fixed_stepsize}
Under  Assumptions~\ref{ass.f}--\ref{hp_mu}, suppose that  \eqref{gamma1gamma2} holds
in case $h_2-\frac{1}{2}\omega>0$ and that 
\begin{equation}\label{choice_alpha_PL}
\alpha\le \min\left\{\frac{\omega}{\gamma_1LM_2}, \frac{1}{2c\theta_1}\right\},
\end{equation}
with $h_2$ given in (\ref{eq.ce_fg_upperbound}) and  $\theta_1$ given in (\ref{theta}).
Then, the expected optimality gap of the iterates generated by Algorithm TRishBB satisfies
\begin{equation}\label{eq.fixed_stepsize_PL}
\E[f(x_{k+1})] - f_*  \xrightarrow{k\rightarrow\infty} \frac{\theta_2}{2c \alpha \theta_1},
\end{equation}
where  $\theta_2$ is given in \eqref{theta}.
\end{theorem}
\begin{proof}
By \eqref{choice_alpha_PL} and \eqref{eq.cases.final}, it follows that
$$
\begin{aligned}
\E_k[f(x_{k+1})] - f(x_k)
\leq&\ -  \alpha \theta_1 \|\nabla f(x_k)\|^2+\theta_2,
\end{aligned}
$$
for all $k$ and $\theta_1,\, \theta_2$ defined in \eqref{theta}.
Therefore, \eqref{eq.PL} implies that
\begin{equation}
\begin{aligned}
\E_k[f(x_{k+1})] - f(x_k)
&\leq - 2c \alpha \theta_1 (f(x_k) - f_*) + \theta_2.
\end{aligned}
\end{equation}
Adding and subtracting $f_*$ on the left-hand side, taking total expectations, and rearranging give
\begin{equation}
\begin{aligned}
\E[f(x_{k+1})] - f_*  &\leq (1 - 2c \alpha \theta_1) (\E[f(x_{k})] - f_*) + \theta_2 \\
&=    \frac{\theta_2}{2c\alpha\theta_1} + (1 - 2c \alpha \theta_1) (\E[f(x_{k})] - f_*) + \theta_2 - \frac{\theta_2}{2c\alpha\theta_1} \\
&=    \frac{\theta_2}{2c\alpha\theta_1} + (1 - 2c \alpha \theta_1) \left(\E[f(x_{k})] - f_* - \frac{\theta_2}{2c\alpha\theta_1}\right).
\end{aligned}
\end{equation}
Since $1 - 2c\alpha\theta_1 \in (0,1)$, this represents a contraction inequality.  Applying the result repeatedly 
through iteration $k$, one obtains the desired result.
\qedhere
\end{proof}


\begin{theorem}\label{th.nonconvex_fixed_stepsize}
Under  Assumptions~\ref{ass.f}, \ref{hp_g} and \ref{hp_mu}, suppose that  \eqref{gamma1gamma2} holds
in case $h_2-\frac{1}{2}\omega>0$ and that 
\begin{equation}\label{choice_alpha}
\alpha\le  \frac{\omega}{\gamma_1LM_2},
\end{equation}
with $h_2$ given in (\ref{eq.ce_fg_upperbound}) and  $\theta_1$ given in (\ref{theta}).
Then the iterates generated by Algorithm TRishBB satisfy
\begin{equation}\label{eq.fixed_stepsize_nonconvex}
\sum_{k=0}^K\E \left[\frac{1}{K} \|\nabla f(x_{k})\|^2 \right] \xrightarrow{K\rightarrow\infty} 
\frac{\theta_2}{\alpha \theta_1},
\end{equation}
where $\theta_2$ is given in \eqref{theta}.
\end{theorem}

\begin{proof}
By \eqref{choice_alpha} and \eqref{eq.cases.final}, it follows that
$$
\begin{aligned}
\E_k[f(x_{k+1})] - f(x_k)
\leq&\ -  \alpha \theta_1 \|\nabla f(x_k)\|^2+\theta_2,
\end{aligned}
$$
for all $k$ and $\theta_1,\, \theta_2$ defined in \eqref{theta}.
Therefore, taking total expectation yields
\begin{equation}
\begin{aligned}
\E[\|\nabla f(x_{k})\|^2]
&\leq \frac{1}{\alpha \theta_1} \left(\E[f(x_k)] - \E[f(x_{k+1})]\right) + \frac{\theta_2}{\alpha \theta_1}.
\end{aligned}
\end{equation}
Summing both sides for $k\in \{1,2,\ldots, K\}$ gives
\begin{equation}
\begin{aligned}
\sum_{k=0}^K\E[\|\nabla f(x_{k})\|^2]&\leq \frac{1}{\alpha \theta_1} \left(f(x_0) - \E[f(x_{K+1})]\right)
+\frac{K\theta_2}{\alpha \theta_1}\nonumber \\
&\leq  \frac{f(x_0) - f_*}{\alpha \theta_1} +\frac{K\theta_2}{\alpha \theta_1}. \nonumber 
 \end{aligned}
\end{equation}
Dividing by $K$ concludes the proof.
\qedhere
\end{proof}
\subsection{Analysis with Unbiased Gradients}
In this subsection
we  assume that $g_k$ is unbiased estimator of the true gradient with bounded variance. Such assumption is summarized below.
\begin{assumption}\label{ass.gH}
For all $k$,  $g_k$ is an unbiased estimator of the gradient $\nabla f(x_k)$, i.e.,  $\E_k[g_k] = \nabla f(x_k)$.
For all $k$, there exist a positive constant $ M_g$ such that
\begin{equation}\label{exp_err}
    \E_k[\|\nabla f(x_k) - g_k\|^2] \leq M_g. 
\end{equation}
\end{assumption}
\noindent Since $g_k$ is unbiased, it follows.
\begin{equation}\label{eq.unbiased_consequence}
  \begin{aligned}
    \E_k[\|\nabla f(x_k) - g_k\|^2]
      &= \E_k[\|\nabla f(x_k)\|^2 - 2\nabla f(x_k)^Tg_k + \|g_k\|^2] \\
      &= -\|\nabla f(x_k)\|^2 + \E_k[\|g_k\|^2].
\end{aligned}
\end{equation}
Hence Assumption (\ref{hp_g}) is satisfied with $\omega=1$, $M_1=M_g$, $M_2=1$ and the convergence results of \Cref{sec_biased} hold. Now   we further characterize the convergence without imposing bounds,  such as inequality (\ref{gamma1gamma2}),  on $\gamma_1$, $\gamma_2$. The first lemma is a technical result.

\begin{lemma}\label{diff_f_g_unbiased}
Suppose that Assumption~\ref{ass.f} holds and
\begin{equation}\label{alpha_mu_final}
        0 < \alpha \leq \frac{\gamma_{2}}{8\gamma_{1}^2L},\qquad \mu_{min}\ge \frac{4\gamma_{1}\alpha}{5}.
\end{equation}
Then, the $k$-th iterate generated by Algorithm TRishBB satisfies
\begin{eqnarray}
& & \E_k[f(x_{k+1})] - f(x_k) \leq   -\frac{1}{16}\gamma_{2} \alpha   \E_k[\|g_k\|^2] + \frac{\gamma_{1}^2}{\gamma_{2}} \alpha \E_k[\|\nabla f(x_k) - g_k\|^2],  \qquad \label{decr_unbiased}\\
& & \E_k [f(x_{k+1})] - f(x_k) 
     \leq  \ - \frac{1}{16}  \gamma_2 \alpha \|\nabla f(x_k)\|^2 + \alpha\left(\frac{\gamma_1^2}{\gamma_2} -  \frac{1}{16}\gamma_2\right)M_g. \qquad \label{eq.for_later}
\end{eqnarray}
\end{lemma}

\begin{proof}
See \Cref{AppB2}.
\end{proof}

In case the objective function satisfies the Polyak-Lojasiewicz condition \eqref{ass.PL}, the 
expected optimality gap is bounded above by a sequence that converges linearly to a constant proportional to $M_g/c$.

\begin{theorem}\label{th.kl.fixed}
Suppose that Assumptions~\ref{ass.f}, \ref{ass.PL} and \ref{ass.gH} hold  and that
\begin{equation}\label{alpha_mu_final}
        0 < \alpha < \min \left \{  \frac{\gamma_{2}}{8\gamma_{1}^2L}, \frac{8}{\gamma_2 c}  \right \},\qquad \mu_{min}\ge \frac{4\gamma_{1}\alpha}{5}.
\end{equation}
 Then,  the sequence generated by Algorithm TRishBB satisfies
\begin{equation}
\E[f(x_{k+1})] - f_*  \xrightarrow{k\rightarrow\infty}  \frac{8\theta_3 M_g}{  \gamma_2 c },
\end{equation}
for all $k$, where  $\theta_3= \left(\frac{\gamma_1^2}{\gamma_2} -  \frac{\gamma_2}{16}  \right)$.

\end{theorem}

\begin{proof}
 
By \eqref{eq.for_later} and  \eqref{eq.PL} we obtain
\begin{equation}
\begin{aligned}
\E_k[f(x_{k+1})] - f(x_k)
&\leq -  \frac{1}{8} \gamma_2 \alpha c  (f(x_k) - f_*) +  \alpha \theta_3M_g.
\end{aligned}
\end{equation}
Taking total expectation, for all $k$ we have 
 
\begin{equation}
\begin{aligned}
\E[f(x_{k+1})] - f_*  &
\le    \frac{8\theta_3M_g}{ \gamma_2  c } +\left  (1 - \frac{1}{8} \gamma_2 \alpha c \right ) (\E[f(x_{k})] - f_*) +  \alpha \theta_3M_g - \frac{8\theta_3M_g}{ \gamma_2  c } \\
&=    \frac{8\theta_3M_g}{ \gamma_2  c } + \left(1 -\frac{1}{8} \gamma_2 \alpha c\right) \left(\E[f(x_{k})] - f_* -\frac{8\theta_3M_g}{ \gamma_2  c } \right).
\end{aligned}
\end{equation}
Since $(1 -\frac{1}{8} \gamma_2 \alpha c) \in (0,1)$, this represents a contraction inequality.  Applying the result repeatedly 
through iteration $k$, one obtains the desired result.
\qedhere
\end{proof}

When $f$ is a nonconvex function, the expected average squared norm of the gradient converges to $\left(\frac{16\gamma_1^2}{ \gamma_2^2} - 1\right)M_g$.
\begin{theorem}\label{th.nonconvex.fixed}
Suppose that Assumptions~\ref{ass.f}, \ref{ass.gH} hold  and that
\begin{equation}\label{alpha_mu_final}
        0 < \alpha \leq \frac{\gamma_{2}}{8\gamma_{1}^2L},\qquad   \mu_{min}\ge \frac{4\gamma_{1}\alpha}{5}.
\end{equation}
Then,  the sequence generated by Algorithm TRishBB satisfies
$$
\E \left[\frac{1}{K} \sum_{k=1}^K \|\nabla f(x_k)\|^2\right]  
     \xrightarrow{K\to\infty}  \left( \frac{16\gamma_1^2}{\gamma_2^2} - 1\right)M_g. 
$$
\end{theorem}

\begin{proof}
Consider (\ref{eq.for_later}). 
Taking total expectation, it follows for all $k$ that
\begin{equation}
    \E[f(x_{k+1})] - \E[f(x_k)] \leq -  \frac{1}{16}\gamma_2 \alpha \E[\|\nabla f(x_k)\|^2] + \alpha\left(\frac{\gamma_1^2}{\gamma_2} -  \frac{1}{16}\gamma_2\right)M_g,
\end{equation}
which implies
\begin{equation}
  \E[\|\nabla f(x_k)\|^2] \leq   \frac{16}{\gamma_2 \alpha} \left(\E[f(x_k)] - \E[f(x_{k+1})]\right) + 
	\left( \frac{16\gamma_1^2}{\gamma_2^2} - 1\right)M_g.
\end{equation}
The proof is concluded by summing this inequality over all $k \in \{1,\dots, K\}$ and using the fact that $f$ is bounded below by $f_*$.
\qedhere
\end{proof}


\section{Numerical Experiments}\label{NumericalSec}
In this section, we present experimental results to compare Algorithms \trishunop, \trishdue and \trishtre with the first-order version of TRish in \cite{curtis2019stochastic}, i.e., Algorithm TRishBB with $H_k=0$, $\forall k$. All experiments were conducted with MATLAB (R2024a) on a Rocky Linux 8.10 (64-bit) server with 256 GiB memory using an NVIDIA A100 PCIe GPU (with 80 GiB on-board memory).

\subsection{Classification Problems}\label{S.classProb}

We consider some ML and DL test problems of the form \eqref{eq1}. 
For binary classification problems using  \texttt{a1a}, \texttt{w1a}, and \texttt{cina} datasets \cite{chang2011libsvm}, we used logistic regression where the single loss function $f_i(x)$ with $x\in \mathbb{R}^n$ in \eqref{eq1} is defined as 
$$f_i(x) = \log\, (1+e^{-b_i(x^Ta_i)}),\quad i=1,\dots,N,\qquad b_i \in \{+1,-1\}.$$
 Note that $x\in \mathbb{R}^n$ and the data sample $a_i \in \mathbb{R}^d$ with the label $b_i$ have the same dimension, i.e., $n = d$. A numerically stable implementation of the logistic sigmoid was employed to prevent numerical overflow.

For multi-classification problems using \texttt{MNIST} \cite{deng2012mnist} and \texttt{CIFAR10} \cite{krizhevsky2009learning} datasets with $C$ classes, we used neural network models coupled with the softmax cross-entropy loss in \eqref{eq1}, whose  single loss function $f_i(x)$ is defined as 
$$f_i(x) = - \sum_{k=1}^{C} (b_i)_k \log (h(a_i; x))_k,\quad i=1,\dots,N,$$
where $a_i \in \mathbb{R}^d$ is a 3-D array (grayscale or RGB image),
 and the categorical label $b_i$ is converted into a one-hot encoded vector; i.e., it  is a binary vector of length $C$, where $C-1$ elements are zero except for the element 
whose index corresponds to the true class that is 
equal to one. The output $h(a_i; \cdot)$ is a prediction provided by the neural network (NN) whose architecture is described in \Cref{tab:nets} in Appendix~\ref{AppC}. The size $n$ of problem \eqref{eq1} using a NN depends on both the number of features $d$ and the depth of the network. We utilized the MATLAB Deep Learning Toolbox, leveraging automatic gradient computation in customized training loops; we used \texttt{forward} and \texttt{crossentropy} for implementing the NN softmax cross-entropy. See \cite[Table2]{krejic2024non} for details on implementing NN classification problems by this toolbox.  


More details on the test datasets and the architecture of the NNs used are available in Appendix~\ref{AppC}. 

\subsection{Experimental Configuration}\label{NumericalSec_configuration}
Algorithm \trishuno was run setting $m=20$ and selecting random and independent mini-batches ${\cal{N}}_k$, thus the gradients $g_k$ are unbiased.  Algorithm \trishdue was run with $\eta=0.9$ and  $m=N_b=\lfloor N/|\mathcal{N}_k|\rfloor$ on the base of the discussion in \Cref{Sub2.3}; the  mini-batches ${\cal{N}}_k$ were selected as disjoint subsets of ${\cal{N}}$, with ${\cal{N}}$ being shuffled at the beginning of each epoch; thus, the resulting gradients are biased. Algorithm \trishtre was run with $m_F=100$ as suggested in \cite{keskar2016adaqn}, its mini-batch selection was  as in \trishunop, and its cyclic iteration parameter $m$ was selected as in \trishduep.

For all $k\geq 0$, the mini-batch size for computing stochastic gradients in these variants was $|\mathcal{N}_k|=64$ for the datasets \texttt{a1a}, \texttt{w1a}, and \texttt{cina} and $|\mathcal{N}_k|=128$ for the datasets \texttt{MNIST} and \texttt{CIFAR10}. The initial BB steplength $\mu_0$ was set to 1 in all the algorithms; moreover, the values $\mu_{\min}$, and $\mu_{\max}$  were set to $10^{-5}$ and $10^{5}$, respectively. 

Regarding the predefined parameters $(\alpha, \gamma_1, \gamma_2)$, our settings were borrowed from those in \cite{curtis2019stochastic} which are accompanied by a detailed discussion. The scalar $\alpha$ takes value in a prefixed set and the triplets $(\alpha, \gamma_1, \gamma_2)$ used in our experiments are given in \Cref{tab:triplex} for each dataset. The parameters $\gamma_1$ and $\gamma_2$ are chosen by a gradient estimation as follows. We run the SG algorithm for one epoch with a fixed learning rate $\ell$ and compute the norm of the stochastic gradient at each iteration; then we set $G$ as the average of the norms computed throughout the epoch. The value $G$ estimates the magnitude of the stochastic gradients and promotes the occurrence of the three cases in 
(\ref{eq.step_TRish_first_order}) and (\ref{eq.step_TRishBB}), i.e., the use of normalized steps. Overall we used 60 and 36 triplets for binary classification and multi-class problems, respectively. 

For the binary and multi-class test problems,  given the triplet $(\alpha, \gamma_1, \gamma_2)$, we ran each algorithm 50 and 10 times using 50 and 10 predetermined random seeds, respectively.  We present the  testing accuracy, i.e.,  the percentage of samples in the testing set  correctly identified; in case of multiple runs, the testing accuracy is averaged over the runs. In the figures shown in the following subsection, we use {\tt ExIJK} to denote the results obtained in the {\tt Ex}periment where $\alpha$, $\gamma_1$, and $\gamma_2$ take the  {\tt I}-th, {\tt J}-th, and {\tt K}-th value indicated in  
\Cref{tab:triplex}, respectively.

In the solution of logistic regression problem, the initial guess was the zero vector, $x_0 = 0$. For solving the NN softmax cross-entropy problems, we used the Glorot (Xavier) initializer \cite{glorot2010understanding} for the weights, and zeros for the biases in the network's layers. All algorithms terminated after five number of epochs, namely when the total number of single gradient evaluations $\nabla f_i$, reached or exceeded $5 N$, where $N$ is the number of training samples. 

\begin{table}[t!]
\centering
\scriptsize
    \caption{Hyper-parameters for the experiments.}
    \label{tab:triplex}
   \begin{tabular}{|p{2cm}|p{7cm}|p{3cm}|}
    \hline
     & Hyper-parameters & $(G, \ell)$ \\
    \hline
    \texttt{ala} \newline \texttt{w1a} \newline \texttt{cina} & 
    $\alpha \in \{10^{-1}, 10^{-\frac{1}{2}}, 1, 10^{\frac{1}{2}}, 10\}$, \newline
    $\gamma_{1} \in \left\{\frac{4}{G}, \frac{8}{G}, \frac{16}{G}, \frac{32}{G}\right\}$, 
    $\gamma_{2} \in \left\{\frac{1}{2G}, \frac{1}{G}, \frac{2}{G}\right\}$ &
    $(0.3477, 0.1)$ \newline
    $(0.0887, 0.1)$ \newline
    $(0.0497, 0.1)$ \\
    \hline
    \texttt{MNIST} \newline \texttt{CIFAR10} & 
    $\alpha \in \{10^{-3}, 10^{-2}, 10^{-1}, 1\}$, \newline
    $\gamma_{1} \in \left\{\frac{4}{G}, \frac{8}{G}, \frac{16}{G}\right\}$, 
    $\gamma_{2} \in \left\{\frac{1}{8G}, \frac{1}{4G}, \frac{1}{2G}\right\}$ &
    $(0.2517, 0.01)$ \newline
    $(1.8375, 0.01)$ \\
    \hline
    \end{tabular}
\end{table}

\subsection{Numerical Results}
We now present and analyze the numerical results obtained. 
Table \ref{tab:2} shows the percentage of unconstrained steps of the form $p_k=-\mu_k g_k$ taken through the iterations within five epochs; the displayed data are averaged over the runs performed. We note that in some cases  this percentage is equal to zero,  thus the sequence of iterates $\{x_k\}$ generated by our TRishBB variants
coincides with the sequence generated by the first-order TRish. We also underline that the percentage of unconstrained steps increases with $\alpha$.  
In fact,  taking into account the form of the steps (\ref{eq.Solpk}), the larger $\alpha$ is, the higher the probability of  taking unconstrained steps $(-\mu_kg_k)$ is; as a consequence the performance of TRishBB variants differ considerably from the performance of TRish for large values of $\alpha$.

In Table \ref{tab:accuracy}, we report,  for every dataset and every value of $\alpha$, the highest testing accuracy obtained by TRish and the variants of TRishBB averaging over the tested values of $\gamma_1$ and $\gamma_2$; the reported values represent the largest testing accuracies among those reached at the end of each epoch.  For each algorithm we also report the ratio ("accuracy ratio") between the smallest and the largest accuracy reported in the corresponding line of the table. We can observe that the highest testing accuracy of the TRishBB variants is comparable to or larger than the accuracy of TRish. For the datasets {\tt a1a},  {\tt w1a}, and  {\tt cina}, the values of the testing accuracy ratio are close to one for \trishdue and \trishtrep, namely the highest testing accuracy obtained over five epochs is not highly sensitive to the value of $\alpha$. In contrast, the testing accuracy reached by the three versions of TRishBB and TRish on  the DL test problems is dependent on $\alpha$. In order to give more insight, in what follows we provide detailed statistics of the runs employing the largest value of $\alpha$
which promotes the use of the BB steplengths.

\begin{table}[t!] 
\centering
\caption{Percentage of BB steplengths taken with different $\alpha$.}
\label{tab:2}
\resizebox{\textwidth}{!}{
{\footnotesize 
\begin{tabular}{|p{3cm}|>{\centering\arraybackslash}m{2cm}|>{\centering\arraybackslash}m{2cm}|>{\centering\arraybackslash}m{2cm}|>{\centering\arraybackslash}m{2cm}|>{\centering\arraybackslash}m{2cm}|}
\hline
\texttt{ala} & $\alpha = 10^{-1}$ & $\alpha = 10^{-\frac{1}{2}}$ & $\alpha=1$ & $\alpha = 10^{\frac{1}{2}}$ & $\alpha = 10\,\,$ \\
\hline
{\trishuno}  & 00.83 & 50.00 & 81.67& 92.50 & 98.33 \\
\hline
{\trishdue}  & 03.17 & 67.46 & 88.89 & 96.03 & 99.21\\
\hline
{\trishtre}  & 00.00& 46.83 & 77.78 & 97.61 & 100 \\
\hline
\end{tabular}
}}
\medskip

\resizebox{\textwidth}{!}{
{\footnotesize 
\begin{tabular}{|p{3cm}|>{\centering\arraybackslash}m{2cm}|>{\centering\arraybackslash}m{2cm}|>{\centering\arraybackslash}m{2cm}|>{\centering\arraybackslash}m{2cm}|>{\centering\arraybackslash}m{2cm}|}
\hline
 \texttt{wla} & $\alpha = 10^{-1}$ & $\alpha = 10^{-\frac{1}{2}}$ & $\alpha=1$ & $\alpha = 10^{\frac{1}{2}}$ & $\alpha = 10\,\,$ \\
\hline
{\trishuno} & 03.80 & 11.96 & 27.72 & 61.62 & 79.35\\
\hline
{\trishdue} & 81.96 & 99.48 & 100 & 100 & 100\\
\hline
{\trishtre} & 53.61 & 78.46 & 97.42 & 100 & 100\\
\hline
\end{tabular}
}}
\medskip

\resizebox{\textwidth}{!}{
{\footnotesize 
\begin{tabular}{|p{3cm}|>{\centering\arraybackslash}m{2cm}|>{\centering\arraybackslash}m{2cm}|>{\centering\arraybackslash}m{2cm}|>{\centering\arraybackslash}m{2cm}|>{\centering\arraybackslash}m{2cm}|}
\hline
\texttt{cina}\, & $\alpha = 10^{-1}$ & $\alpha = 10^{-\frac{1}{2}}$ & $\alpha=1$ & $\alpha = 10^{\frac{1}{2}}$ & $\alpha = 10\,\,$ \\
\hline
{\trishuno} & 14.11 & 60.08 & 88.04 & 90.59 & 89.65\\
\hline
{\trishdue} & 71.48 & 99.74 & 100 & 100 & 100\\
\hline
{\trishtre} & 23.66 & 44.63 & 58.8 &  87.60& 99.87\\
\hline
\end{tabular}
}}
\medskip

\resizebox{\textwidth}{!}{
{\footnotesize 
\begin{tabular}{|p{3.5cm}|>{\centering\arraybackslash}m{2.5cm}|>{\centering\arraybackslash}m{2.5cm}|>{\centering\arraybackslash}m{2.5cm}|>{\centering\arraybackslash}m{2.5cm}|}
\hline
\texttt{MNIST} & $\alpha = 10^{-3}$ & $\alpha = 10^{-2}$ & $\alpha=10^{-1}$ & $\alpha = 1$\\
\hline
{\trishuno} & 00.00 & 00.36 & 66.67 & 96.10\\
\hline
{\trishdue} & 00.00 & 00.00 & 00.13 & 87.88\\
\hline
{\trishtre} & 00.00 &  00.00&  00.00& 15.23\\
\hline
\end{tabular}
}}
\medskip

\resizebox{\textwidth}{!}{
{\footnotesize 
\begin{tabular}{|p{3.5cm}|>{\centering\arraybackslash}m{2.5cm}|>{\centering\arraybackslash}m{2.5cm}|>{\centering\arraybackslash}m{2.5cm}|>{\centering\arraybackslash}m{2.5cm}|}
\hline
\texttt{CIFAR10} & $\alpha = 10^{-3}$ & $\alpha = 10^{-2}$ & $\alpha=10^{-1}$ & $\alpha = 1$\\
\hline
{\trishuno} & 00.00 & 00.00 & 24.93 & 81.62 \\
\hline
{\trishdue} & 00.00 & 00.00 & 00.00 &  10.80\\
\hline
{\trishtre} &  00.00&  00.00& 00.00 &  00.78\\
\hline
\end{tabular}
}}
\end{table}

\begin{table}[t!] 
\centering
\caption{Highest average testing accuracy with different $\alpha$ (in \%) and the accuracy ratio.}
\label{tab:accuracy}
\resizebox{\textwidth}{!}{
{\small 
\begin{tabular}{|p{2.5cm}|>{\centering\arraybackslash}m{1.8cm}|>{\centering\arraybackslash}m{1.8cm}|>{\centering\arraybackslash}m{1.8cm}|>{\centering\arraybackslash}m{1.8cm}|>{\centering\arraybackslash}m{2cm}|>{\centering\arraybackslash}m{3cm}|}
\hline
\texttt{ala} & $\alpha = 10^{-1}$ & $\alpha = 10^{-\frac{1}{2}}$ & $\alpha=1$ & $\alpha = 10^{\frac{1}{2}}$ & $\alpha = 10$ & \text{accuracy ratio}\,\, \\
\hline
{\trishuno}  &83.76 & 83.55 & 83.09 & 82.54 & 82.40 & 0.9837 \\
\hline
{\trishdue}  &83.87 & 83.67 & 83.52 & 83.52 & 83.52 & 0.9958\\
\hline
{\trishtre}  & 83.75 &  83.49 & 83.40 & 83.40  & 83.40 &   0.9958\\
\hline
{TRish}      & 83.75 & 83.51 & 81.55 & 79.50 & 79.06 & 0.9440\\
\hline
\end{tabular}
}}
\medskip

\resizebox{\textwidth}{!}{
{\small 
\begin{tabular}{|p{2.5cm}|>{\centering\arraybackslash}m{1.8cm}|>{\centering\arraybackslash}m{1.8cm}|>{\centering\arraybackslash}m{1.8cm}|>{\centering\arraybackslash}m{1.8cm}|>{\centering\arraybackslash}m{2cm}|>{\centering\arraybackslash}m{3cm}|}
\hline
\texttt{wla} & $\alpha = 10^{-1}$ & $\alpha = 10^{-\frac{1}{2}}$ & $\alpha=1$ & $\alpha = 10^{\frac{1}{2}}$ & $\alpha = 10$ & \text{accuracy ratio}\,\, \\
\hline
{\trishuno}  & 89.61 & 89.67 & 89.57 & 89.23 & 88.83 & 0.9906\\
\hline
{\trishdue}  & 89.59 & 89.58 & 89.57 & 89.57 & 89.57 & 0.9998\\
\hline
{\trishtre} & 89.60 & 89.60 & 89.51 &  89.51 &  89.51& 0.9989\\
\hline
{TRish}  & 89.62 & 89.67 & 89.57 & 89.30 & 88.99 & 0.9924\\
\hline
\end{tabular}
}}
\medskip

\resizebox{\textwidth}{!}{
{\small 
\begin{tabular}{|p{2.5cm}|>{\centering\arraybackslash}m{1.8cm}|>{\centering\arraybackslash}m{1.8cm}|>{\centering\arraybackslash}m{1.8cm}|>{\centering\arraybackslash}m{1.8cm}|>{\centering\arraybackslash}m{2cm}|>{\centering\arraybackslash}m{3cm}|}
\hline
\texttt{cina} & $\alpha = 10^{-1}$ & $\alpha = 10^{-\frac{1}{2}}$ & $\alpha=1$ & $\alpha = 10^{\frac{1}{2}}$ & $\alpha = 10$ & \text{accuracy ratio}\,\, \\
\hline
{\trishuno}  & 91.70 & 91.63 & 91.37 & 90.84 & 90.29 & 0.9846\\
\hline
{\trishdue}  & 91.49 & 91.45 & 91.45 & 91.45 & 91.45 & 0.9996\\
\hline
{\trishtre}  & 91.59 &  91.56&  91.18&  91.18& 91.18 & 0.9955\\
\hline
{TRish}      & 91.64 & 91.68 & 90.63 & 87.93 & 87.51 & 0.9545\\
\hline
\end{tabular}
}}
\medskip

\resizebox{\textwidth}{!}{
{\small 
\begin{tabular}{|p{3cm}|>{\centering\arraybackslash}m{2.2cm}|>{\centering\arraybackslash}m{2.2cm}|>{\centering\arraybackslash}m{2.2cm}|>{\centering\arraybackslash}m{2.2cm}|>{\centering\arraybackslash}m{3cm}|}
\hline
\texttt{MNIST} & $\alpha = 10^{-3}$ & $\alpha = 10^{-2}$ & $\alpha=10^{-1}$ & $\alpha = 1$ & \text{accuracy ratio}
\\
\hline
{\trishuno}  & 89.74 & 96.83 & 98.70 & 98.72 &  0.9090\\
\hline
{\trishdue}  & 87.27 & 96.98 & 98.95 & 98.81 & 0.8820 \\
\hline
{\trishtre} & 87.21& 96.95 &  98.89&  98.59& 0.8819\\
\hline
{TRish}  & 87.21 & 96.95 & 98.87 & 98.42 & 0.8821\\
\hline
\end{tabular}
}}
\medskip

\resizebox{\textwidth}{!}{
{\small 
\begin{tabular}{|p{3cm}|>{\centering\arraybackslash}m{2.2cm}|>{\centering\arraybackslash}m{2.2cm}|>{\centering\arraybackslash}m{2.2cm}|>{\centering\arraybackslash}m{2.2cm}|>{\centering\arraybackslash}m{3cm}|}
\hline
\texttt{CIFAR10} & $\alpha = 10^{-3}$ & $\alpha = 10^{-2}$ & $\alpha=10^{-1}$ & $\alpha = 1$ & \text{accuracy ratio}
\\
\hline
{\trishuno}  & 30.86 & 44.72 & 59.19 & 65.41 & 0.4718 \\
\hline
{\trishdue}  & 31.15 & 44.95 & 60.74 & 65.98 & 0.4721\\
\hline
{\trishtre}  & 31.09 & 45.06 & 60.28 & 64.89 &  0.4792\\
\hline
{TRish}      & 31.08 & 45.04 & 60.27 & 64.66 & 0.4807\\
\hline
\end{tabular}
}}
\end{table}

\subsubsection{Datasets {\tt a1a},  {\tt w1a}, and  {\tt cina}}  \label{S.a1aw1acina}
We present the averaged accuracies obtained with $\alpha=10$  for each of the five epochs and for each  triplet $(\alpha, \gamma_1, \gamma_2)=(10, \gamma_1, \gamma_2)$. For this value of $\alpha$, our algorithms  take the unconstrained step in the large majority of the iterations, and thus differ from TRish considerably, see Table \ref{tab:2}.

\begin{figure}
    \vspace*{-50pt}
    \centering
    \begin{adjustbox}{width=1.0\linewidth, center}
    \includegraphics{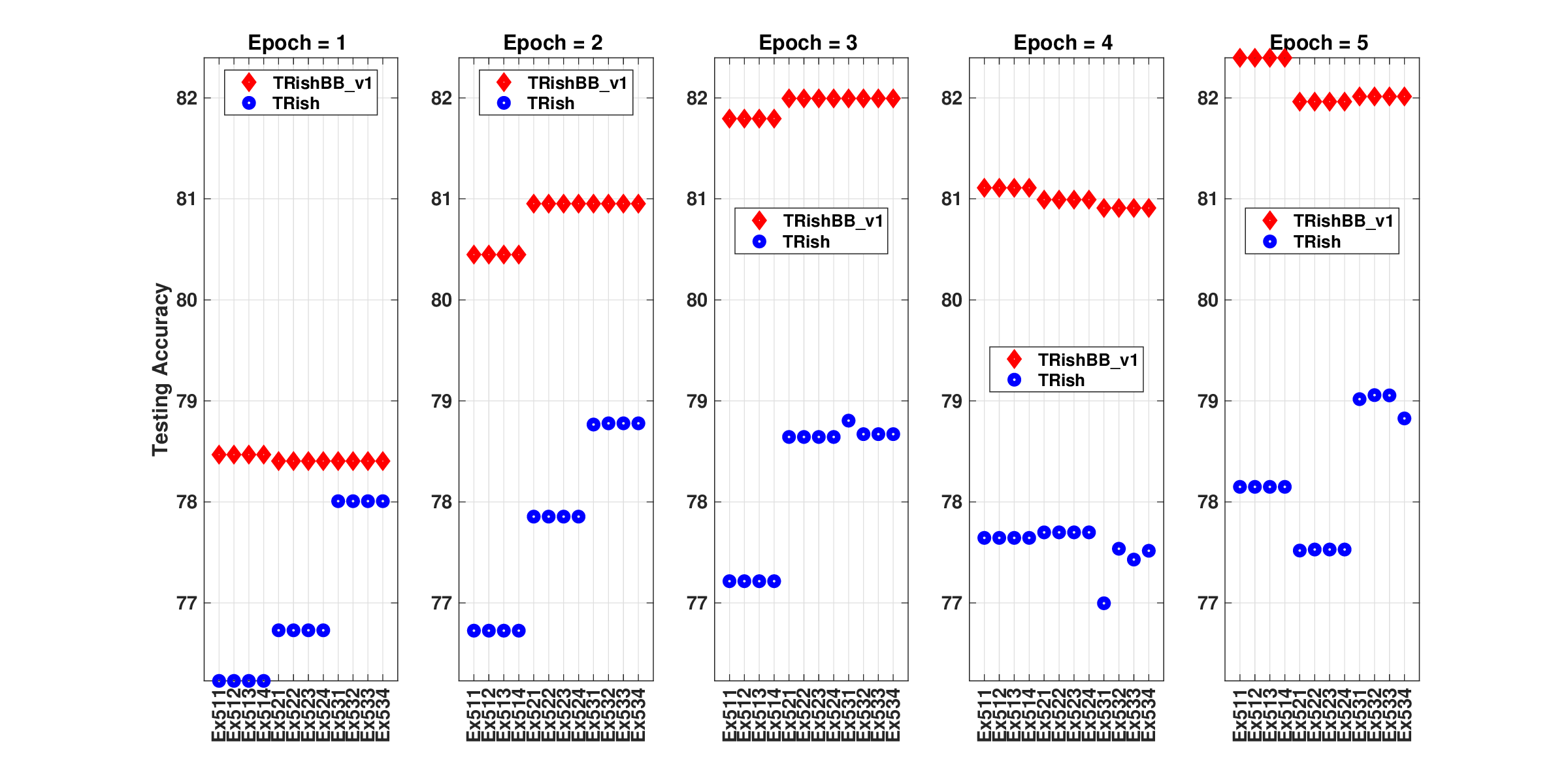}
    \end{adjustbox}
    
    \centering
    \begin{adjustbox}{width=1.0\linewidth, center}
\includegraphics{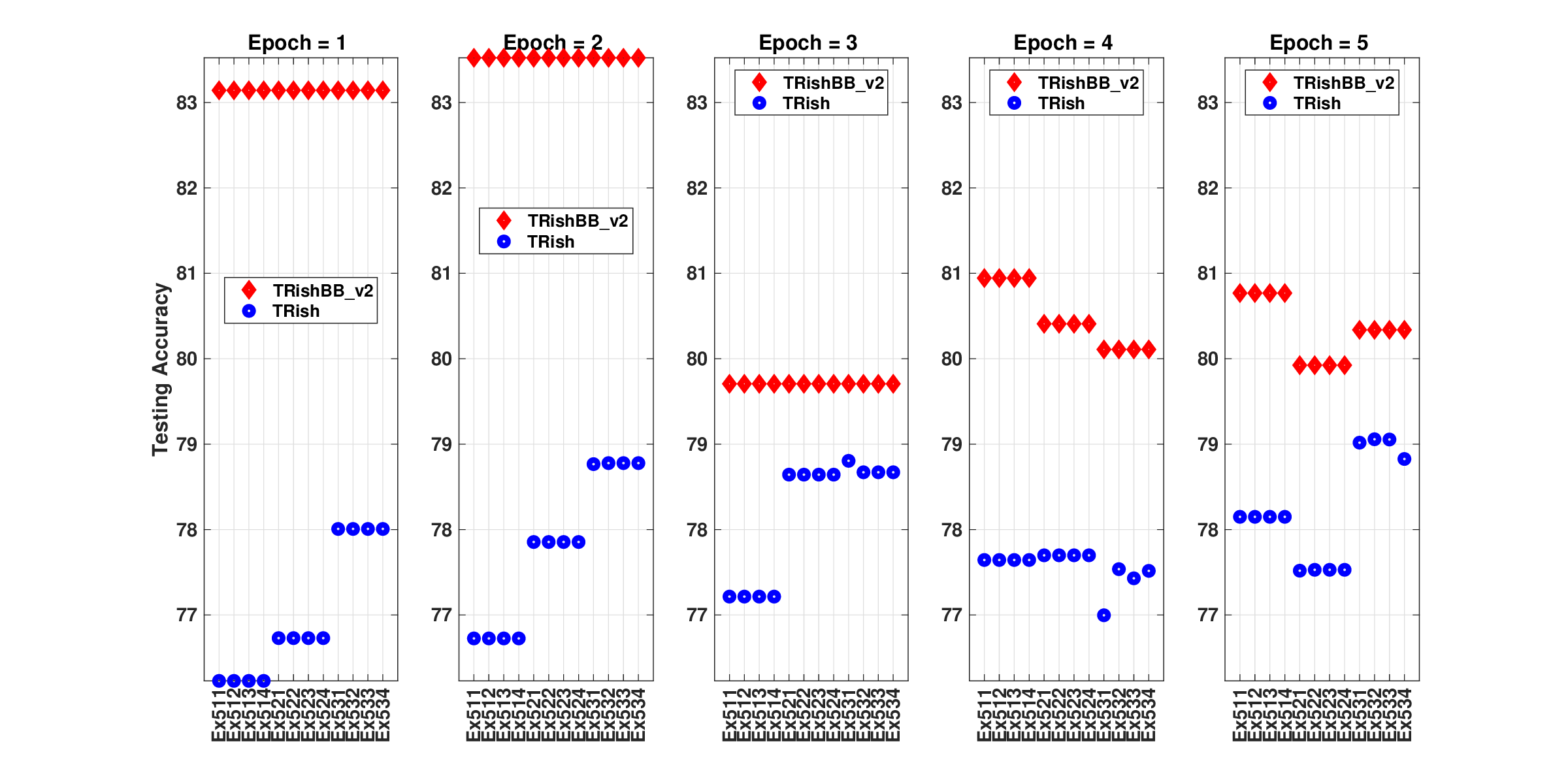}
    \end{adjustbox}
  
    \centering
    \begin{adjustbox}{width=1.0\linewidth, center}
    \includegraphics{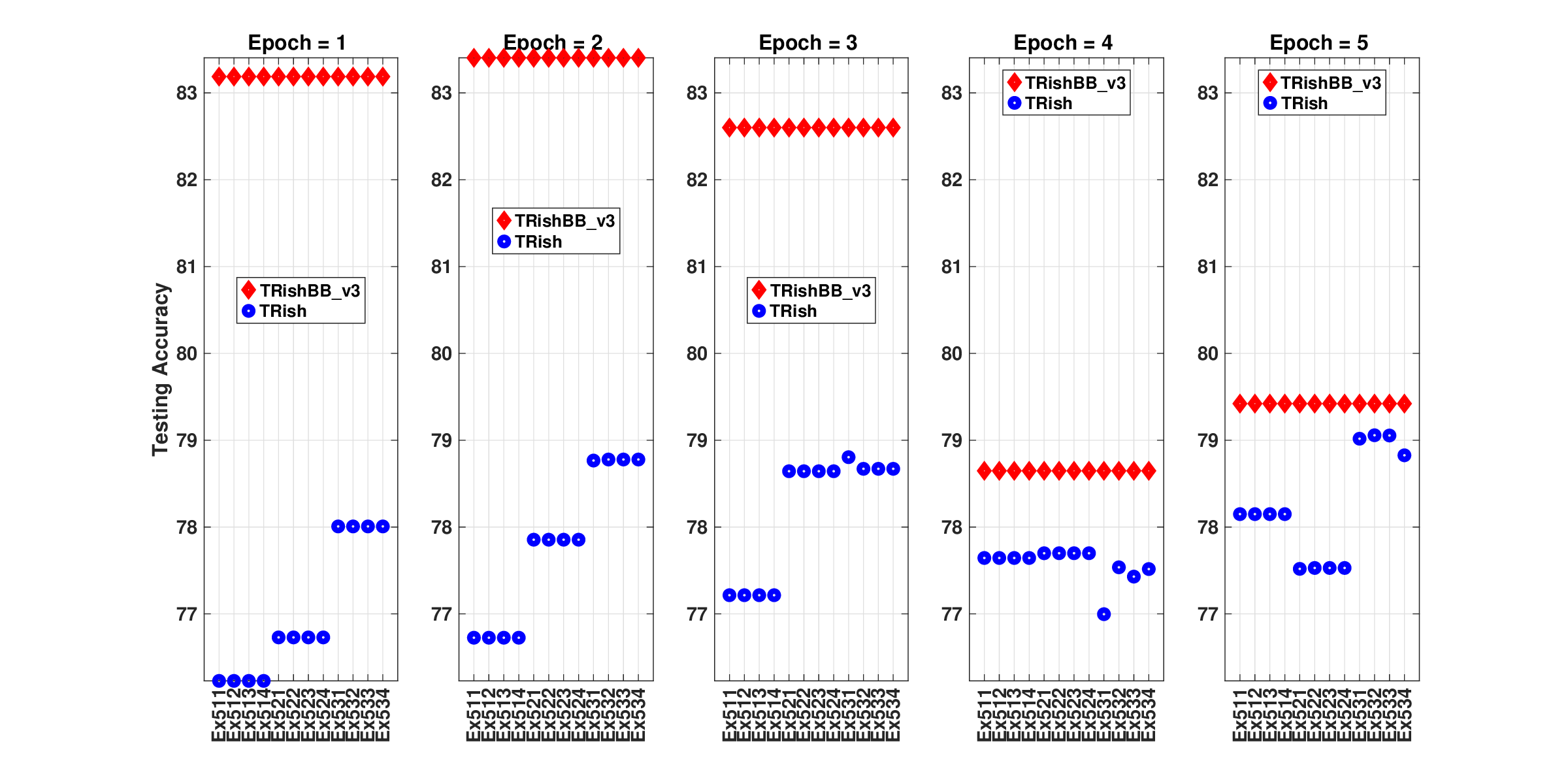}
    \end{adjustbox}
    \caption{\footnotesize{\texttt{a1a}: Average accuracy, $(\alpha, \gamma_1, \gamma_2)=(10, \gamma_1, \gamma_2)$, $|\mathcal{N}_k|=64$. 
    Top: \trishunop,   $m=20$.
    Middle:  \trishduep,  $m=25$.
    Bottom:  \trishtrep,  $m=25$. }}
    \label{fig:ala}
\end{figure}

\Cref{fig:ala} illustrates the results obtained on the dataset {\tt a1a}. The results displayed indicate that the highest accuracy is achieved by \trishdue and \trishtre in the first two epochs and that all TRishBB variants outperform TRish. Remarkably the testing accuracy reached by TRishBB is quite insensitive to the tested values of $\gamma_1$ and $\gamma_2$ contrary to TRish.


\begin{figure}
\vspace*{-50pt}
    \centering
    \begin{adjustbox}{width=1.0\linewidth, center}
    \includegraphics{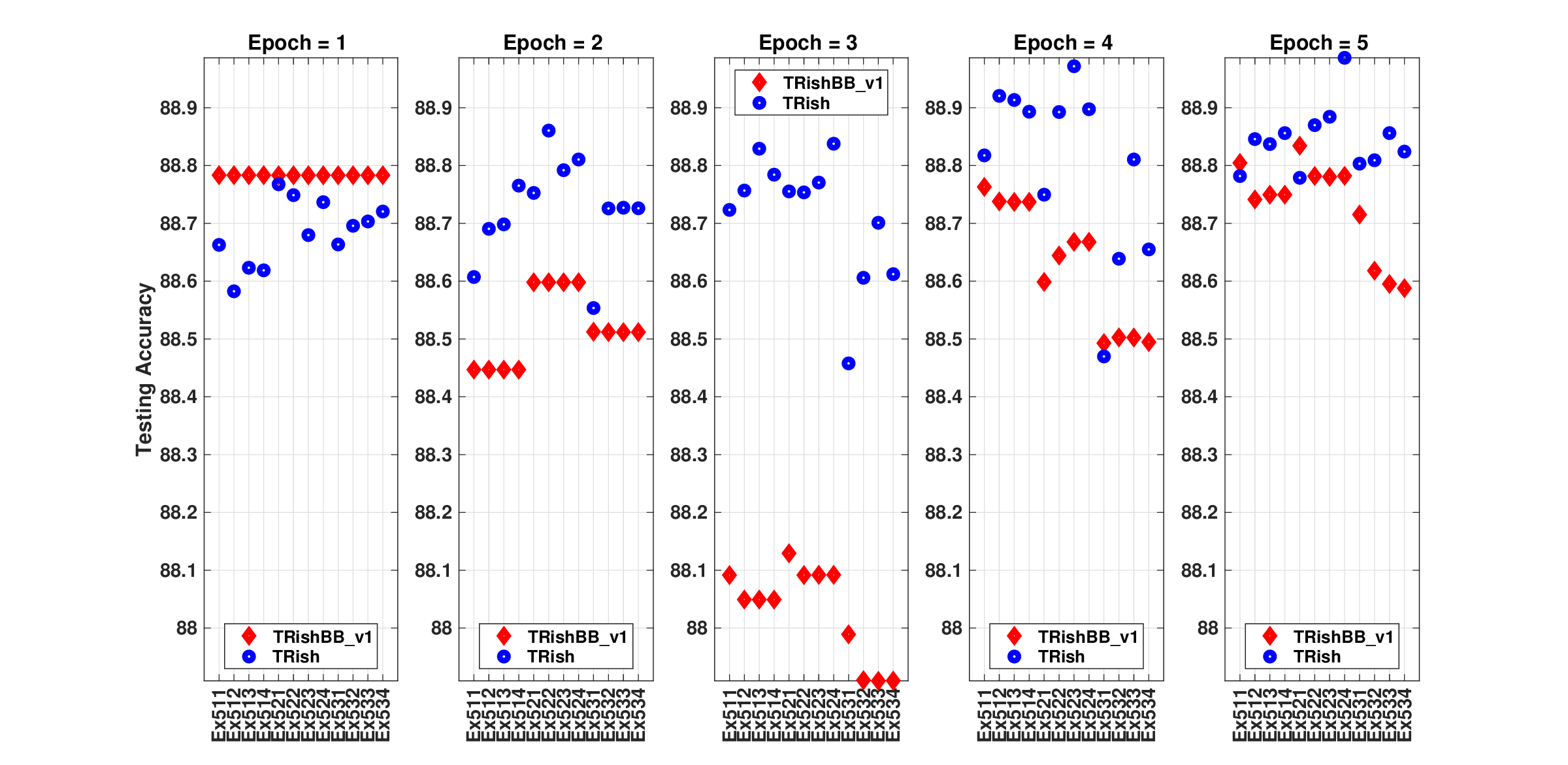}
    \end{adjustbox}    
    \centering
    \begin{adjustbox}{width=1.0\linewidth, center}
   \includegraphics{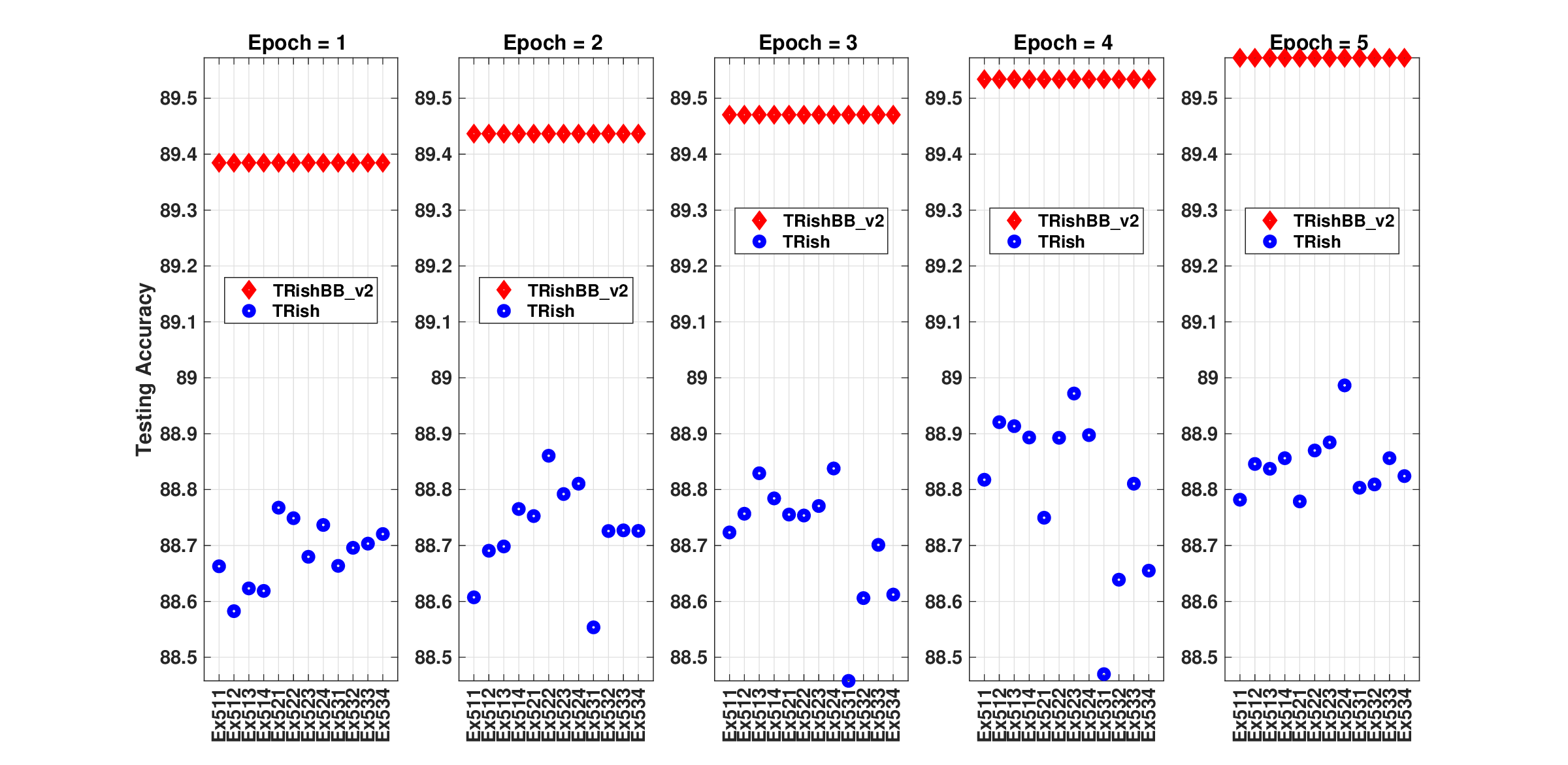}
    \end{adjustbox}    
    \centering
    \begin{adjustbox}{width=1.0\linewidth, center}
    \includegraphics{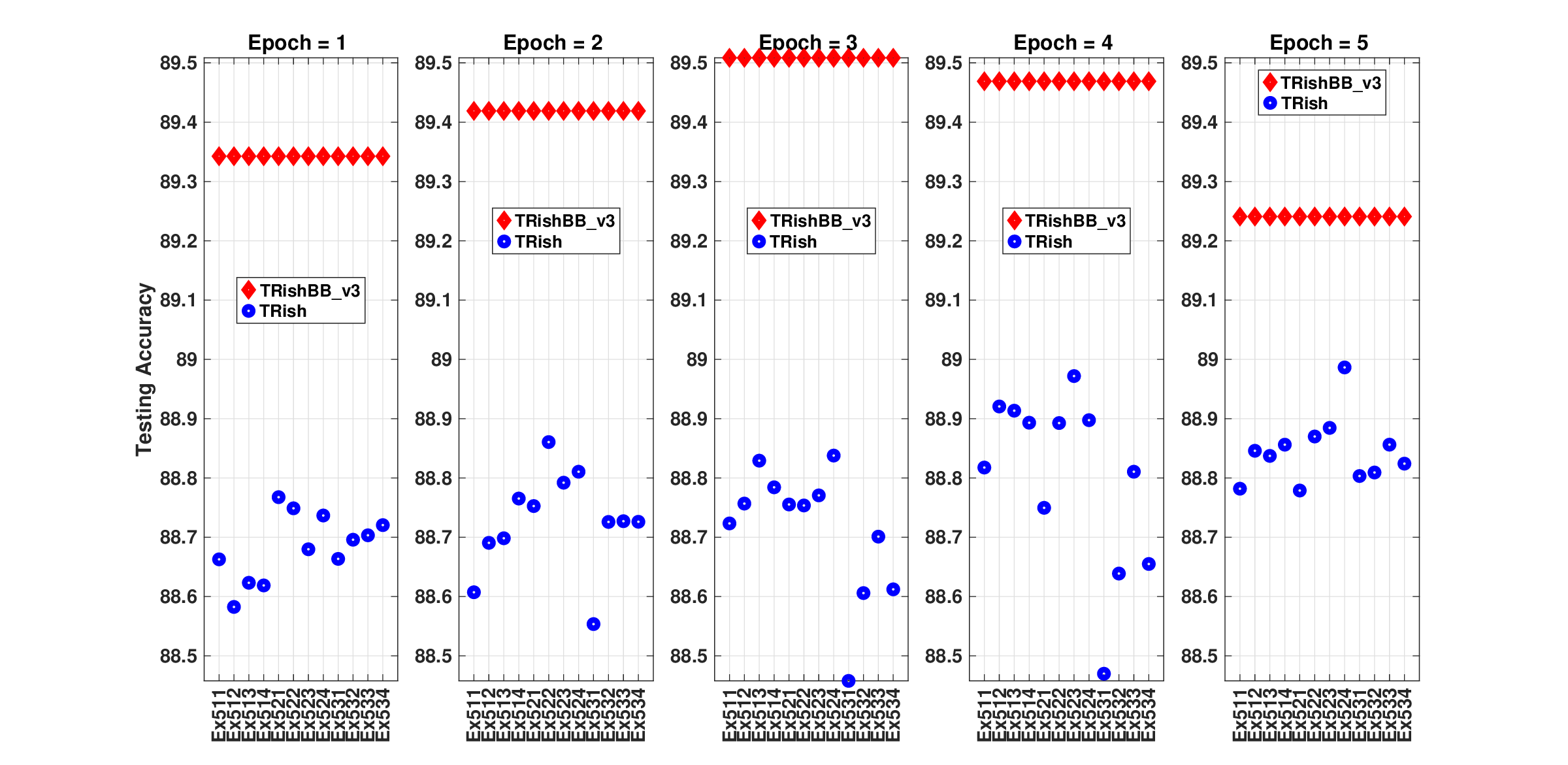}
    \end{adjustbox}    
    \caption{\footnotesize{ \texttt{w1a}: Average accuracy, $(\alpha, \gamma_1, \gamma_2)=(10, \gamma_1, \gamma_2)$, $|\mathcal{N}_k|=64$. Top: \trishunop,  $m=20$.
    Middle: \trishduep, $m=38$.
    Bottom: \trishtrep, $m=38$. }}
    \label{fig:wla}
\end{figure}

\Cref{fig:wla} concerns the dataset {\tt w1a}. Algorithms \trishdue and \trishtre were very effective and their performance is insensitive to the triplet  used. We note that the steplength $\mu_k$ was selected at every iteration by these algorithms, see Table \ref{tab:2}. On the other hand, the steplength $\mu_k$ was selected  approximately $79\%$ of the iterations in  \trishuno and its performance is comparable to the performance of TRish.

\begin{figure}
\vspace*{-50pt}
    \centering
    \begin{adjustbox}{width=1.0\linewidth, center}
    \includegraphics{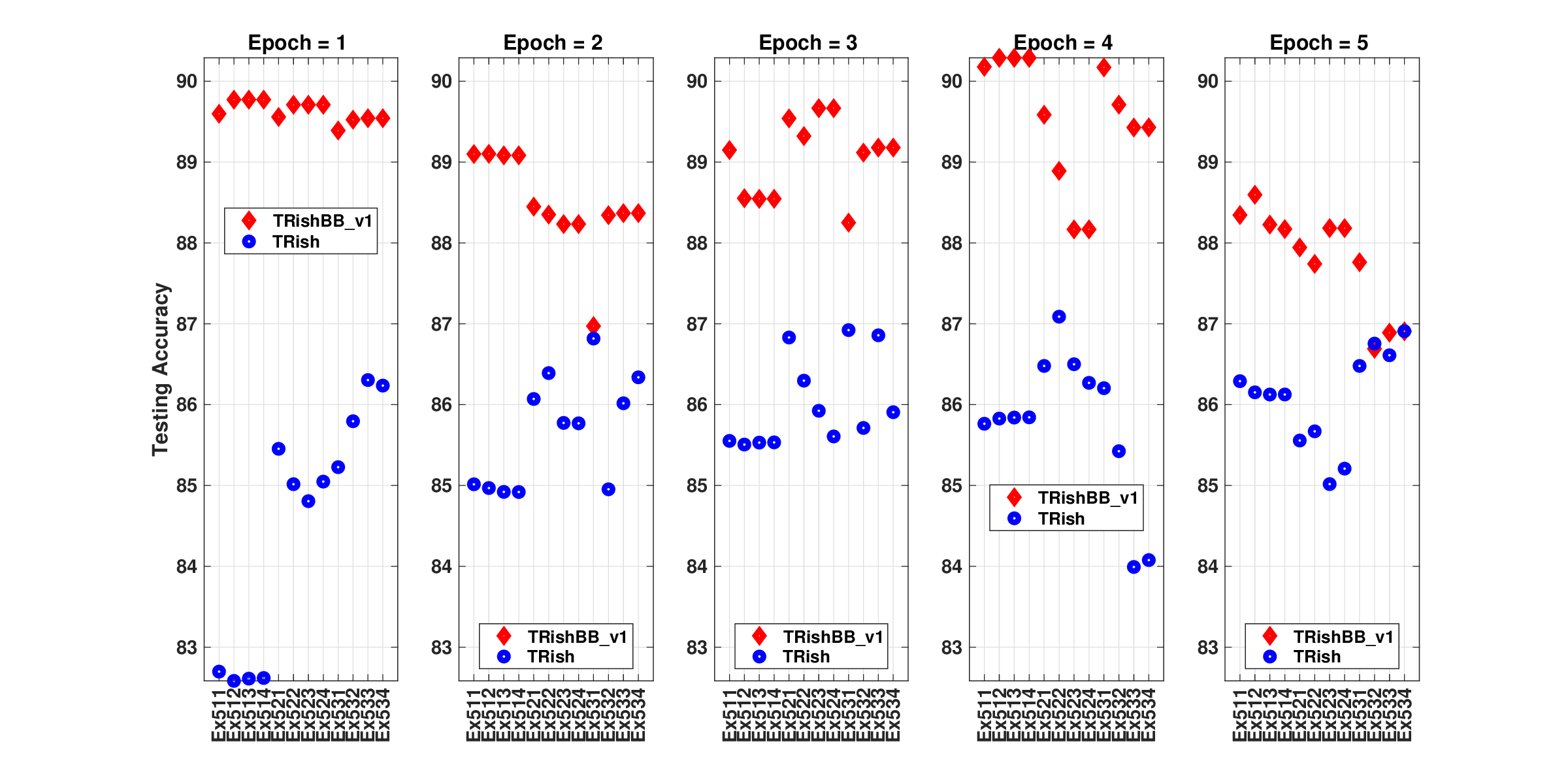}
    \end{adjustbox}    
    \centering
    \begin{adjustbox}{width=1.0\linewidth, center}  \includegraphics{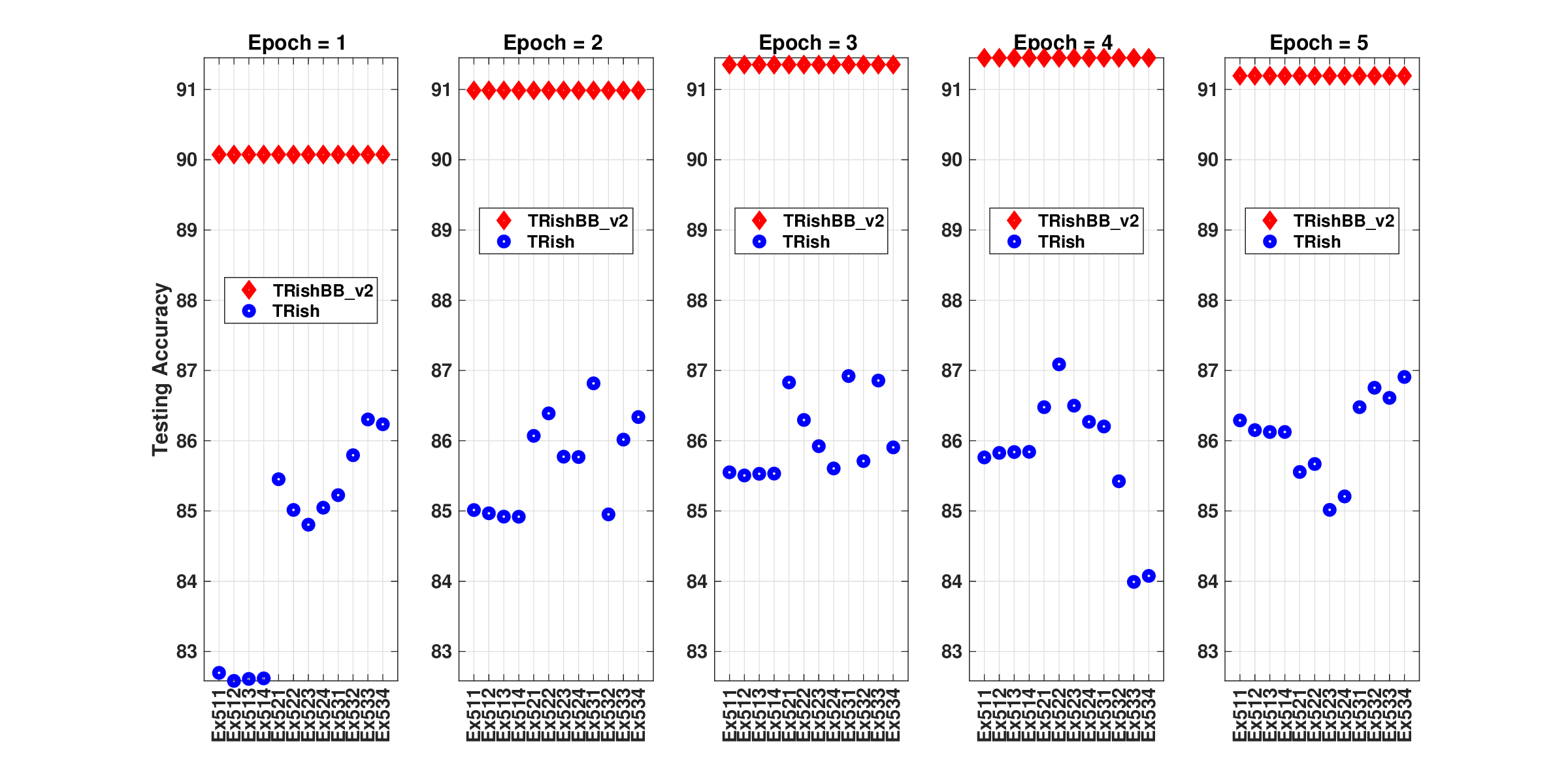}
    \end{adjustbox}    
    \centering
    \begin{adjustbox}{width=1.0\linewidth, center}
    \includegraphics{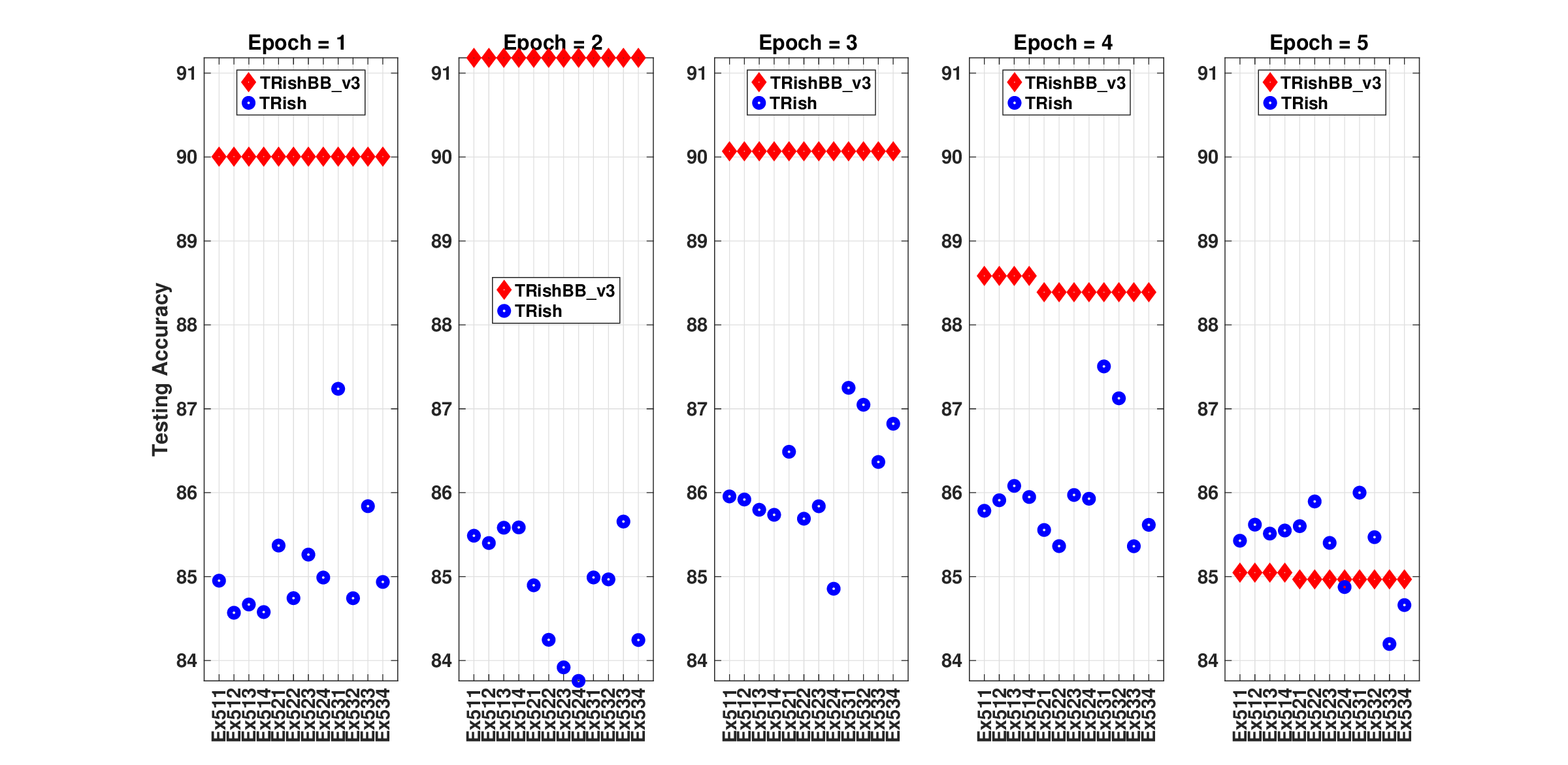}
    \end{adjustbox}    
    \caption{\footnotesize{ \texttt{cina}: Average accuracy, $(\alpha, \gamma_1, \gamma_2)=(10, \gamma_1, \gamma_2)$,  $|\mathcal{N}_k|=64$. 
    Top: \trishunop, $m=20$;  Middle: \trishduep, $m=156$.
    Bottom: \trishtrep, $m=156$.}}
    \label{fig:cina}
\end{figure}

Results for the dataset {\tt cina} are displayed  in  \Cref{fig:cina}. 
From Table \ref{tab:2} we know that \trishuno selected the steplength $\mu_k$ approximately at the $89\%$ of the iterations,  \trishdue selected  $\mu_k$ at all iterations and \trishtre selected $\mu_k$ around at the  $99.87\%$ of all iterations. 
\trishuno generally achieves higher accuracy than TRish; both algorithms  are quite sensitive to the choice of $\gamma_1$ and $\gamma_2$ but the dependence of \trishuno by the parameters is milder. In contrast, \trishdue and \trishtre exhibit (almost) steady behavior with respect to the parameters and achieve high accuracy.

\begin{figure}[t]
    \centering
    \includegraphics[width=1.05\textwidth, height=5cm]{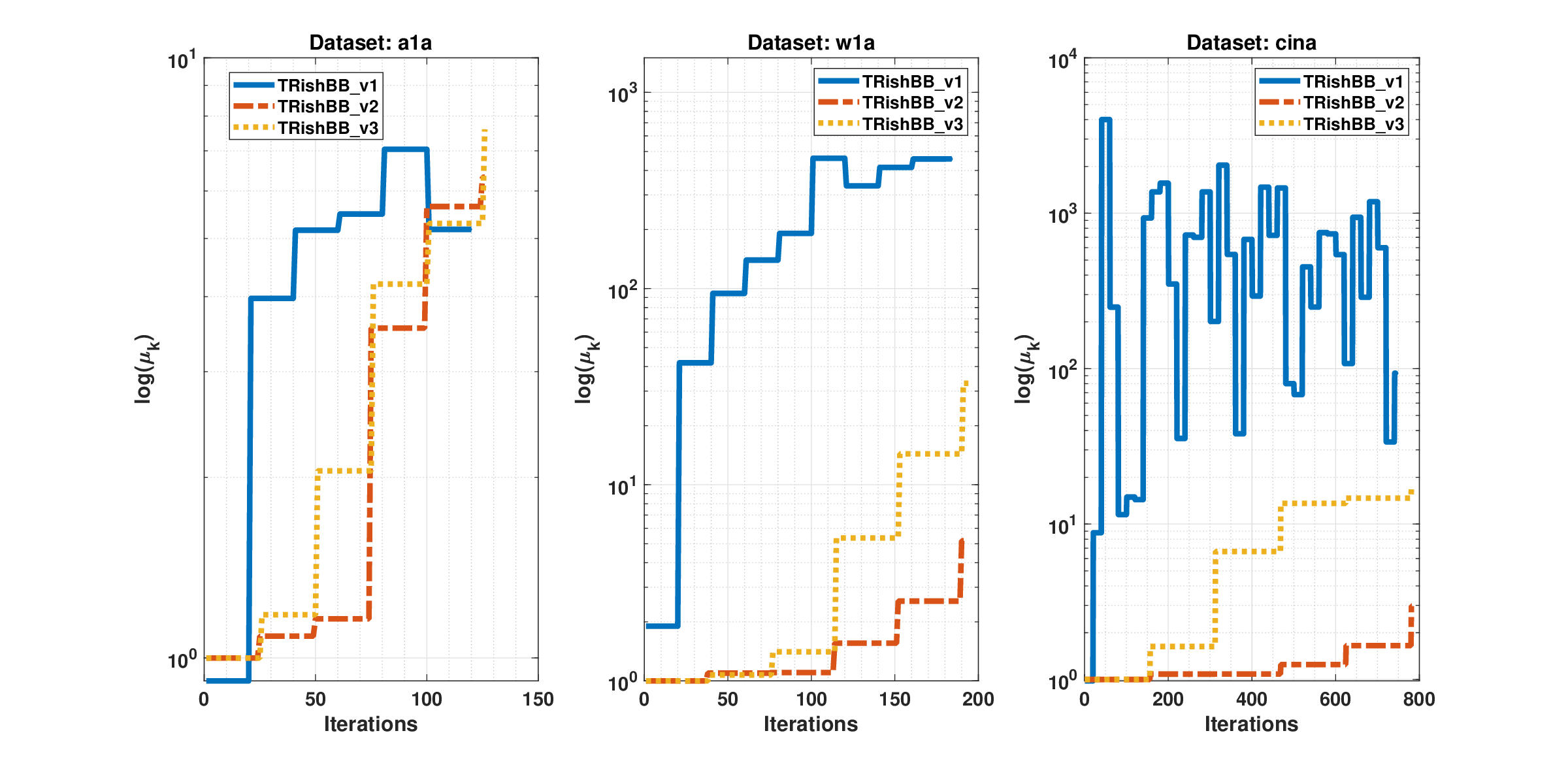}
    
    \vspace{\baselineskip}
    
    \includegraphics[width=0.75\textwidth, height=5cm]{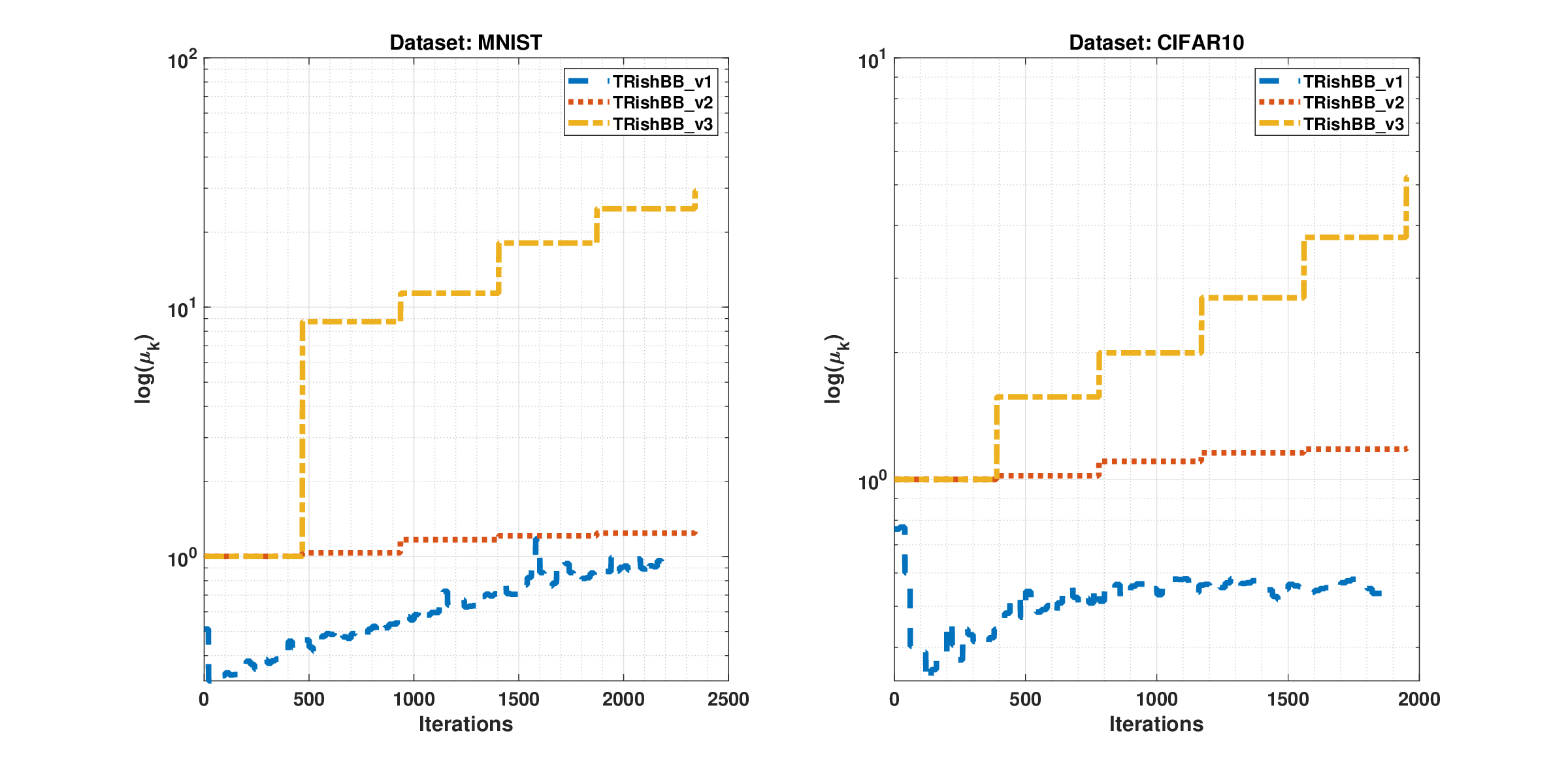}
    
    \caption{\label{fig:muA_alpha5} Average values of $\mu_k$
    with triplets $(\alpha,\gamma_1, \gamma_2)$ along the iterations. Top: \texttt{a1a}, \texttt{w1a} and \texttt{cina} with $\alpha=10$. Bottom:  \texttt{MNIST} and \texttt{CIFAR10} with $\alpha = 1$.} 
\end{figure}

The results described above indicate that  exploiting stochastic BB step\-lengths  and  guidelines on stochastic quasi-Newton methods improves the overall performance of TRish.
\Cref{fig:muA_alpha5} shows the average values of the BB parameters built with $\alpha=10$ and  varying $\gamma_1$ and $\gamma_2$ for each dataset. It is evident that the rule for generating $\mu_k$ in \trishdue and \trishtre provides smaller values of it,  enhancing the selection of unconstrained steps $(-\mu_k g_k)$. Algorithm \trishdue provides the best accuracy and its performance is quite insensitive to the triplet used.

\subsubsection{ Datasets {\tt MNIST} and {\tt CIFAR10}} \label{S.MnistCifar10}
We present the averaged accuracies obtained for each of the five epochs performed in the experiments and for each  triplet $(\alpha, \gamma_1, \gamma_2)=(1, \gamma_1, \gamma_2)$. 
As  shown in Table \ref{tab:2}, the variants of TRishBB employed the unconstrained trust-region solution less frequently than for datasets  \texttt{a1a}, \texttt{w1a} and \texttt{cina},  but  the percentage of unconstrained steps increases with the value of $\alpha$.

\begin{figure}
\vspace*{-50pt}
    \centering
    \begin{adjustbox}{width=1.\linewidth, center}
    \includegraphics{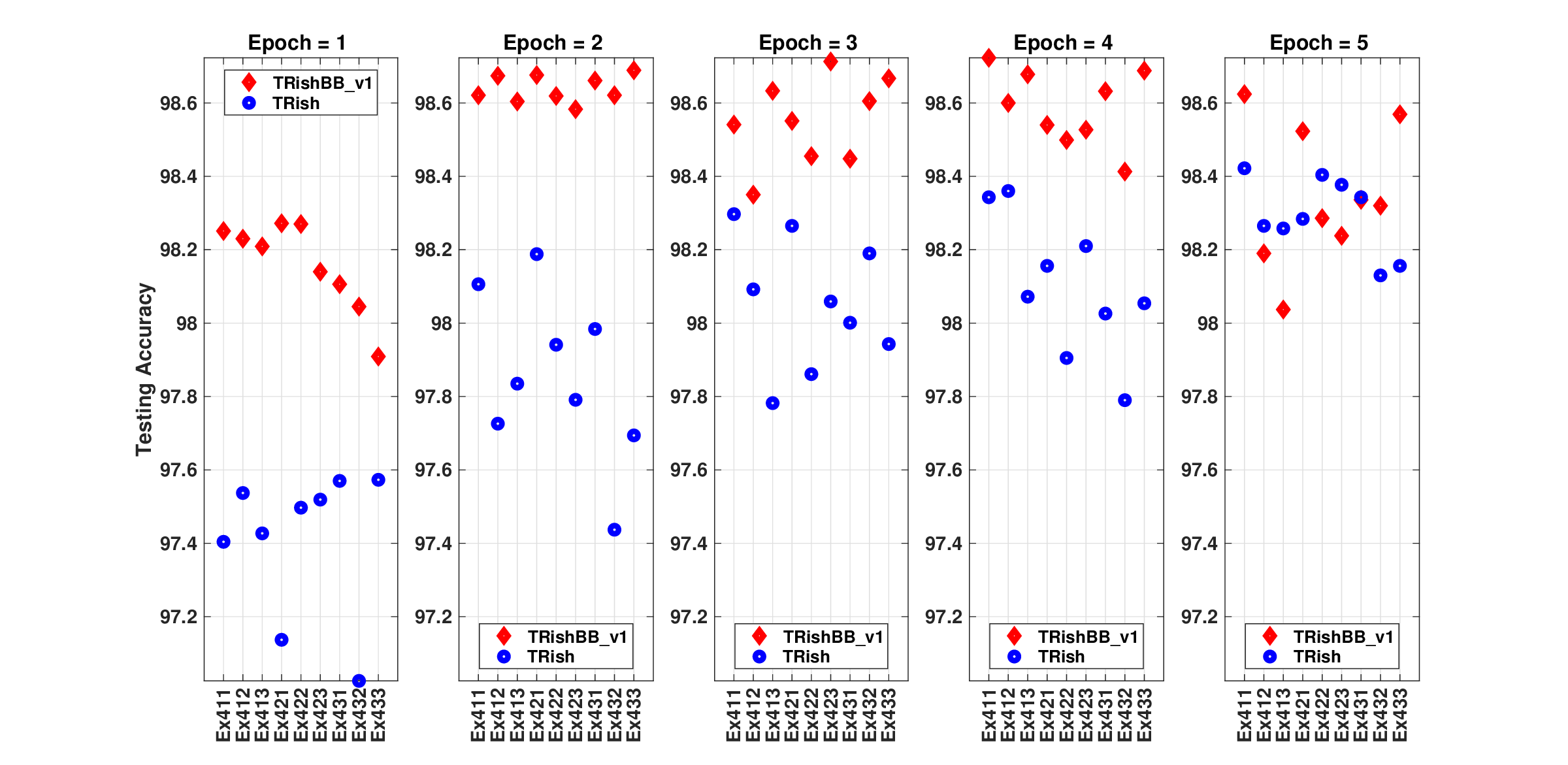}
    \end{adjustbox}    

    \centering
    \begin{adjustbox}{width=1.0\linewidth, center}
    \includegraphics{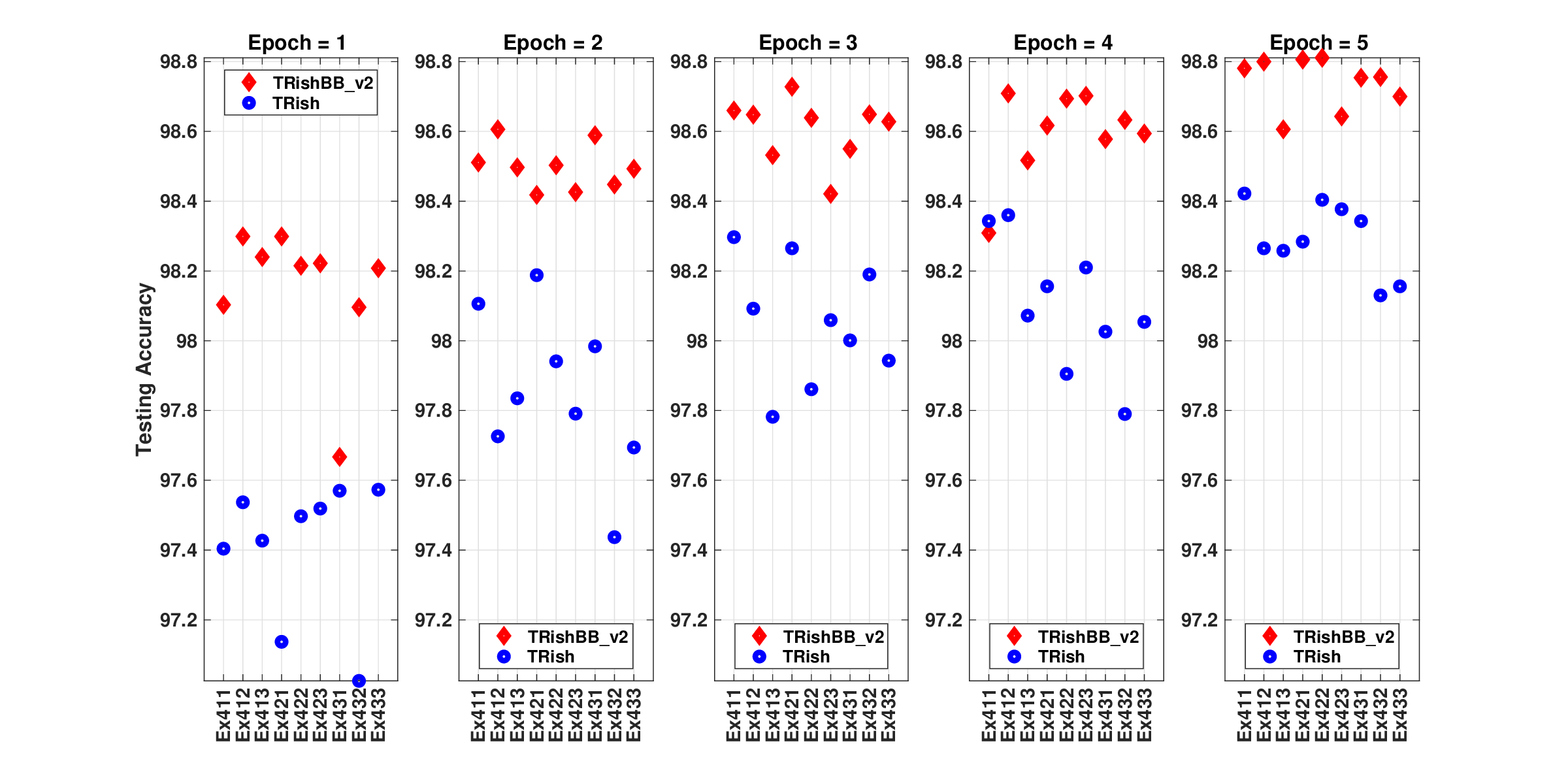}
    \end{adjustbox}    
    \centering
        \begin{adjustbox}{width=1.0\linewidth, center}
    \includegraphics{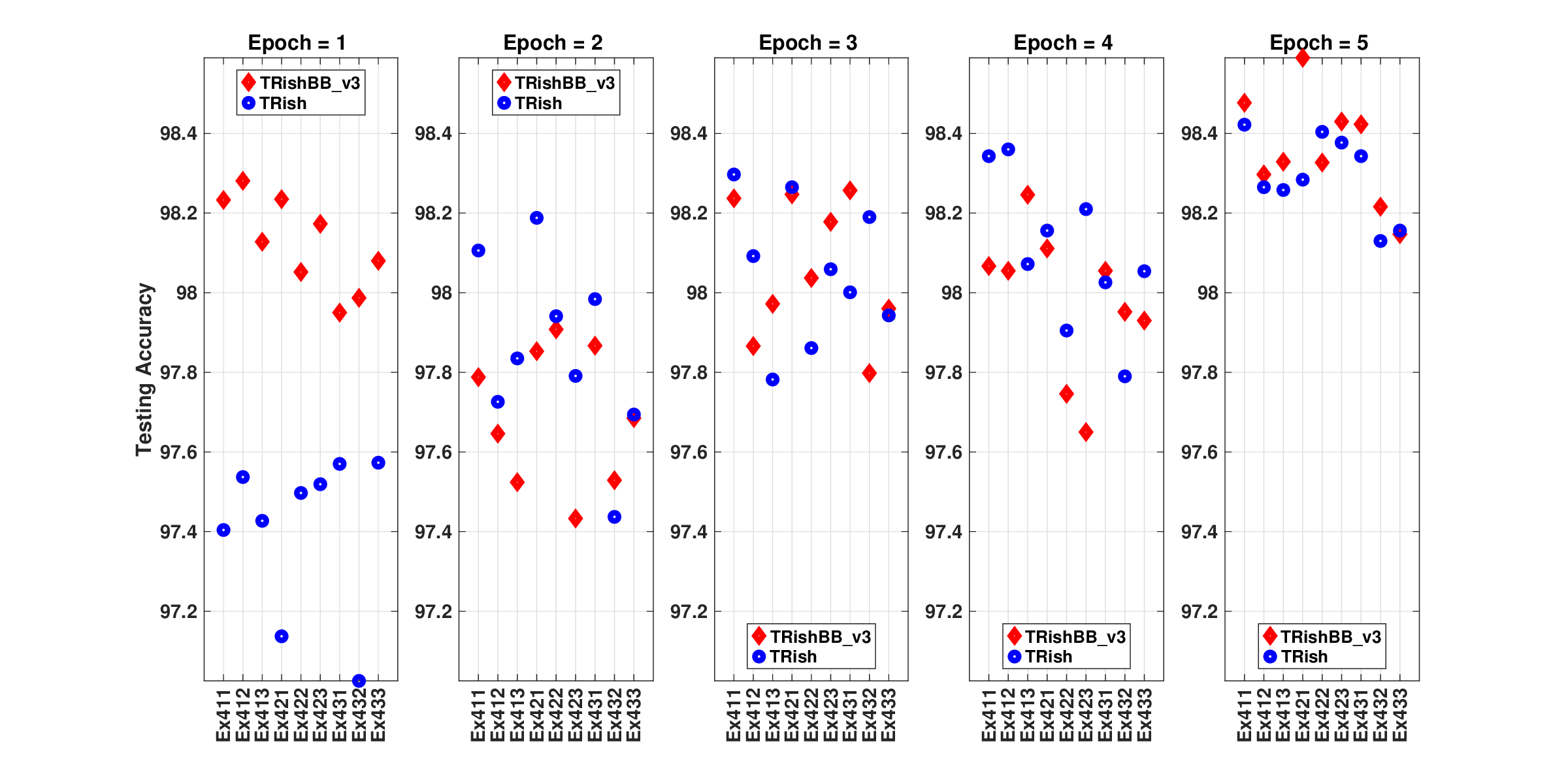}
    \end{adjustbox}    
    \caption{\footnotesize{ \texttt{MNIST}: Testing accuracy, $(\alpha, \gamma_1, \gamma_2)=(1, \gamma_1, \gamma_2)$, $|\mathcal{N}_k|=128$. Top: \trishunop, $m=20$. Middle: \trishduep, $m=468$. Bottom: \trishtrep, $m=468$ }}
    \label{fig:mnist}
\end{figure}

 \Cref{fig:mnist} displays the averaged accuracy obtained with the dataset \texttt{MNIST}. \trishdue compares well with   \trishuno even though it employed a smaller percentage  of unconstrained steps. Both \trishuno and \trishdue are considerably less sensitive to the values of $\gamma_1$ and $\gamma_2$ compared to TRish. \trishtre exhibits performance similar to that of TRish, which is not surprising, as it selects the unconstrained step in only approximately 15.23\% of the iterations. Nevertheless, \trishtre outperforms TRish at the first epoch due to the predominance of BB steplengths in that phase; specifically, Table~\ref{tab:sol1} shows the number of BB steplengths taken at the end of the first epoch compared to the total number taken at the end of the fifth epoch. 

\begin{figure}
\vspace{-50pt}
    \centering
    \begin{adjustbox}{width=1.0\linewidth, center}
    \includegraphics{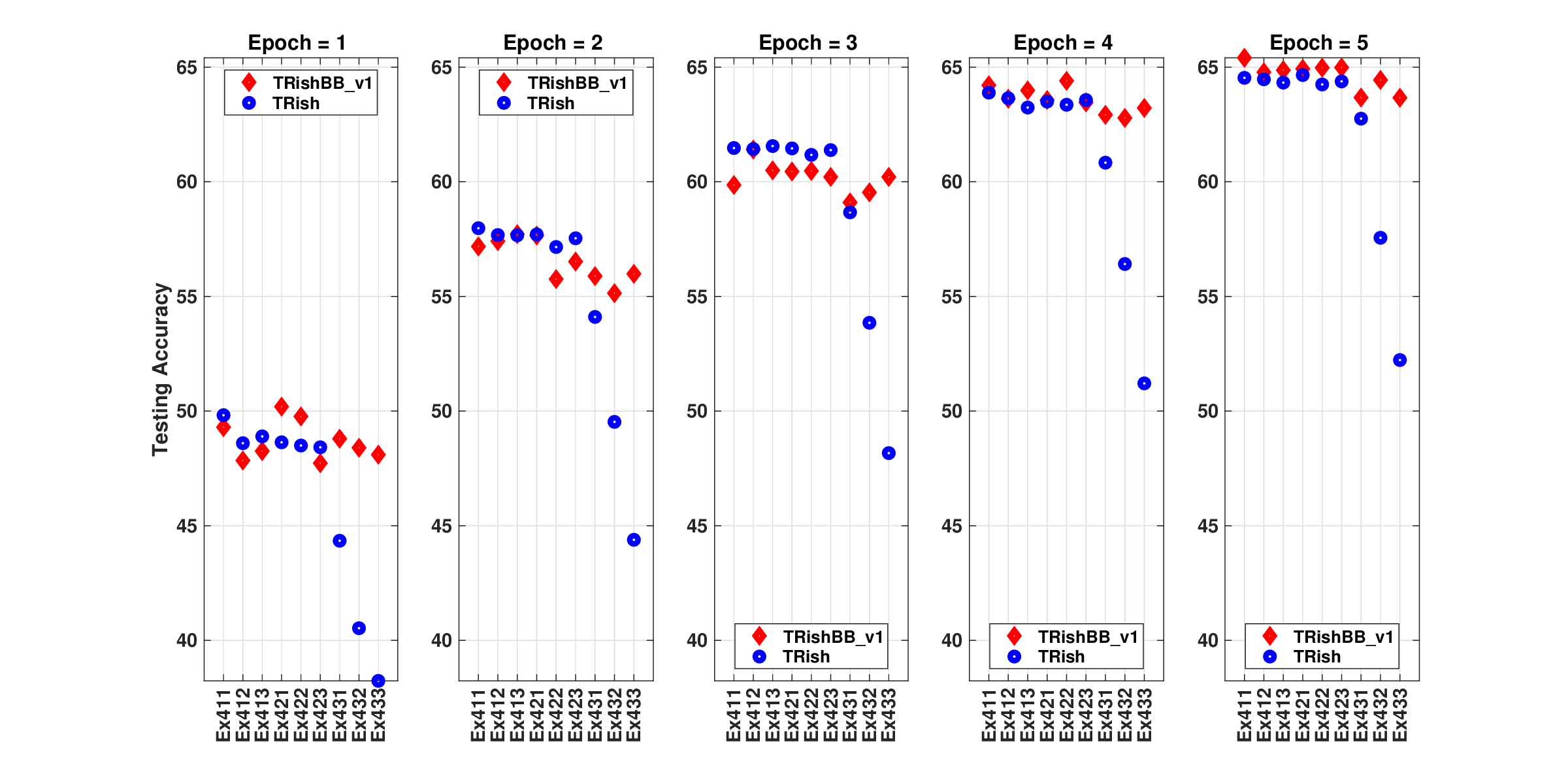}
    \end{adjustbox}    
    \centering
    \begin{adjustbox}{width=1.0\linewidth, center}
    \includegraphics{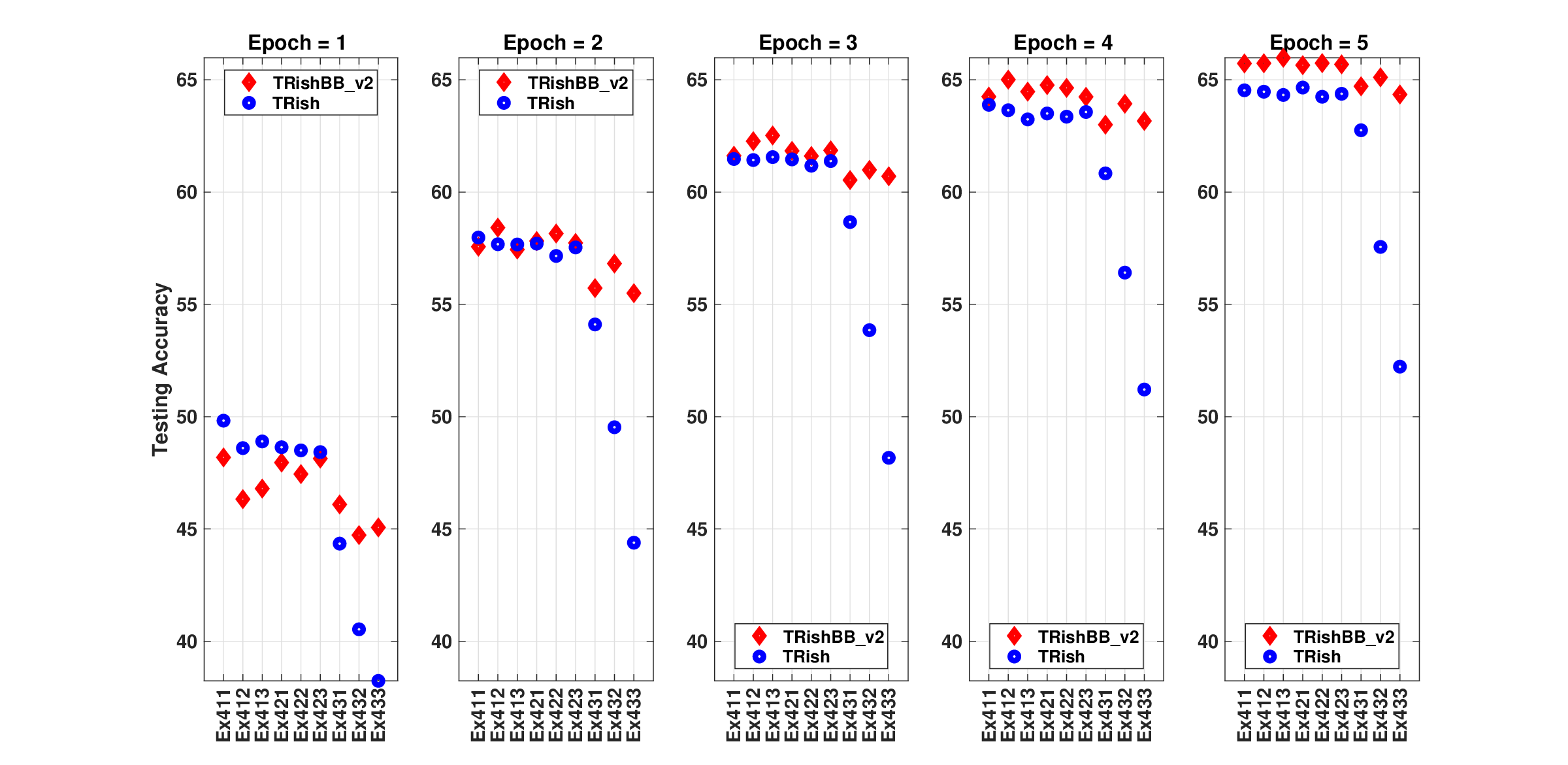}
    \end{adjustbox}    
    \centering
        \begin{adjustbox}{width=1.0\linewidth, center}
    \includegraphics{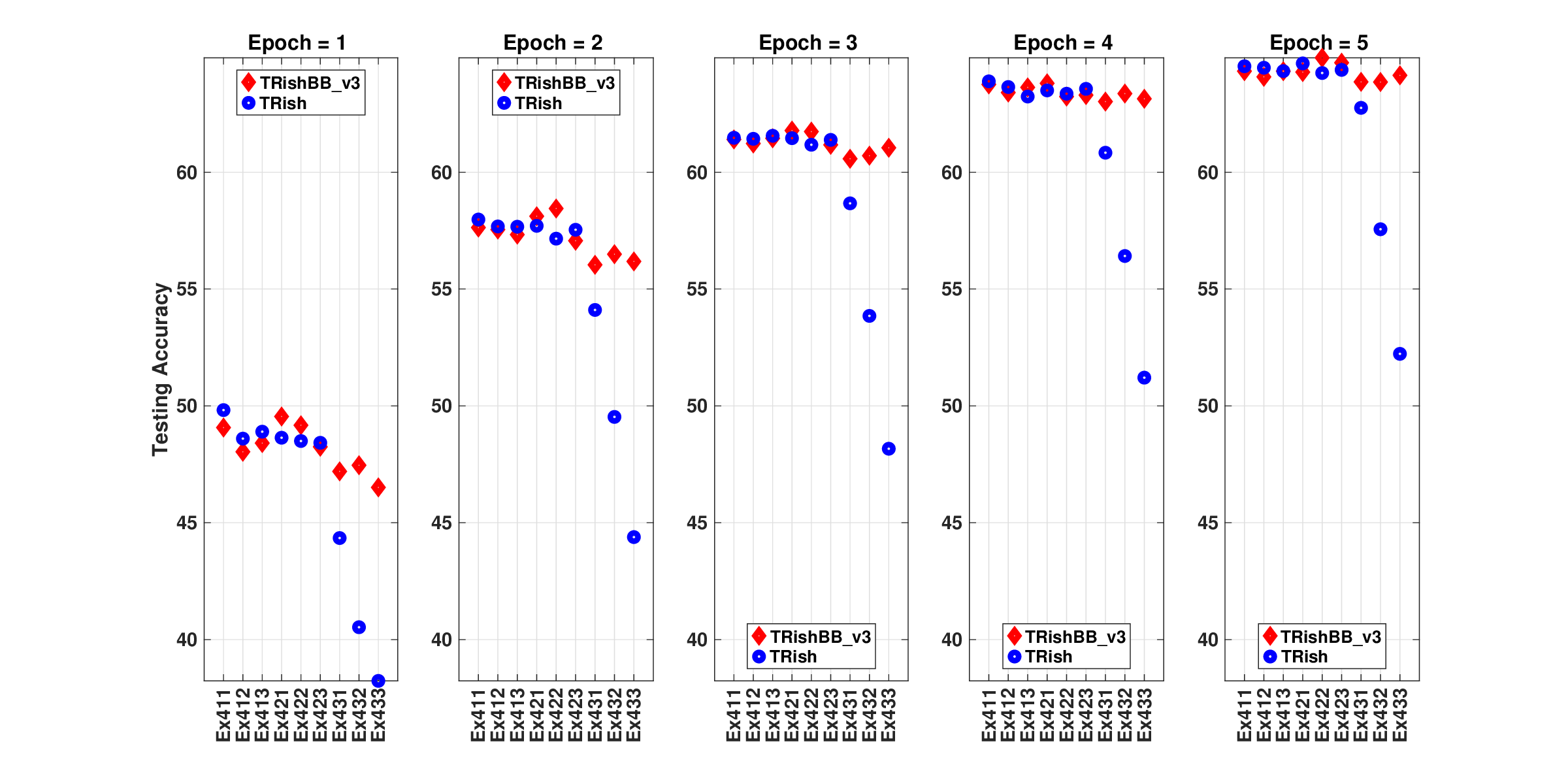}
    \end{adjustbox}    
    \caption{\footnotesize{ \texttt{CIFAR10}: Testing accuracy, $(\alpha, \gamma_1, \gamma_2)=(1, \gamma_1, \gamma_2)$, $|\mathcal{N}_k|=128$. Top: \trishunop, $m=20$. Middle: \trishduep, $m=390$. Bottom: \trishtrep, $m=390$. }}
    \label{fig:cifar}
\end{figure}

\Cref{fig:cifar} shows the averaged accuracies with the dataset \texttt{CIFAR10}. The main observation is that the proposed TRishBB variants  are less sensitive than TRish to the  choice of $\gamma_1$ and $\gamma_2$, see e.g., experiments \texttt{Ex431}, \texttt{Ex432} and \texttt{Ex433}.
Since accuracies of TRishBB variants and TRish are increasing during the epochs, we focus on the fifth epoch and observe that \trishdue performs better than all other algorithms. \trishuno and \trishtre performed better than TRish limited to some choices of $\gamma_1$ and $\gamma_2$.

\Cref{fig:muA_alpha5}  shows that the values of $\mu_k$ generated by \trishuno are  smaller than \trishdue and \trishtrep, consistently with Table \ref{tab:2}. Although \trishuno takes many BB steplengths, \trishdue stands out as the best-performing variant as in the previous set of experiments.

\begin{table}[t]
\centering
{\footnotesize
\caption{The number of BB steplengths $\mu_k$ taken by \trishtre on \texttt{MNIST} at the end of first epoch vs.  the total number taken at the end of the fifth epoch.}
\label{tab:sol1}
\begin{tabular}{|l|*{9}{c|}}
\hline
\textbf{Experiments} & \textbf{\texttt{E411}} & \textbf{\texttt{E412}} & \textbf{\texttt{E413}} & \textbf{\texttt{E421}} & \textbf{\texttt{E422}} & \textbf{\texttt{E423}} & \textbf{\texttt{E431}} & \textbf{\texttt{E432}} & \textbf{\texttt{E433}} \\
\hline
{End of $\text{Epoch} ~1$} &  359&  356&  357& 357 & 361 & 357 & 355 & 348 & 350 \\
\hline
{End of $\text{Epoch} ~5$} & 359 & 356 & 358 & 357 & 361 & 360 & 355 &  348& 357 \\
\hline
\end{tabular}
}
\end{table}


\subsubsection{Average testing loss obtained by \trishdue} 
\begin{figure}[t]
\centering
    \begin{subfigure}[t]{\textwidth}\centering
        {\includegraphics[width=0.32\linewidth]{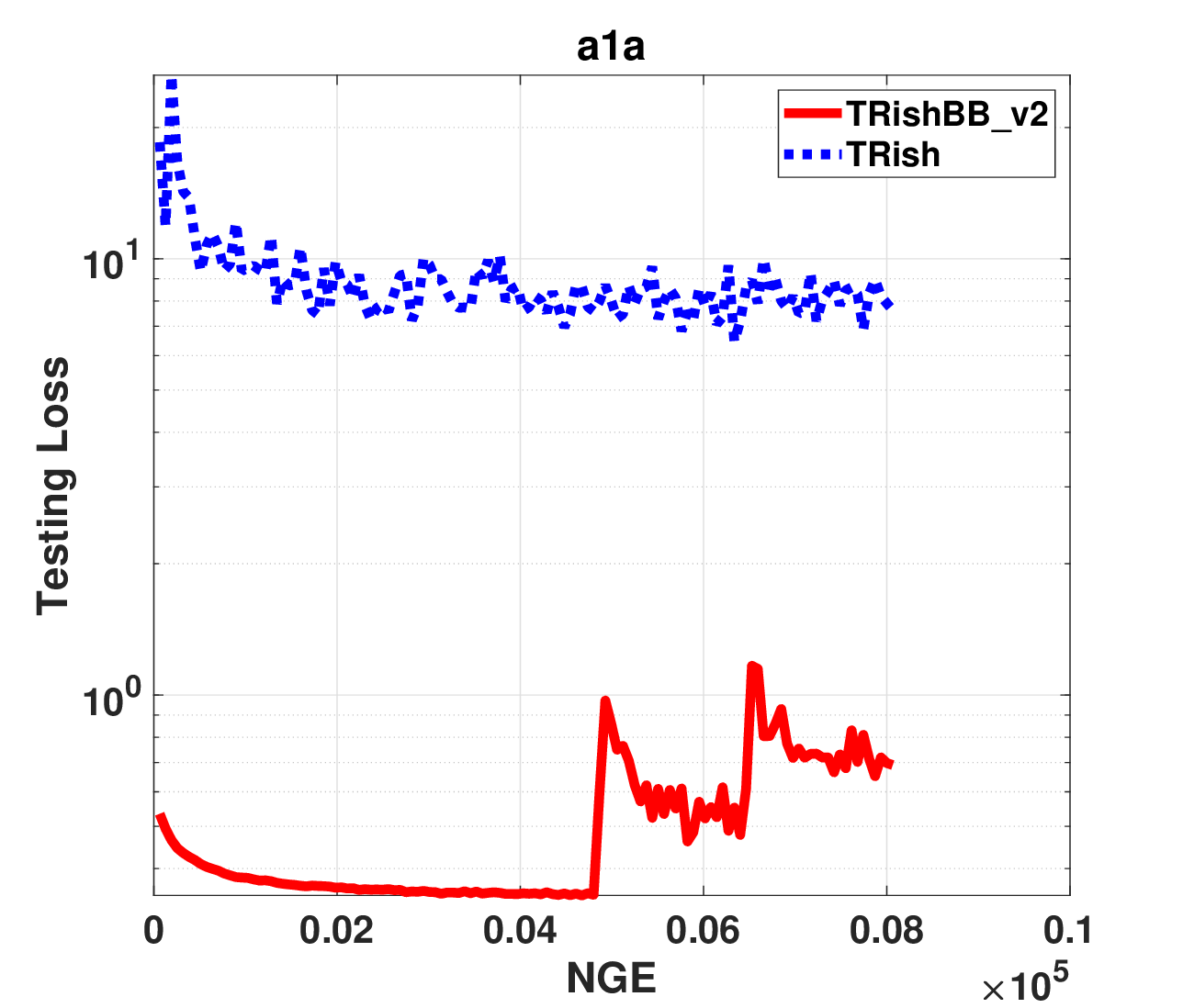}}
        {\includegraphics[width=0.32\linewidth]{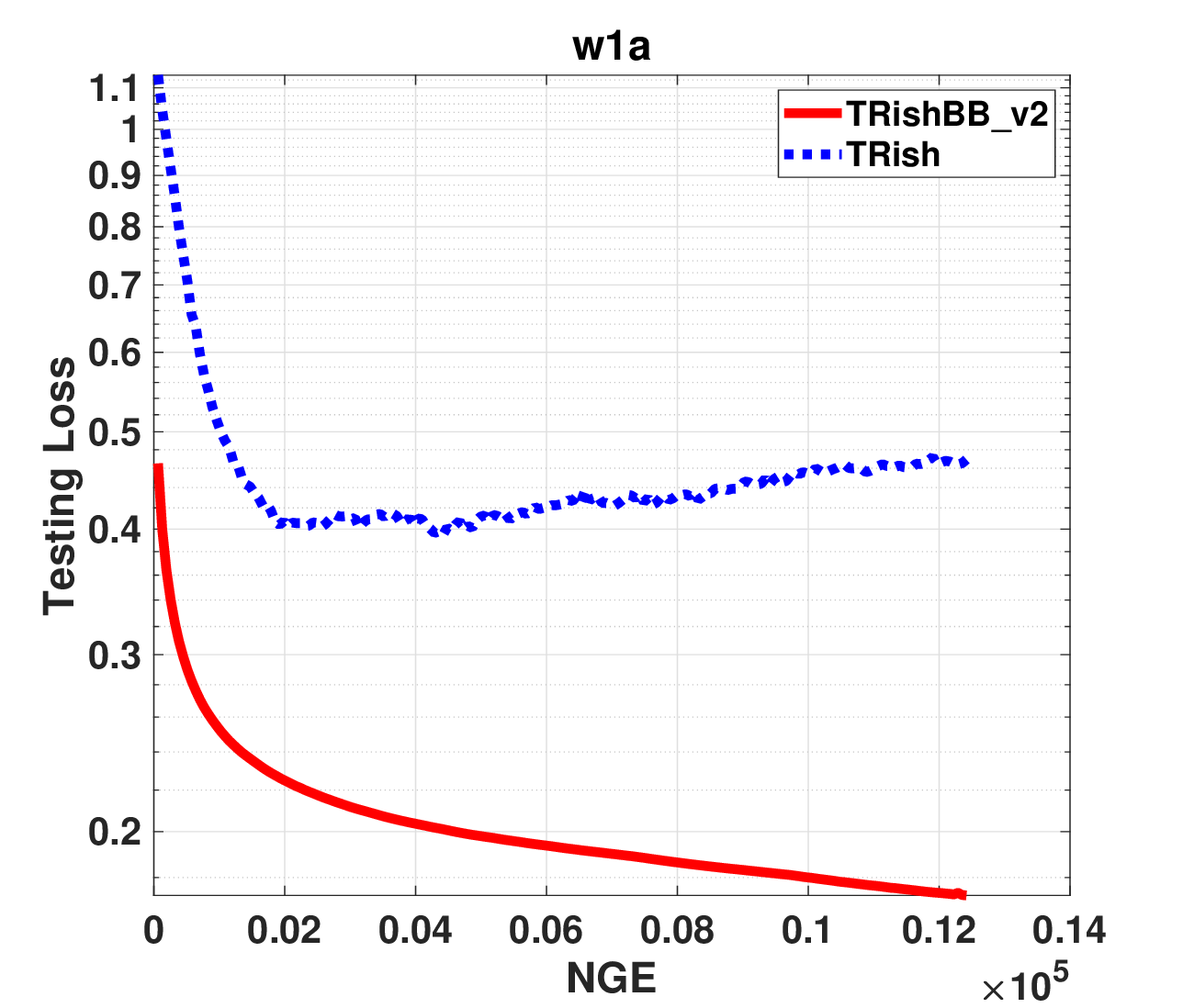}}
        {\includegraphics[width=0.32\linewidth]{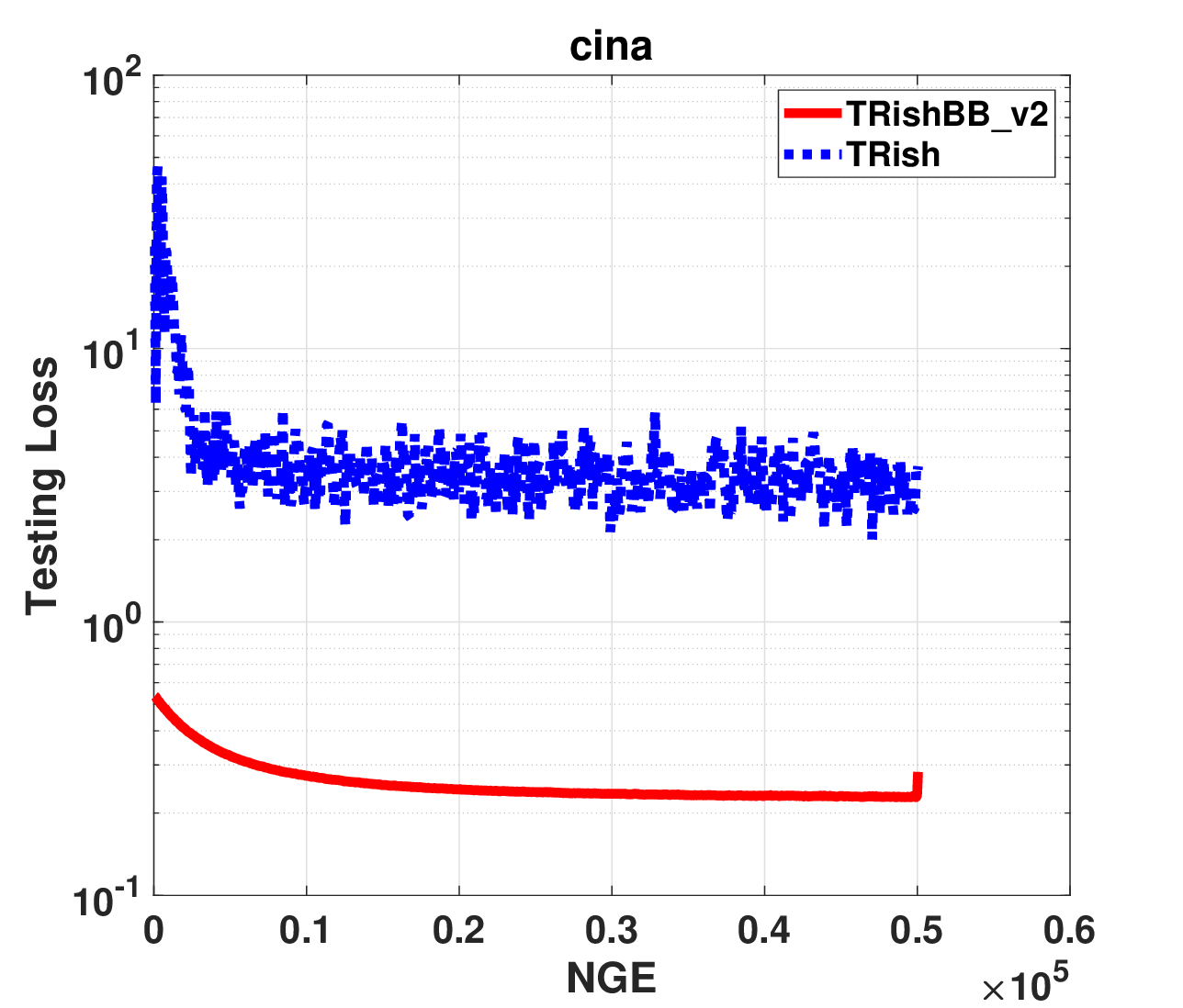}}
        {\includegraphics[width=0.32\linewidth]{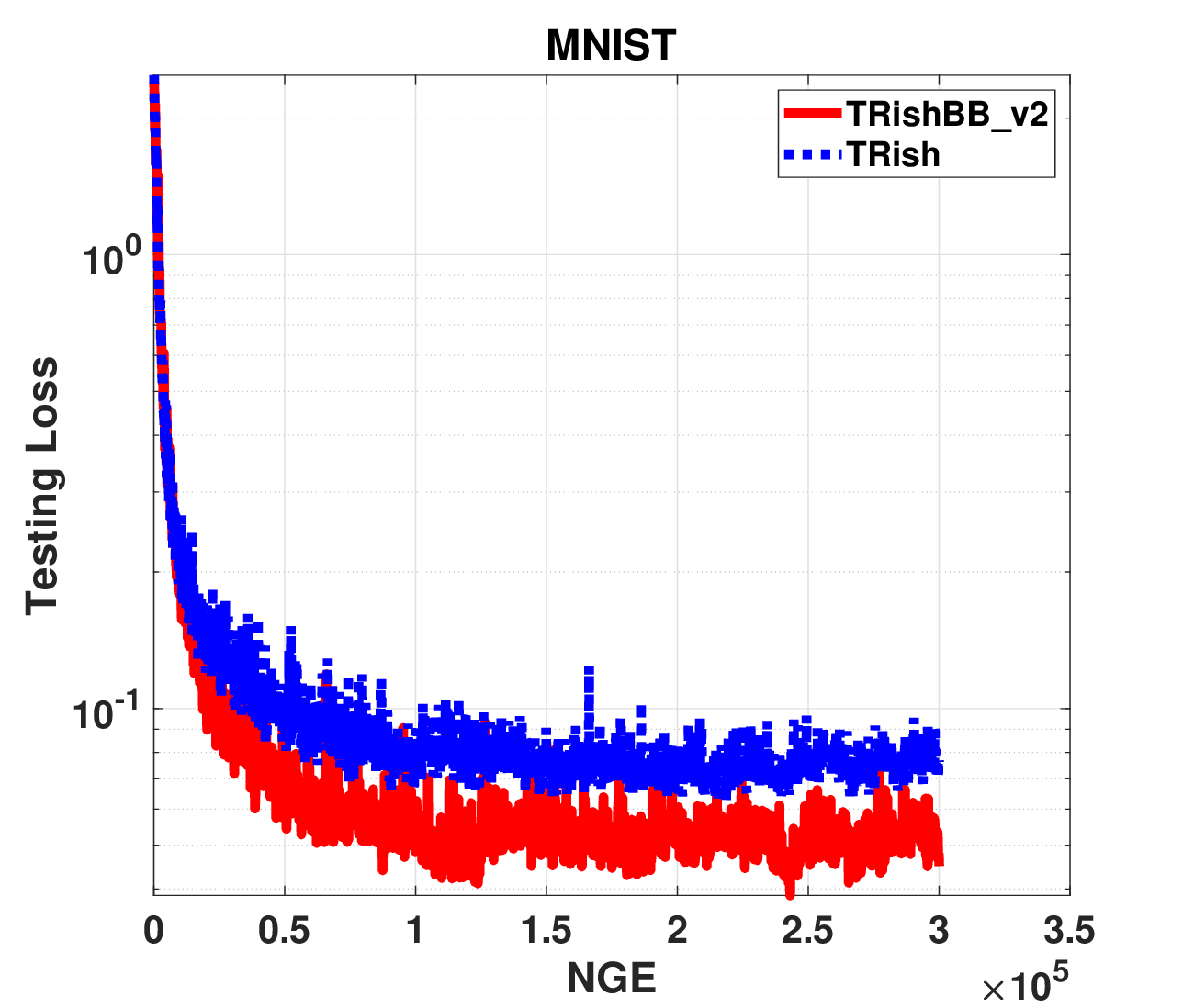}}
        {\includegraphics[width=0.32\linewidth]{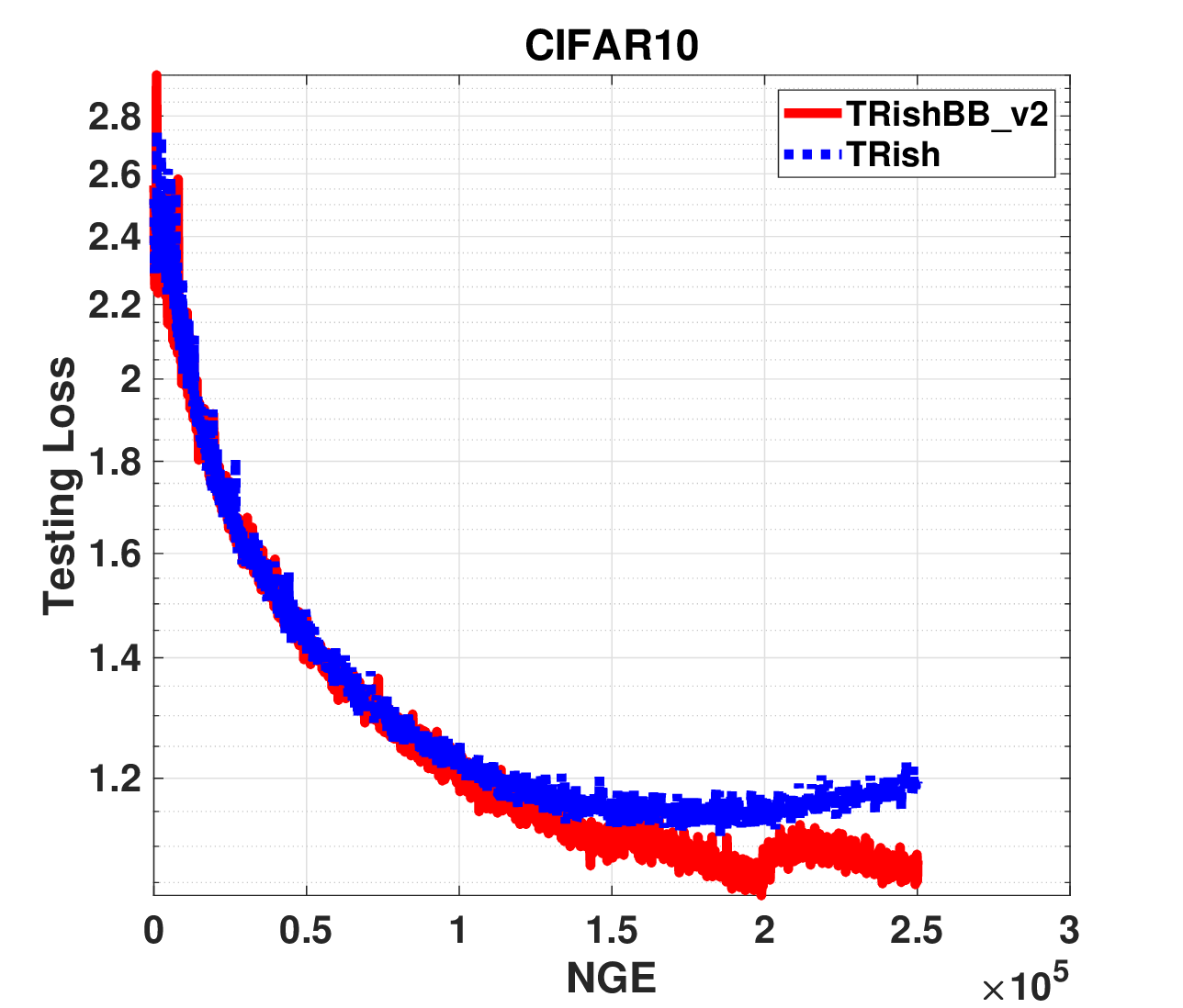}}
    \end{subfigure}
    \caption{\label{fig:avg_testing_loss} Average testing loss of \trishdue vs. number of gradient evaluations (NGE) compared to TRish across five datasets.}
\end{figure}

Sections~\ref{S.a1aw1acina} and~\ref{S.MnistCifar10} indicate that \trishdue is the overall best variant of TRishBB. In this section, we present the average testing loss vs. the number of gradient evaluations  obtained with both \trishdue and TRish over five epochs.
We consider specific configurations \texttt{ExIJK} for the triplet $(\alpha, \gamma_1, \gamma_2)$ described in Section~\ref{NumericalSec_configuration}. The results shown in Figure~\ref{fig:avg_testing_loss} 
refer to the highest average testing accuracy of these algorithms in Table~\ref{tab:accuracy} for the largest value of $\alpha$. Specifically, \trishdue achieves its best performance with configuration \texttt{Ex511} on \texttt{a1a}, \texttt{w1a}, and \texttt{cina}, \texttt{Ex422} on \texttt{MNIST}, and \texttt{Ex413} on \texttt{CIFAR10}, while TRish attains its highest accuracy with \texttt{Ex532} on \texttt{a1a}, \texttt{Ex524} on \texttt{w1a}, \texttt{Ex531} on \texttt{cina}, \texttt{Ex411} on \texttt{MNIST}, and \texttt{Ex421} on \texttt{CIFAR10}.
The figure shows that our method reaches lower testing loss values earlier than TRish and that in the solution of \texttt{a1a}, \texttt{w1a}, and \texttt{cina} 
the value of the testing loss at termination is largely different and in favour of \trishduep.

\subsubsection{Comparison to other algorithms}\label{S.ComOther}
As discussed in Sections \ref{Sub2.3} and \ref{Sub2.4},  \trishdue borrows features from  {SGD-BB} method given in \cite{tan2016barzilai}, while \trishtre borrows ideas from {AdaQN} Algorithm given in \cite{keskar2016adaqn}. In this section, we compare the performance of these algorithms. For \trishdue and \trishtrep, all hyper-parameters are set as described in Section~\ref{NumericalSec_configuration}.

\paragraph{\textbf{\trishdue vs. {SGD-BB}}}{

We implemented {SGD-BB} based on Algorithm~3 in~\cite{tan2016barzilai} allowing mini-batch size greater than one. We set the cardinality of $\mathcal{N}_k$, the scalar $\beta$ and the cyclic update parameter $m$ as described in Section \ref{NumericalSec_configuration} for the experiments with \trishduep.

We remark that the original SGD-BB algorithm did not apply thresholding for $\mu_k$ but our experiments demonstrated that thresholding is necessary for robustness, along with fine tuning for the initial value $\mu_0$. Specifically, initially, we ran {SGD-BB} with $\mu_0 = 1$ which was the setting used in \trishdue but the  BB steplength generated  in SGD-BB were too large causing failures in the fist epoch. To mitigate this issue, we ran the algorithm with $\mu_0 = 10^{-2}$ for \texttt{a1a}, \texttt{w1a}, and \texttt{cina}, and with $\mu_0 = 10^{-3}$ for \texttt{MNIST} and \texttt{CIFAR10}. Nonetheless, these smaller values for $\mu_0$ were still insufficient to provide good results in the next epochs. Hence, we applied the thresholding of $\mu_k$ as in \trishdue and verified that setting $\mu_{\min}=10^{-5}$ and $\mu_{\max} = 10^{-1}$ provided convergence for all datasets.

Table \ref{tab:accuracy2} compares the performance of {SGD-BB} against \trishdue performed with the largest values of $\alpha$ indicated in Sections~\ref{S.a1aw1acina} and~\ref{S.MnistCifar10}. The highest average testing accuracy achieved by \trishdue per epoch is close to or better than the accuracy achieved by tuned SGD-BB. Indeed, in SGD-BB a careful adjustment of both the initial stepsize $\mu_0$ and the threshold $\mu_{\max}$ was necessary to obtain good performance and the average testing accuracy may differ greatly from one epoch to the other.

\begin{table}[t!] 
\centering
\caption{Highest average testing accuracy per epoch obtained by \trishdue with the largest $\alpha$ and by SGD-BB.}
\label{tab:accuracy2}
\resizebox{\textwidth}{!}{
{\small 
\begin{tabular}{
|p{5.8cm}|>
{\centering\arraybackslash}m{1.8cm}|>
{\centering\arraybackslash}m{1.8cm}|>
{\centering\arraybackslash}m{1.8cm}|>
{\centering\arraybackslash}m{1.8cm}|>
{\centering\arraybackslash}m{1.8cm}|}
\hline
\texttt{ala} & 
$\text{epoch}=1$ & 
$\text{epoch}=2$ & 
$\text{epoch}=3$ & 
$\text{epoch}=4$ & 
$\text{epoch}=5$ \\
\hline
{\trishdue}  &  83.14 &    83.52 &    79.71 &    80.94 &    80.77  \\
\hline
{SGD-BB}     &  82.63 &    81.73 &    60.03 &    61.76 &    54.57  \\ 
\hline
\end{tabular}
}}
\medskip

\resizebox{\textwidth}{!}{
{\small 
\begin{tabular}{
|p{5.8cm}|>
{\centering\arraybackslash}m{1.8cm}|>
{\centering\arraybackslash}m{1.8cm}|>
{\centering\arraybackslash}m{1.8cm}|>
{\centering\arraybackslash}m{1.8cm}|>
{\centering\arraybackslash}m{1.8cm}|}
\hline

\texttt{wla} & 
$\text{epoch}=1$ & 
$\text{epoch}=2$ & 
$\text{epoch}=3$ & 
$\text{epoch}=4$ & 
$\text{epoch}=5$ \\
\hline
{\trishdue}  & 89.38 &    89.44 &    89.47 &    89.53 &    89.57  \\
\hline
{SGD-BB}     &  89.38 &    89.46 &    89.54 &   89.45 &    87.74  \\
\hline
\end{tabular}
}}
\medskip

\resizebox{\textwidth}{!}{
{\small 
\begin{tabular}{
|p{5.8cm}|>
{\centering\arraybackslash}m{1.8cm}|>
{\centering\arraybackslash}m{1.8cm}|>
{\centering\arraybackslash}m{1.8cm}|>
{\centering\arraybackslash}m{1.8cm}|>
{\centering\arraybackslash}m{1.8cm}|}
\hline
\texttt{cina} & 
$\text{epoch}=1$ & 
$\text{epoch}=2$ & 
$\text{epoch}=3$ & 
$\text{epoch}=4$ & 
$\text{epoch}=5$ \\
\hline
{\trishdue}  & 90.07 &    90.99 &    91.35 &    91.45 &    91.20  \\
\hline
{SGD-BB}     & 89.52 &    91.01 &    91.32 &    83.99 &    49.21 \\
\hline
\end{tabular}
}}
\medskip

\resizebox{\textwidth}{!}{
{\small 
\begin{tabular}{
|p{5.8cm}|>
{\centering\arraybackslash}m{1.8cm}|>
{\centering\arraybackslash}m{1.8cm}|>
{\centering\arraybackslash}m{1.8cm}|>
{\centering\arraybackslash}m{1.8cm}|>
{\centering\arraybackslash}m{1.8cm}|}
\hline
\texttt{MNIST} & 
$\text{epoch}=1$ & 
$\text{epoch}=2$ & 
$\text{epoch}=3$ & 
$\text{epoch}=4$ & 
$\text{epoch}=5$ \\
\hline
{\trishdue}  & 98.30 &    98.61 &    98.73 &    98.71 &    98.81  \\
\hline
{SGD-BB}     &  93.44 &    95.01 &    73.08 &    98.26 &    98.70  \\
\hline
\end{tabular}
}}
\medskip

\resizebox{\textwidth}{!}{
{\small 
\begin{tabular}{
|p{5.8cm}|>
{\centering\arraybackslash}m{1.8cm}|>
{\centering\arraybackslash}m{1.8cm}|>
{\centering\arraybackslash}m{1.8cm}|>
{\centering\arraybackslash}m{1.8cm}|>
{\centering\arraybackslash}m{1.8cm}|}
\hline
\texttt{CIFAR10} & 
$\text{epoch}=1$ & 
$\text{epoch}=2$ & 
$\text{epoch}=3$ & 
$\text{epoch}=4$ & 
$\text{epoch}=5$ \\
\hline
{\trishdue}  & 48.18 &    58.42 &    62.52 &    65.01 &    65.98  \\
\hline
{SGD-BB}     &  20.54 &    27.49 &    19.57 &    51.67 &    57.83 \\
\hline
\end{tabular}
}}
\end{table}

}

\paragraph{\textbf{\trishtre vs. {AdaQN}}}{

Our experiments with {AdaQN} were performed using  
\texttt{adaQN.m} from the package \texttt{minSQN}, available at \url{https://github.com/keskarnitish/minSQN/blob/master/%2Bmethods/adaQN.m}. As suggested therein and in the reference paper \cite{keskar2016adaqn}, we set the hyper-parameters

We set the limited-memory size of the L-BFGS Hessian approximations $H_k$ to $20$, and the Fisher information storage size $m_F$ to $100$. AdaQN provides two strategies for initializing $H_0$. We adopted the default option,  RMSprop-like initialization, which showed better performance in our experiments compared to the Adagrad-like strategy.

The cyclic iteration parameter $L$ in {AdaQN}, corresponding to our parameter $m$ in  \trishtrep, needs fine tuning and we observed that it is highly sensitive to the choice of random seed during hyper-parameter tuning. We fixed its value in the set $L\in  \{2, 5, 10, 20, 64, m\}$ where the parameter $m$ is the same as in \trishtrep. The remaining hyper-parameter values were adopted from~\cite{keskar2016adaqn} and 
the implementation mentioned above. Following this implementation, the stepsize $\alpha$ in {AdaQN} was sampled from a log-uniform distribution over the interval $[10^{-6}, 10^2]$. We tuned {AdaQN} using 20 different random seeds for the values of \(L\) and \(\alpha\), as part of the standard procedure in the \texttt{adaQN.m} routine, with each run lasting 5 epochs. Using the resulting optimal hyper-parameters \(L\) and \(\alpha\), we obtained the results presented in Table~\ref{tab:accuracy3} where {AdaQN} is compared with \trishtre using the largest $\alpha$ indicated in Sections~\ref{S.a1aw1acina} and~\ref{S.MnistCifar10}. 

\begin{table}[t!] 
\centering
\caption{Highest average testing accuracy per epoch obtained by  \trishtre  with  the largest $\alpha$ and by {AdaQN} with its tuned $(L, \alpha)$.}
\label{tab:accuracy3}
\resizebox{\textwidth}{!}{
{\small 
\begin{tabular}{
|p{5.8cm}|>
{\centering\arraybackslash}m{1.8cm}|>
{\centering\arraybackslash}m{1.8cm}|>
{\centering\arraybackslash}m{1.8cm}|>
{\centering\arraybackslash}m{1.8cm}|>
{\centering\arraybackslash}m{1.8cm}|}
\hline
\texttt{ala} & 
$\text{epoch}=1$ & 
$\text{epoch}=2$ & 
$\text{epoch}=3$ & 
$\text{epoch}=4$ & 
$\text{epoch}=5$ \\
\hline
{\trishtre}  & 83.19 &    83.40 &    82.60 &    78.65 &    79.42 \\
\hline
{AdaQN} $(L=10,\, \alpha=1.13\times 10^{-3})$    &  76.62 &    78.83 &    80.34 &    80.52 &    80.48  \\
\hline
\end{tabular}
}}
\medskip

\resizebox{\textwidth}{!}{
{\small 
\begin{tabular}{
|p{5.8cm}|>
{\centering\arraybackslash}m{1.8cm}|>
{\centering\arraybackslash}m{1.8cm}|>
{\centering\arraybackslash}m{1.8cm}|>
{\centering\arraybackslash}m{1.8cm}|>
{\centering\arraybackslash}m{1.8cm}|}
\hline

\texttt{wla} & 
$\text{epoch}=1$ & 
$\text{epoch}=2$ & 
$\text{epoch}=3$ & 
$\text{epoch}=4$ & 
$\text{epoch}=5$  \\
\hline
{\trishtre}  & 89.34 &    89.42 &    89.51 &    89.47 &    89.24  \\
\hline
{AdaQN}  ($L=20, \, \alpha =8.41 \times 10^{-3} $)   &  48.27 &    47.78 &    47.65 &    47.70 &    47.57  \\
\hline
\end{tabular}
}}
\medskip

\resizebox{\textwidth}{!}{
{\small 
\begin{tabular}{
|p{5.8cm}|>
{\centering\arraybackslash}m{1.8cm}|>
{\centering\arraybackslash}m{1.8cm}|>
{\centering\arraybackslash}m{1.8cm}|>
{\centering\arraybackslash}m{1.8cm}|>
{\centering\arraybackslash}m{1.8cm}|}
\hline
\texttt{cina} & 
$\text{epoch}=1$ & 
$\text{epoch}=2$ & 
$\text{epoch}=3$ & 
$\text{epoch}=4$ & 
$\text{epoch}=5$ \\
\hline
{\trishtre}  & 90.00 &   91.18 &   90.07 &   88.58 &   85.05 \\
\hline
{AdaQN} ($L=5, \, \alpha =7.21\times 10^{-3}$)    &  88.73 &    87.76 &    87.92 &    88.25 &    88.14  \\
\hline
\end{tabular}
}}
\medskip

\resizebox{\textwidth}{!}{
{\small 
\begin{tabular}{
|p{5.8cm}|>
{\centering\arraybackslash}m{1.8cm}|>
{\centering\arraybackslash}m{1.8cm}|>
{\centering\arraybackslash}m{1.8cm}|>
{\centering\arraybackslash}m{1.8cm}|>
{\centering\arraybackslash}m{1.8cm}|}
\hline
\texttt{MNIST} & 
$\text{epoch}=1$ & 
$\text{epoch}=2$ & 
$\text{epoch}=3$ & 
$\text{epoch}=4$ & 
$\text{epoch}=5$ \\
\hline
{\trishtre}   & 98.28 &   97.91 &   98.26 &   98.25 &   98.59  \\
\hline
{AdaQN} ($L=468, \, \alpha = 6.50\times10^{-4}$)   &  98.13 &    98.34 &    98.66 &    98.82 &    98.86  \\
\hline
\end{tabular}
}}
\medskip

\resizebox{\textwidth}{!}{
{\small 
\begin{tabular}{
|p{5.8cm}|>
{\centering\arraybackslash}m{1.8cm}|>
{\centering\arraybackslash}m{1.8cm}|>
{\centering\arraybackslash}m{1.8cm}|>
{\centering\arraybackslash}m{1.8cm}|>
{\centering\arraybackslash}m{1.8cm}|}
\hline
\texttt{CIFAR10} & 
$\text{epoch}=1$ & 
$\text{epoch}=2$ & 
$\text{epoch}=3$ & 
$\text{epoch}=4$ & 
$\text{epoch}=5$ \\
\hline
{\trishtre}  & 49.55 &    58.45 &   61.78 &    63.80 &   64.89  \\
\hline
{AdaQN}  ($L=64, \, \alpha = 4.73\times 10^{-4}$)    &   54.65 &    60.65 &    62.99 &    66.02 &    66.46  \\
\hline
\end{tabular}
}}
\end{table}

The table shows that the highest average testing accuracy achieved by \trishtre per epoch over \texttt{a1a}, \texttt{w1a}, and \texttt{cina} is generally larger than the accuracy achieved by AdaQN, in particular for earlier epochs.
On the other hand, the performance of both \trishtre and AdaQN is close to each other for \texttt{MNIST} while AdaQN outperforms \trishtre in the solution of  \texttt{CIFAR10}. In \texttt{MNIST} and  \texttt{CIFAR10} problems, \trishtre  could meet rarely the iterations at which the BB steplength is computed, see  Table~\ref{tab:2}; in addition, the good performance of {AdaQN} is motivated by the use of an L-BFGS-based Hessian approximation  coupled with a step acceptance criterion based on sampled function evaluations.

\subsubsection{ The cyclic iteration
parameter $m$ in \trishtre}\label{S.choiceM}

We conclude our experimental analysis, focusing on the cyclic iteration parameter $m$ in \trishtrep. In all experiments conducted we have fixed $m = N_b$ with $N_b = \lfloor N/|\mathcal{N}_k| \rfloor$ ($N_b$ equal to 25, 38, 156, 468, and 390 for \texttt{a1a}, \texttt{w1a}, \texttt{cina}, \texttt{MNIST}, and \texttt{CIFAR10}, respectively). To assess the sensitivity of \trishtre to the hyper-parameter $m$, we  report in  Figures~\ref{fig:ML(m)}--\ref{fig:DL(m)} additional results obtained with $m \in \{ 5,   20,64, N_b\}$, i.e., the values used for tuning Ada-QN; to improve readability in the figures we have eliminated the results corresponding to $m=2$ and $m=10$ because $m=5$ sufficiently represents the impact of small values of $m$. These figures show that for large values of $m$, $m=64$ and $m=N_b$, the performance is not highly sensitive to the choice of $\gamma_1$ and $\gamma_2$ and that
the setting $m=N_b$ provides high accuracy at some epochs for all tests, making unnecessary to tune a further hyper-parameter.

\section{Conclusion}\label{ConSec}
 We introduced the Algorithm TRishBB applicable to  finite-sum minimization problems. It is  a stochastic trust-region method that combines the TRish methodology with stochastic Barzilai-Borwein steplengths to inject cheap second-order information.  We presented theoretical analysis of TRishBB and devised three practical variants differing in the computation of the BB steplength. \trishuno mimics the implementation of a standard BB rule,  \trishdue employs an accumulated BB parameter per iteration, \trishtre retrieves the second-order information by the accumulated  Fisher Information matrix. 
Numerical results demonstrated that, for large values of $\alpha$, \trishdue  and \trishtre are quite insensitive to the hyper-parameters  and performs better than the first-order TRish method.   
We also  observe that they compare well with two stochastic  quasi-Newton methods from the literature. Taking into account that \trishtre is applicable as long as it is feasible the use of the empirical Fischer Information,  \trishdue appears to be a robust and general algorithm.

Regarding future investigations, it is worthy to further  explore the definition of stochastic BB parameters,  including stochastic diagonal BB rules, as suggested by \cite{franchini2023diagonal}.

\begin{figure}[H]
\vspace{-50pt}
    \centering
    \begin{adjustbox}{width=1.0\linewidth, center}
    \includegraphics{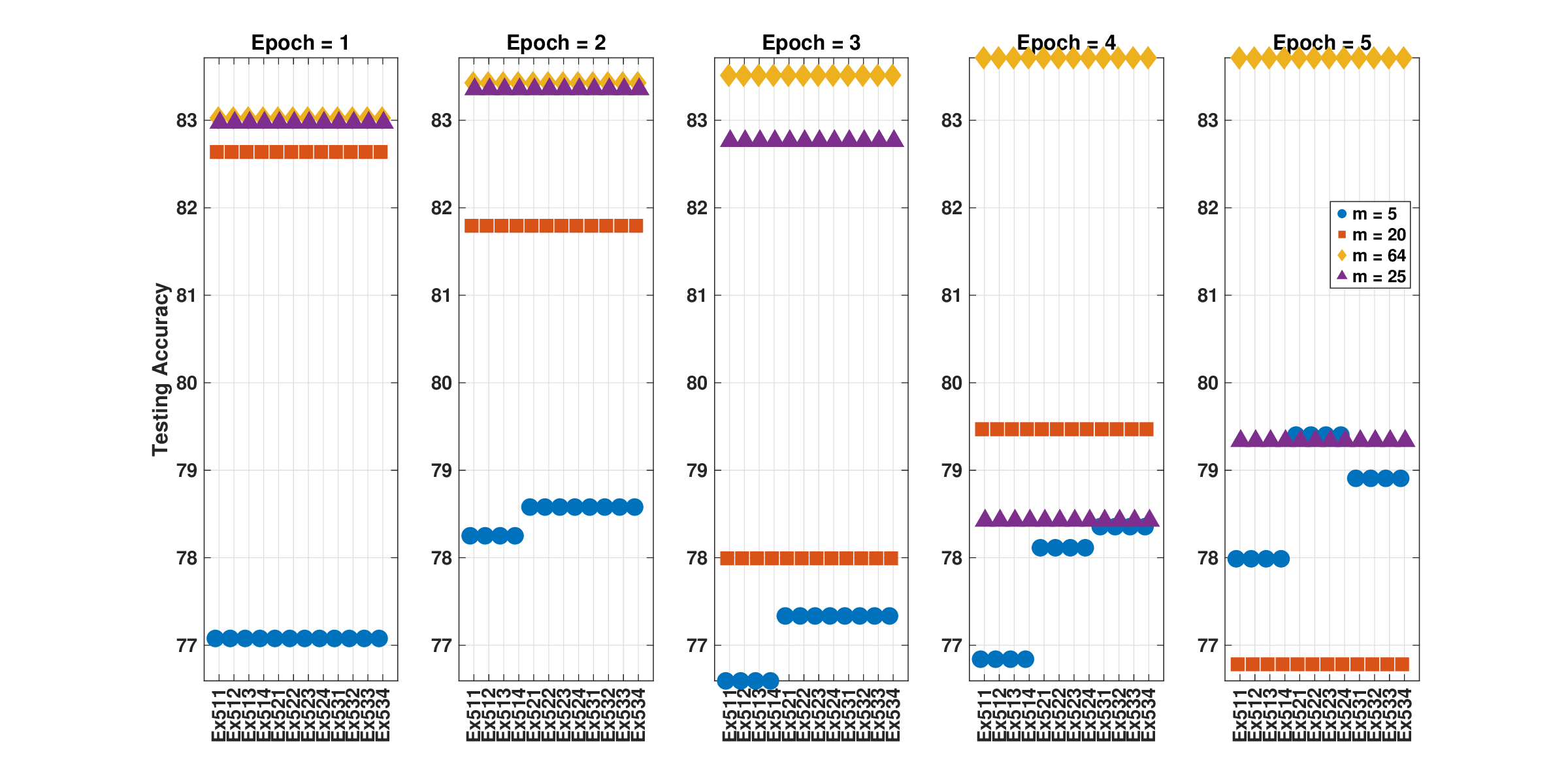}
    \end{adjustbox}    
  
    \centering
    \begin{adjustbox}{width=1.0\linewidth, center}
    \includegraphics{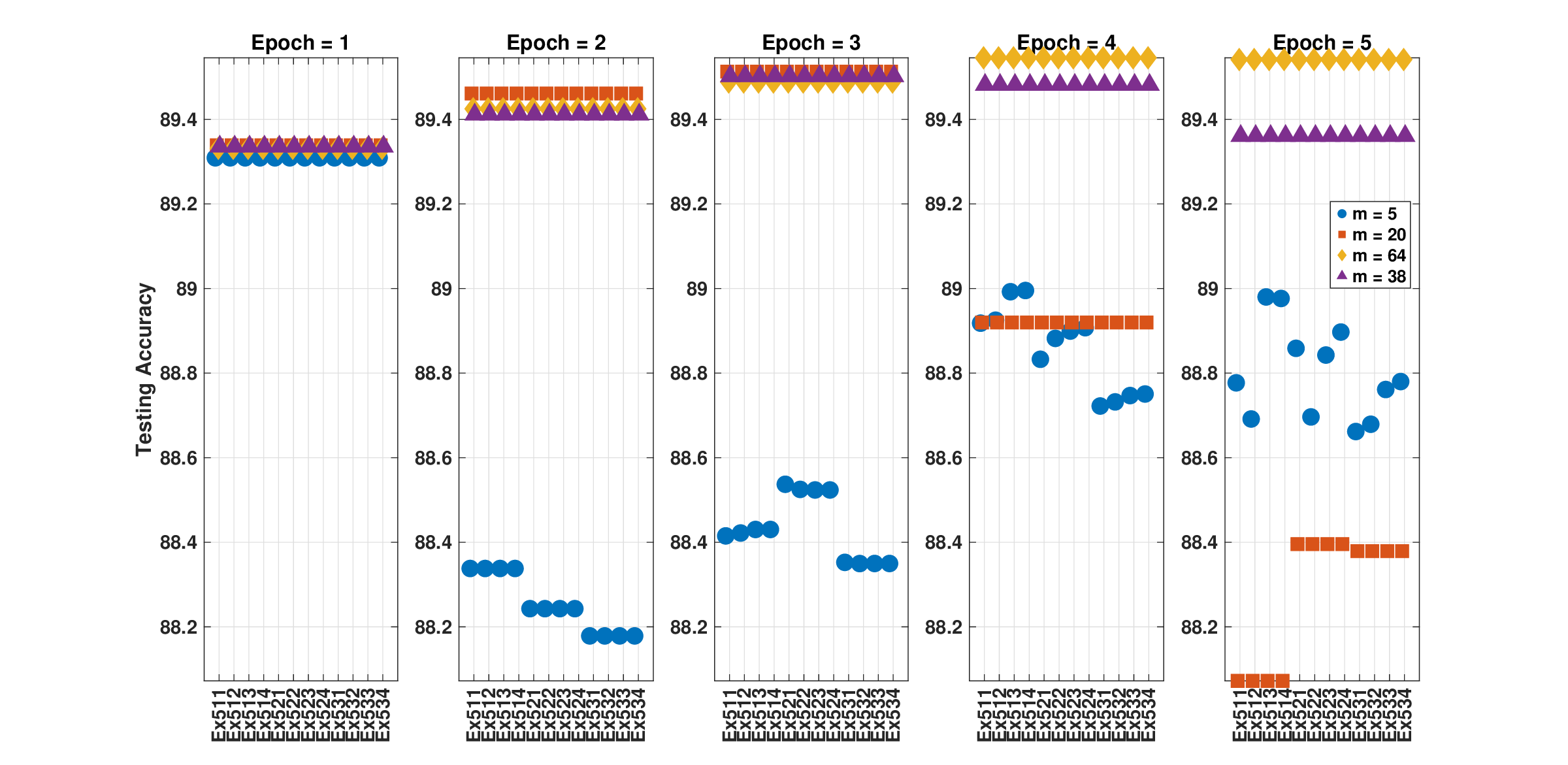}
    \end{adjustbox}    

    \centering
    \begin{adjustbox}{width=1.0\linewidth, center}
    \includegraphics{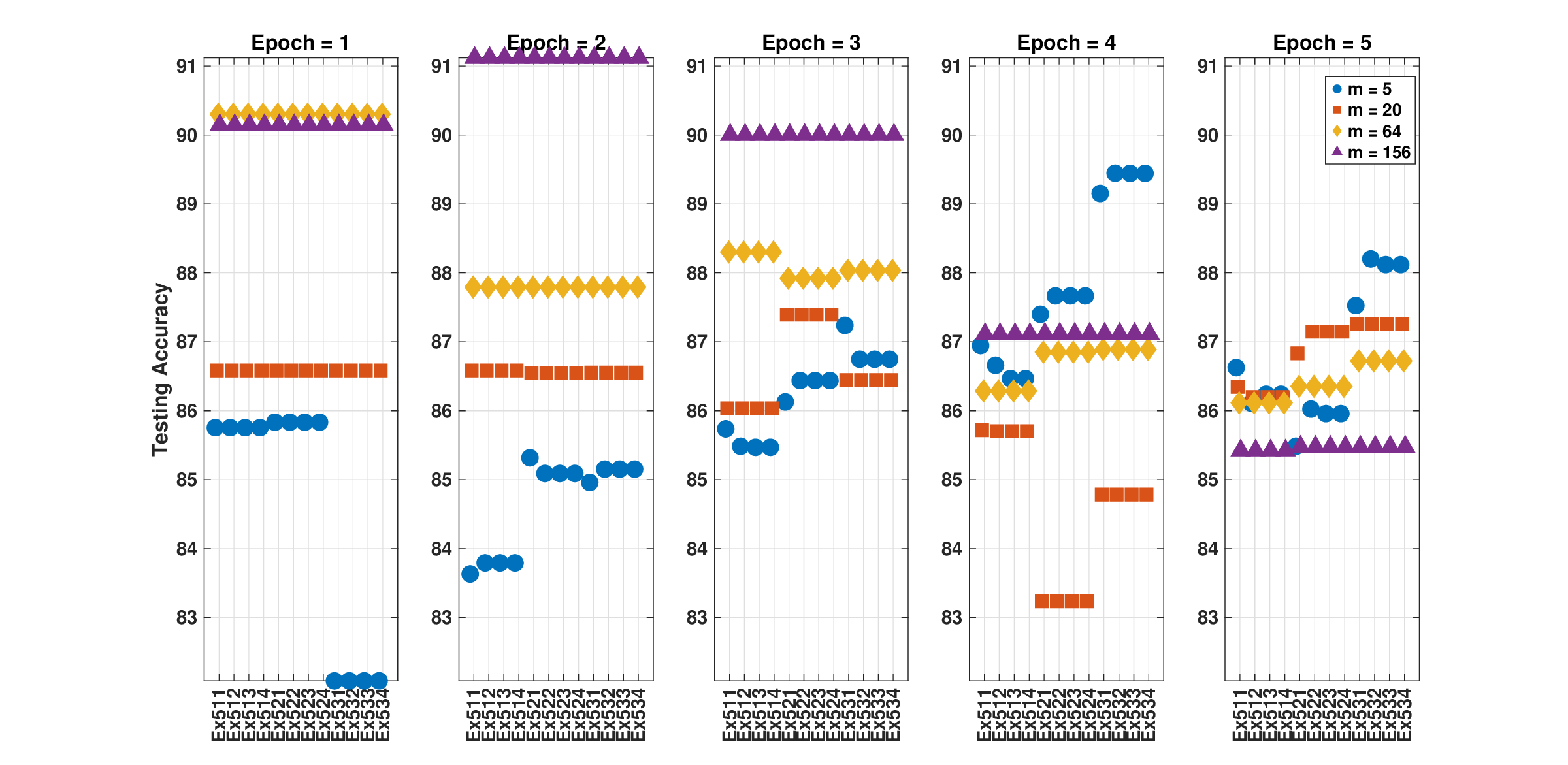}
    \end{adjustbox}    
    \caption{\footnotesize{Testing average accuracy of \trishtre with  $(\alpha, \gamma_1, \gamma_2)=(10, \gamma_1, \gamma_2)$, $m \in\{5,20,64,N_b\}$. 
    Top: \texttt{a1a}, $N_b=25$. Middle: \texttt{w1a}, $N_b=38$. Bottom: \texttt{cina}, $N_b=156$.}}
    \label{fig:ML(m)}
\end{figure}

\begin{figure}[H]
\vspace{-50pt}
    \centering
    \begin{adjustbox}{width=1.0\linewidth, center}
    \includegraphics{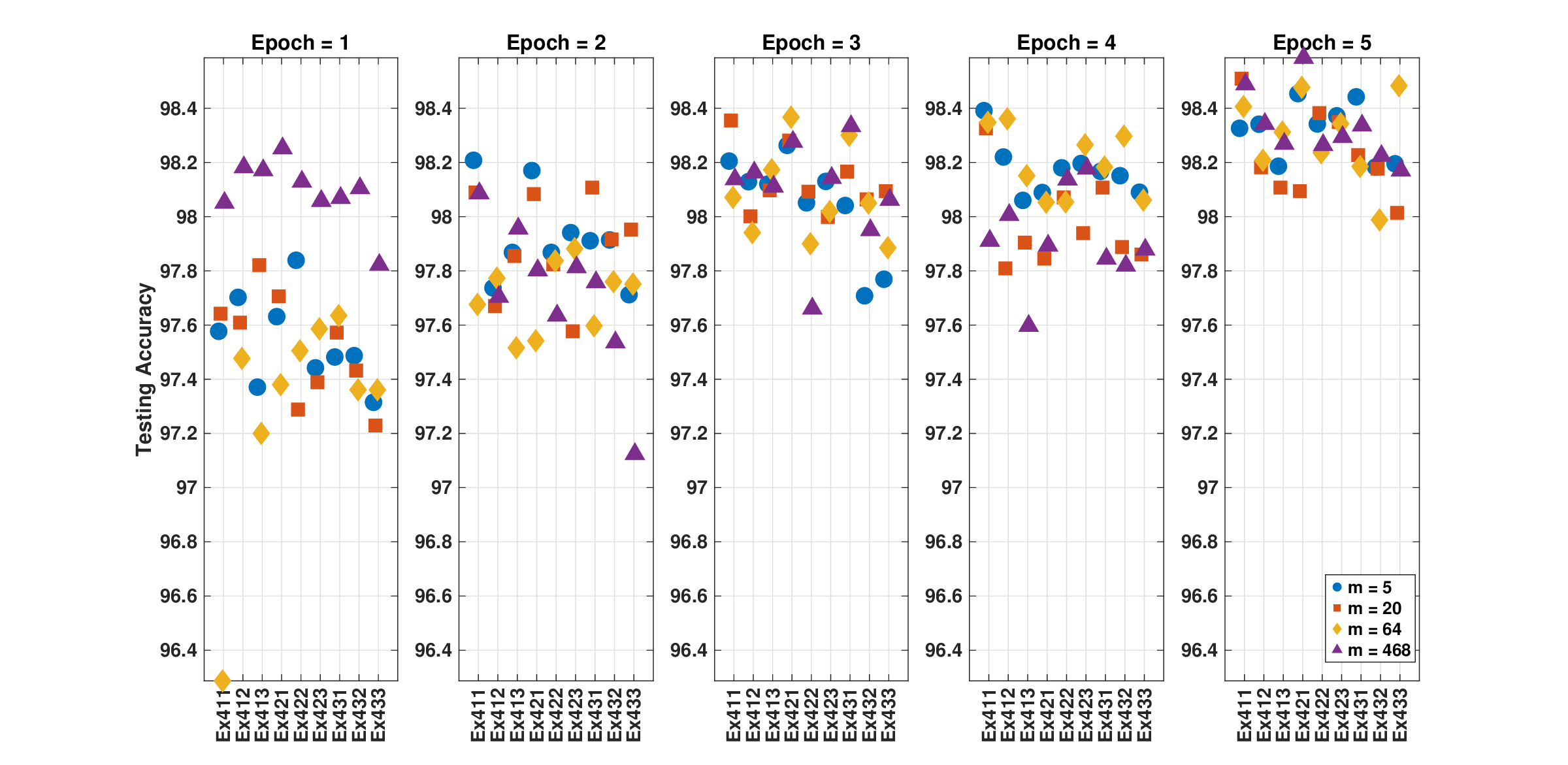}
    \end{adjustbox}    
    \centering
    \begin{adjustbox}{width=1.0\linewidth, center}
    \includegraphics{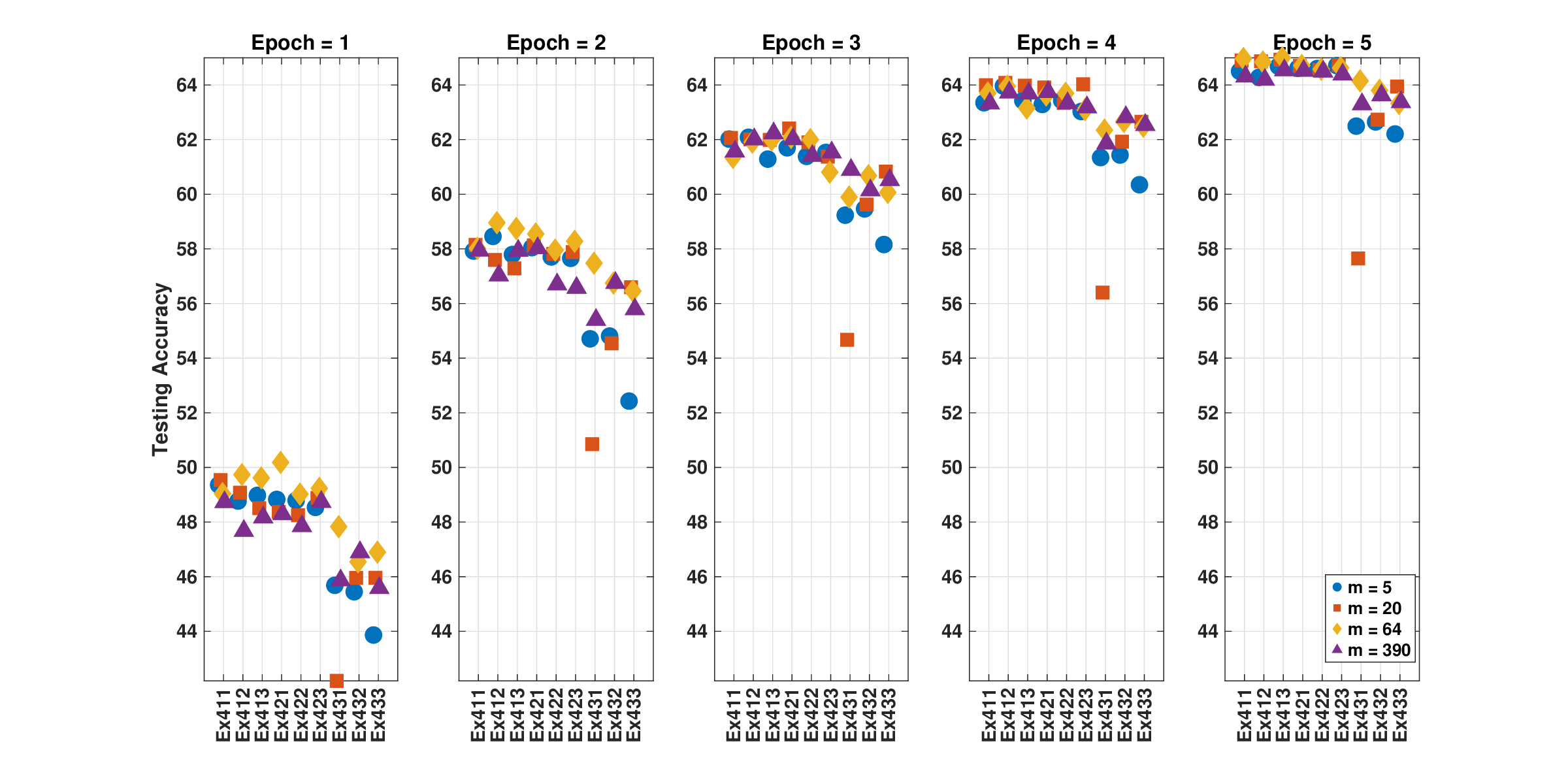}
    \end{adjustbox}    
    \caption{\footnotesize{Testing average accuracy of \trishtre with $(\alpha, \gamma_1, \gamma_2)=(1, \gamma_1, \gamma_2)$, $m\in\{5,20,64,N_b\}$. Top:  \texttt{MNSIT}, $N_b=468$. Bottom: \texttt{CIFAR10}, $N_b=390$.}}
    \label{fig:DL(m)}
\end{figure}

\appendix
\renewcommand{\appendixname}{}

\section{Some intermediate results}\label{AppA}
We present some theoretical results required in \Cref{AppB}.
\begin{lemma}\label{bound_diff_f}
Suppose Assumption~\ref{ass.f} holds, then
\begin{equation}\label{diff_f}
f(x_{k+1})-f(x_k)\le\underbrace{g_k^T p_k+\frac{1}{2\mu_k}\|p_k\|^2}_{m_k(p_k)}+(\nabla f(x_k)-g_k)^T p_k+\frac{L}{2}\|p_k\|^2.
\end{equation}
\end{lemma}


\begin{proof}
By \eqref{eq.g_Lip}, we have
$$
f(x_{k+1})-f(x_k) \leq \nabla f(x_k)^Tp_k + \frac{L}{2} \|p_k\|^2.
$$
Adding, subtracting $g_k^T p_k$, and using the fact that $\frac{1}{2\mu_k}\|p_k\|^2\ge 0$ conclude the proof.
\qedhere

\end{proof}


\begin{lemma}\label{bound_model}
Suppose Assumption~\ref{ass.f} holds, then
\begin{equation}\label{bound_mk}
m_k(p_k)= g_k^T p_k+\frac{1}{2\mu_k}\|p_k\|^2\le -\Delta_k\|g_k\|+\frac{1}{2\mu_k}\Delta_k^2.
\end{equation}
\end{lemma}
\begin{proof}
If $p_k=-\mu_kg_k$, then
$m_k(p_k)=-\thalf  \mu_k\|g_k\|^2 $.
If $p_k=-\frac{\Delta_k}{\|g_k\|}g_k$, then $m_k(p_k)=-\Delta_k\|g_k\|+\frac{1}{2\mu_k}\Delta_k^2$.
Since $m_k(-\mu g_k)$ is the minimum value attainable by $m_k$, the proof is concluded.
\qedhere
\end{proof}
\section{Proofs from \Cref{ConvergenceSec}}\label{AppB}
\subsection{\sc{Proof of Lemma \ref{lem.basic}}}\label{AppB1}
\begin{proof}
Given the step taken by the \Cref{alg.TRishBB}, i.e., $p_k = x_{k+1} - x_k$, the inequality \eqref{eq.g_Lip} implies
  \begin{equation}
    f(x_{k+1}) = f(x_k + p_k) \leq f(x_k) + \nabla f(x_k)^T p_k + \frac{L}{2} \|p_k\|^2.
\end{equation}
Thus, considering the events that Case 1--3 occur, i.e., $C_{i,k}$ for $i=1,2,3$, the law of total probability gives
  \begin{align}
    & \E_k[f(x_{k+1})] - f(x_k)
    \leq \E_k[\nabla f(x_k)^Tp_k] + \frac{L}{2} \E_k[\|p_k\|^2] \nonumber \\
      =  &  \sum_{i=1}^3 \Pbk[C_{i,k}]\E_k[\nabla f(x_k)^Tp_k | C_{i,k}] 
       + \frac{L}{2} \sum_{i=1}^3 \Pbk[C_{i,k}]\E_k[\|p_k\|^2 | C_{i,k}]. \label{eq.trick1}
\end{align}
Further, 
\begin{eqnarray}
\E_k[\nabla f(x_k)^Tp_k | C_{i,k}]&=&\Pbk[\nosbb|C_{i,k}]\E_k[\nabla f(x_k)^Tp_k | C_{i,k}\cap\nosbb] \nonumber\\
& + & \Pbk[\sbb|C_{i,k}]\E_k[\nabla f(x_k)^Tp_k | C_{i,k}\cap\sbb]  \label{expcki_1} ,
\end{eqnarray}
and
\begin{eqnarray}
\E_k[\|p_k\|^2 | C_{i,k}]&=&  \Pbk[\nosbb|C_{i,k}]\E_k[  \|p_k\|^2 | C_{i,k}\cap\nosbb]\nonumber \\
& + &  \Pbk[\sbb|C_{i,k}]\E_k[  \|p_k\|^2 | C_{i,k}\cap\sbb]  \label{expcki_2},
\end{eqnarray}
for $i=1,2, 3$.  Now, we analyze the events separately.
\smallskip

\noindent \textit{\sc{Event} $ {C_{1,k}}$.} Consider the event $\cunonosbb\eqdef C_{1,k}\cap\nosbb$. It holds $p_k = -\gamma_{1} \alpha g_k$, 
and \cite[Equations (8), (9)]{curtis2019stochastic} gives
\begin{align}
\E_k[\nabla f(x_k)^Tp_k | \cunonosbb] \nonumber 
& \le \ -\gamma_{1} \alpha \E_k\left[\nabla f(x_k)^Tg_k | \cunonosbb\right] \\
& + (\gamma_{1} - \gamma_{2}) 
\alpha \Pbk\left[\ev | \cunonosbb\right] \E_k\left[\nabla f(x_k)^Tg_k |\cunonosbb \cap \ev\right], \nonumber
\end{align}
and
\begin{equation*}
    \E_k\left[\|p_k\|^2 | \cunonosbb\right] = \gamma_{1}^2 \alpha^2 \E_k\left[\|g_k\|^2 | \cunonosbb\right].
\end{equation*}
In the event $\cunosbb\eqdef C_{1,k}\cap\sbb$, it holds $p_k = -\mu_k g_k$ with $\mu_k\le  \gamma_{1}\alpha$,
and paralleling \cite[Equations (8) and (9)]{curtis2019stochastic} results in
\begin{align}
\E_k[\nabla f(x_k)^Tp_k | \cunosbb] & = \E_k[-\mu_k\nabla f(x_k)^Tg_k |\cunosbb]\nonumber \\
& \le (\gamma_{1}\alpha - \mu_{min})  \Pbk[\ev | \cunosbb] \E_k[\nabla f(x_k)^Tg_k |\cunosbb \cap \ev] \nonumber\\
& -\gamma_{1} \alpha \E_k[\nabla f(x_k)^Tg_k | \cunosbb], \nonumber
\end{align}
and
\begin{equation*}
    \E_k[\|p_k\|^2 |\cunosbb] \le \gamma_{1}^2 \alpha^2 \E_k[\|g_k\|^2 | \cunosbb].
\end{equation*}
Now, we combine  \eqref{expcki_1}, \eqref{expcki_2} and the previous results. We take into account that 
$(\gamma_{1}-  \gamma_{2})\alpha\le 
(\gamma_{1}\alpha- \mu_{min})$ by Assumption \ref{hp_mu} and $$\Pbk(\ev|\cunonosbb)\le \Pbk(\ev|C_{1,k}),
\qquad \Pbk(\ev|\cunosbb)\le \Pbk(\ev|C_{1,k}).$$
Thus, we get
\begin{align*}
&\,  \E_k[\nabla f(x_k)^Tp_k | C_{1,k}]
\le\,    (\gamma_1\alpha-\mu_{min})\Pbk[\sbb|C_{1,k}]\Pbk(\ev|\cunosbb)\E_k[\nabla f(x_k)^Tg_k | \cunosbb\cap \ev]
\Bigr)\\
\, -&\gamma_{1}\alpha \Bigl(\Pbk[\nosbb|C_{1,k}]\E_k[\nabla f(x_k)^Tg_k | \cunonosbb] +
\Pbk[\sbb|C_{1,k}]\E_k[\nabla f(x_k)^Tg_k|\cunosbb]\Bigr)\\
\, +&
(\gamma_{1}\alpha - \mu_{min}) \Bigl(\Pbk[\nosbb|C_{1,k}]\Pbk[\ev|\cunonosbb]\E_k[\nabla f(x_k)^Tg_k | \cunonosbb\cap \ev],
\end{align*}
i.e.,
\begin{align*}
\E_k[\nabla f(x_k)^Tp_k | C_{1,k}] &\le (\gamma_{1}\alpha - \mu_{min})   \Pbk(\ev|C_{1,k})\E_k[\nabla f(x_k)^Tg_k | C_{1,k}\cap \ev] \\
& -\gamma_{1}\alpha  \E_k[\nabla f(x_k)^Tg_k | C_{1,k} ],
\end{align*} 
and
\begin{equation}\label{case1_final2}
    \E_k[\|p_k\|^2 |C_{1,k}] \le \gamma_{1}^2 \alpha^2 \E_k[\|g_k\|^2 | C_{1,k}].
\end{equation}


\noindent \textit{\sc{Event} $ {C_{2,k}}$.} The event $C_{2,k}$ occurs when $\|g_k\|^{-1} \leq \gamma_{1}$ and $\|g_k\|^{-1} \geq \gamma_{2}$.
In the event $\cduenosbb\eqdef C_{2,k}\cap\nosbb$, it holds
$p_k=-\alpha g_k/\|g_k\|$  and  \cite[Equations (10) and (11)]{curtis2019stochastic} gives
  \begin{align}
\E_k[\nabla f(x_k)^Tp_k | \cduenosbb] &\le-\gamma_{1}\alpha\E_k\left[\nabla f(x_k)^Tg_k |  \cduenosbb\right] \nonumber \\
& + (\gamma_{1} - \gamma_{2}) \alpha \Pbk
\left[\ev |  \cduenosbb\right] \E_k\left[\nabla f(x_k)^Tg_k |  \cduenosbb \cap \ev\right]\nonumber 
\end{align}
and 
\begin{equation*} 
\Ek\left[\|p_k\|^2 |  \cduenosbb\right] = \alpha^2 \leq \gamma_{1}^2 \alpha^2 \Ek\left[\|g_k\|^2 |  \cduenosbb\right].
\end{equation*}
Assumption \ref{hp_mu} ensures that $\mu_{min}\le \frac{\alpha}{\|g_k\|}$ and in the event $ \cduesbb\eqdef C_{2,k}\cap \sbb$, it holds $\mu_{min}\le \mu_k\le \frac{\alpha}{\|g_k\|}$.
Then, 
\begin{align}
\E_k\left[\nabla f(x_k)^Tp_k | \cduesbb\right ] 
&\le -\mu_{min}\Pbk\left[\ev | \cduesbb\right] \E_k\left[ \nabla f(x_k)^Tg_k  \bigg| \cduesbb \cap \ev\right] \nonumber\\
&- \alpha 
 \Pbk\left[\overline\ev | \cduesbb\right] \E_k\left[\frac{\nabla f(x_k)^Tg_k}{\|g_k\|} \bigg|\cduesbb \cap \overline\ev\right] \nonumber\\
&\leq -\mu_{min}  \Pbk\left[\ev| \cduesbb\right] \E_k\left[\nabla f(x_k)^Tg_k | \cduesbb \cap \ev\right] \nonumber\\
&- \gamma_{1} \alpha 
\Pbk\left[\overline\ev | \cduesbb\right ] \E_k\left [\nabla f(x_k)^Tg_k | \cduesbb \cap \overline\ev\right ] \nonumber\\
&= -\gamma_{1}\alpha\E_k[\nabla f(x_k)^Tg_k |  \cduesbb] \nonumber\\
&+ (\gamma_{1}\alpha - \mu_{min})   \Pbk\left [\ev | \cduesbb\right] \E_k\left[\nabla f(x_k)^Tg_k |  \cduesbb \cap \ev\right], \nonumber 
\end{align}
and 
\begin{equation*}
\Ek[\|p_k\|^2 | \cduesbb] \le \alpha^2 \leq \gamma_{1}^2 \alpha^2 \Ek[\|g_k\|^2 | \cduesbb].
\end{equation*}
Now, we combine  \eqref{expcki_1}, \eqref{expcki_2} and the previous results. We take into account that 
$(\gamma_{1}-  \gamma_{2})\alpha\le 
(\gamma_{1}\alpha- \mu_{min})$ by Assumption \ref{hp_mu} and that 
$$\Pbk\left(\ev| \cduenosbb\right)\le \Pbk\left(\ev|C_{2,k}\right),\qquad \Pbk\left(\ev| \cduesbb\right)\le \Pbk\left(\ev|C_{2,k}\right).$$ 
It follows
\begin{align*}
& \E_k[\nabla f(x_k)^Tp_k | C_{2,k}]
\le  -\gamma_{1}\alpha \Bigl(\Pbk[\nosbb|C_{2,k}]\E_k\left[\nabla f(x_k)^Tp_k |  \cduenosbb\right]\\
 & +
\Pbk\left[\sbb|C_{2,k}\right]\E_k\left[\nabla f(x_k)^Tp_k |  \cduesbb\right]\Bigr) \\
& + (\gamma_{1}\alpha - \mu_{min}) \Bigl( 
\Pbk\left[\nosbb|C_{2,k}\right]\Pbk\left[\ev|C_{2,k}\cap \nosbb\right]
\E_k\left[\nabla f(x_k)^Tp_k |  \cduenosbb\cap \ev\right] \\
& + \Pbk[\sbb|C_{2,k}]\Pbk\left[\ev| \cduesbb\right]\E_k\left[\nabla f(x_k)^Tp_k |  \cduesbb\cap \ev\right]\Bigr), 
\end{align*}
and finally
\begin{align}
\E_k[\nabla f(x_k)^Tp_k | C_{2,k}] \le & (\gamma_{1}\alpha - \mu_{min})  \Pbk(\ev|C_{2,k})\E_k[\nabla f(x_k)^Tp_k | C_{2,k}\cap \ev] \nonumber\\
 - & \gamma_{1}\alpha  \E_k[\nabla f(x_k)^Tp_k | C_{2,k} ].
\end{align} \label{case2_final1}
Further,
\begin{equation}\label{case2_final2}
\Ek[\|p_k\|^2 | C_{2,k}] \le \alpha^2 \leq \gamma_{1}^2 \alpha^2 \Ek[\|g_k\|^2 | C_{2,k}].
\end{equation}

\noindent \textit{\sc{Event} $ {C_{3,k}}$.}
Following \cite[Equations (12) and(13)]{curtis2019stochastic} and the lines of the proof for the event $C_{1,k}$, we obtain 
\begin{align}
  \E_k[\nabla f(x_k)^Tp_k | C_{3,k}]\le & (\gamma_{1}\alpha- \mu_{min}) \Pbk(\ev|C_{3,k})\E_k[\nabla f(x_k)^Tg_k | C_{3,k}\cap \ev] \nonumber  \\ 
  - & \gamma_{1}\alpha  \E_k[\nabla f(x_k)^Tg_k | C_{3,k} ], \label{case3_final1}
\end{align}
and 
\begin{equation}\label{case3_final2}
    \E_k[\|p_k\|^2 | C_{3,k} ] \le \gamma_{1}^2 \alpha^2 \E_k[\|g_k\|^2 | C_{3,k}].
\end{equation}
 Combining \eqref{eq.trick1}, \eqref{case3_final1} and \eqref{case3_final2}, the claim follows.
\qedhere
\end{proof}
\subsection{\sc{Proof of Lemma \ref{diff_f_g_unbiased}}}\label{AppB2}
\begin{proof} Regarding  (\ref{decr_unbiased}), we analyze the different events  ${C_{1,k}}$, $i=1,2,3$, separately.
\smallskip

\noindent \textit{\sc{Event} ${C_{1,k}}$.} Since $\Delta_k=\gamma_{1}\alpha\|g_k\|$, the rightmost inequality in \eqref{alpha_mu_final}
and \eqref{bound_mk} give
$$
m_k(p_k)\le -\gamma_{1}\alpha\left(1-\thalf \frac{\gamma_{1}\alpha}{\mu_k}\right )\|g_k\|^2 
\le -  \frac 3  8 \gamma_{1}\alpha \|g_k\|^2.
$$
Using $\|p_k\| \leq \gamma_{1}\alpha\|g_k\|$, Lemma~\ref{bound_diff_f}, the leftmost inequality in \eqref{alpha_mu_final} and  $\gamma_{1} \geq \gamma_{2}$, we find that
\begin{align}
&\ f(x_{k+1}) - f(x_k) \nonumber \\
&\leq -  \frac 3  8 \gamma_{1}\alpha \|g_k\|^2 + \gamma_{1}\alpha \|\nabla f(x_k) - g_k\|\|g_k\| + \thalf L
\gamma_{1}^2 \alpha^2 \|g_k\|^2 \nonumber \\
&\leq -  \frac 3  8  \gamma_{1}\alpha \|g_k\|^2 +  \gamma_{1}\alpha \|\nabla f(x_k) - g_k\|\|g_k\| + \frac{1}{16}
  \gamma_{1}  \alpha \|g_k\|^2. \nonumber 
\end{align}
\noindent
Then, by $\gamma_{1} \geq \gamma_{2}$ and
$$0 \leq \left(\thalf \|g_k\| - \|\nabla f(x_k) - g_k\|\right)^2 = \tfrac14\|g_k\|^2 - \|\nabla f(x_k) - g_k\|\|g_k\| +\|\nabla f(x_k) - g_k\|^2,$$  
we obtain
\begin{equation}
\begin{aligned}
& \ f(x_{k+1}) - f(x_k) \\ 
&\leq   -\frac{1}{16} \gamma_{1}\alpha \|g_k\|^2 + \gamma_{1}\alpha
\|\nabla f(x_k) - g_k\|^2 \\
&\leq  -\frac{1}{16}\gamma_{2} \alpha   \|g_k\|^2 + \frac{\gamma_{1}^2}{\gamma_{2}} \alpha \|\nabla f(x_k) - g_k\|^2.
\end{aligned}
\end{equation}

\noindent \textit{\sc {Event}  ${C_{2,k}}$.} Since $\Delta_k= \alpha$, using \eqref{bound_mk}, 
the rightmost inequality in \eqref{alpha_mu_final} and
$1/\gamma_{1}\le \|g_k\|\le 1/\gamma_{2}$ gives
$$
m_k(p_k)\le -\alpha\|g_k\|+\thalf \frac{\alpha^2}{\mu_k} \le 
-\alpha\|g_k\|+ \frac 5 8 \frac{\alpha}{\gamma_1} \le  -\frac 3 8 \alpha\|g_k\| 
\le  -\frac 3 8 \gamma_{2}\alpha \|g_k\|^2.
$$
Using $1/\gamma_{1}\le \|g_k\|\le 1/\gamma_{2}$, Lemma~\ref{bound_diff_f}, and the leftmost inequality in \eqref{alpha_mu_final}, we find that
\begin{equation*}
\begin{aligned}
\ f(x_{k+1}) - f(x_k) \leq&\  -\frac 3 8 \gamma_{2}\alpha \|g_k\|^2 + \|\nabla f(x_k) - g_k\| \alpha + 
\thalf L\alpha^2 \\
\leq&\  -\frac 3 8 \gamma_{2}\alpha \|g_k\|^2 + \ \|\nabla f(x_k) - g_k\|\alpha + \thalf \gamma_{1}^2 \alpha^2 L\|g_k\|^2\\
\leq&\  -\frac 3 8 \gamma_{2}\alpha \|g_k\|^2 + \ \|\nabla f(x_k) - g_k\|\alpha + \frac{1}{16} \gamma_{2} \alpha\|g_k\|^2.
\end{aligned}
\end{equation*}
Since $$0 \leq \frac{\gamma_{2}}{\gamma_{1}^2}\left(\thalf - \frac{\gamma_{1}^2}{\gamma_{2}} \|\nabla f(x_k) - g_k\|\right)^2 = \frac{\gamma_{2}}{4\gamma_{1}^2} - \|\nabla f(x_k) - g_k\| + \frac{\gamma_{1}^2}{\gamma_{2}} \|\nabla f(x_k) - g_k\|^2,$$
we obtain
\begin{equation*}
\begin{aligned}
 f(x_{k+1}) - f(x_k) 
        \leq \  -\frac{1}{ 16}\gamma_{2} \alpha \|g_k\|^2 +   
				\frac{\gamma_{1}^2}{\gamma_{2}} \alpha\|\nabla f(x_k) - g_k\|^2. 
\end{aligned}
\end{equation*}

\noindent \textit{\sc{Event} ${C_{3,k}}$.} The proof trivially follows what we argued for the event $ {C_{1,k}}$ as we have  $\gamma_{1} \geq \gamma_{2}$. 

 Summarizing the analysis above, inequality \eqref{decr_unbiased} holds.
Finally, regarding (\ref{eq.for_later}), due to \eqref{exp_err} and \eqref{eq.unbiased_consequence}, 
inequality \eqref{decr_unbiased} yields
  \begin{align}
        &\ \E_k [f(x_{k+1})] - f(x_k) \nonumber \\
\leq&\ -  \frac{1}{16} \gamma_2 \alpha (\|\nabla f(x_k)\|^2 + \E_k[\|\nabla f(x_k) - g_k\|^2]) + \frac{\gamma_1^2}{\gamma_2} \alpha \E_k [\|\nabla f(x_k) - g_k\|^2] \nonumber \\
       =&\ -  \frac{1}{16}  \gamma_2 \alpha \|\nabla f(x_k)\|^2  + \alpha\left(\frac{\gamma_1^2}{\gamma_2} -   \frac{1}{16} \gamma_2\right)\E_k[\|\nabla f(x_k) - g_k\|^2] \nonumber
\end{align}
which give the thesis.
\qedhere
\end{proof}

\section{Classification Problems}\label{AppC}
\begin{table}[t!]
\centering
\scriptsize
    \caption{Classification test problems.}
    \label{tab:testProblems}
    \begin{tabular}{|p{1.6cm}|c|c|c|c|p{2.4cm}|p{1.4cm}|}
        \hline
        \textbf{Dataset}&  \textbf{$\mathbf{(N, {N}_t)}$} & \textbf{$\mathbf{d}$} & \textbf{$\mathbf{n}$} & \textbf{Type} & \textbf{$\mathbf{f_i(x)}$} & \textbf{NN} \\
        \hline
        \texttt{a1a} \cite{chang2011libsvm} & $(1605,29351)$ & $123$ & $123$ & Binary &  Logistic regression &  - \\
        \hline
        \texttt{w1a} \cite{chang2011libsvm} & $(2477,47272)$ & $300$ & $300$ & Binary & Logistic regression &  - \\
        \hline
        \texttt{cina} \cite{chang2011libsvm} & $(10000,6033)$ & $132$ & $132$ & Binary & Logistic regression &  - \\
        \hline
        \texttt{MNIST} \cite{deng2012mnist} & $(60000, 10000)$ & $28\times 28\times 1$ & $431030$ & Multi & Cross-entropy & \texttt{NET-1} \\
        \texttt{CIFAR10} \cite{krizhevsky2009learning} & $(50000, 10000)$ & $32\times 32\times 3$ & $ 545098$ & Multi & Cross-entropy & \texttt{NET-2} \\
        \hline
    \end{tabular}
\end{table}

\Cref{tab:testProblems} summarizes the five datasets used in our experiments for both binary and multi-class classification tasks in Section~\ref{S.classProb}. Each test dataset with $C$ possible classes is divided into training and testing sets, denoted as $\{(a_i,b_i)\,| i\in \mathcal{N}\}$ and $\{(a_i,b_i)\, | i\in \mathcal{N}_t\}$, with $|\mathcal{N}|=N$ and $|\mathcal{ N}_t|={N}_t$ samples, respectively; the values of $N$ and ${N}_t$ are reported in the table.  The vector  $a_i \in \mathbb{R}^d$ represents the $i$-th data with $d$ features while $b_i$ is the corresponding categorical label associated. 

The dataset \texttt{cina} was normalized such that each feature vector has components in the interval $[0,1]$. For \texttt{MNIST} and \texttt{CIFAR10}, besides the zero-one rescaling, we implemented z-score normalization to have zero mean and unit variance. Precisely, all $N$ training images, (3-D arrays of size $\sigma_1 \times \sigma_2\times \sigma_2 $ where $\sigma_1$, $\sigma_2$ and $\sigma_3$ denote the height, width, and number of channels of the images, respectively), undergo normalization by subtracting the mean (a 3-D array) and dividing by the standard deviation  (a 3-D array)  component-wise. Data  in the testing set are also normalized using the same statistics as in the training data.
 The architecture of the neural network (NN) is described in \Cref{tab:nets}.\footnote{The compound $Conv(5,20,1,0){/}ReLU{/}MaxPool(2,2,0)$ indicates a convolutional layer ($Conv$) using $20$ filters of size $5\times 5$, stride $1$, and no padding, a nonlinear activation function ($ReLu$) and, a 2-D max-pooling layer with a channel of size $2\times2$, stride 2 and no padding. The syntax $FC(C){/}Softmax$ denotes a fully connected layer of $C$ neurons followed by the $softmax$ layer.} For \texttt{MNIST} dataset with grayscaled images in $C=10$ classes (digits 0--9), and \texttt{CIFAR10} with RGB images in $C=10$ categories (e.g., dogs and ships), we used \texttt{NET-1} and \texttt{NET-2} respectively. 

\begin{table}[t!]
\centering
\scriptsize
\caption{NN architecture.}
\label{tab:nets}
\begin{tabular}{|c|p{9.7cm}|}
    \hline
    \texttt{NET-1} & $Input(d)$\newline
    $Conv(5,20,1,0){/}ReLU{/}MaxPool(2\times2,2,0)$\newline
    $Conv(5,50,1,0){/}ReLU{/}MaxPool(2\times2,2,0)$\newline
    $FC(500){/}ReLU$\newline
    $FC(C){/}Softmax$ 
    \\ 
    \hline
    \texttt{NET-2} & $Input(d)$\newline
    $Conv(3,32,1,same){/}ReLU{/}AveragePool(2\times2,2,0)$\newline
    $Conv(3,64,1,same){/}ReLU{/}AveragePool(2\times2,2,0)$\newline
    $FC(128){/}ReLU$\newline
    $FC(C){/}Softmax$
    \\ 
    \hline
    \end{tabular}
\end{table}

\section*{Acknowledgments} 
The authors are members of the INdAM-GNCS Research Group. This research was partially supported by INDAM-GNCS through Progetti di Ricerca 2023, by PNRR - Missione 4 Istruzione e Ricerca - Componente C2 Inve\-stimento 1.1, Fondo per il Programma Nazionale di Ricerca e Progetti di Rilevante Interesse Nazionale (PRIN) funded by the European Commission under the NextGeneration EU programme, project ``Advanced optimization METhods for automated central veIn Sign detection in multiple sclerosis from magneTic resonAnce imaging (AMETISTA)'',  code: P2022J9SNP, MUR D.D. financing decree n. 1379 of 1st September 2023 (CUP E53D23017980001), project ``Numerical Optimization with Adaptive Accuracy and Applications to Machine Learning'', code: 2022N3ZNAXMUR D.D. financing decree n. 973 of 30th June 2023 (CUP B53D23012670006), and by Partenariato esteso FAIR ``Future Artificial Intelligence Research'' SPOKE 1 Human-Centered AI. Obiettivo 4, Project ``Mathematical and Physical approaches to innovative Machine Learning technologies (MaPLe)''. 

\section*{Competing interests declarations}       
The authors have no relevant interests to disclose.

\section*{Data availability statements}
The datasets utilized in this research are publicly accessible and commonly employed benchmarks in the field of Machine Learning and Deep Learning, see \cite{chang2011libsvm, deng2012mnist}.





\end{document}